%% file: approximatelastversion.tex
\documentclass[11pt]{article}
\usepackage[latin1]{inputenc}
\usepackage[english]{babel}
\usepackage{hyperref}
\usepackage{amsmath,amsfonts,amscd,amssymb}
\usepackage[all]{xy}
\usepackage{graphicx,epsfig}
\usepackage{makeidx}
\usepackage{color}
\usepackage[active]{srcltx} 
\newtheorem{thm}{Theorem}[subsection]
\newtheorem{cor}[thm]{Corollary}
\newtheorem{lem}[thm]{Lemma}
\newtheorem{prop}[thm]{Proposition}
\newtheorem{deft}[thm]{Definition}
\newtheorem{rek}[thm]{Remark}
\newtheorem{conj}[thm]{Conjecture}
\let\noi=\noindent

\let\sur=\overline

\newcommand{\dsp}{\displaystyle}
\def\Nu{\mathcal{V}} 
\def\N{\mathbb{N}} 
\def\Z{\mathbb{Z}} 
\def\C{\mathbb{C}} 
\def\Q{\mathbb{Q}} 
\def\F{\mathbb{F}}
\def\sper{\mbox{Sper}}

\def\gr{\mbox{gr}}
\def\init{\mbox{in}}
\def\notin{\mbox{$\in$ \hspace{-.8em}/}} 
\def\notsubset{\mbox{$\subset$ \hspace{-.9em}/}} 
\def\sgn{\mbox{sgn}} 
\makeindex
\title{APPROXIMATE ROOTS OF A VALUATION AND THE PIERCE--BIRKHOFF CONJECTURE}
\author{F. Lucas\\
D{\'e}partement de Math{\'e}matiques/CNRS UMR 6093\\
Universit{\'e} d'Angers\\
2, bd Lavoisier\\
49045 Angers C{\'e}dex, France \and
J. Madden\\
Department of Mathematics\\
Louisiana State University at Baton Rouge,\\
Baton Rouge, LA, USA \and
D. Schaub\\
D{\'e}partement de Math{\'e}matiques/CNRS UMR 6093\\
Universit{\'e} d'Angers\\
2, bd Lavoisier\\
49045 Angers C{\'e}dex, France \and
M. Spivakovsky\\
Inst. de Math\'ematiques de Toulouse/CNRS 5219\\
Universit\'e Paul Sabatier\\
118, route de Narbonne\\
31062 Toulouse cedex 9, France.}

\date{}

\begin{document}
\maketitle

\abstract{In this paper, we construct an object, called a system of
  approximate roots of a valuation, centered in a regular local ring, which
describes the fine structure of the valuation (namely, its valuation ideals and
the graded algebra). We apply this construction to valuations associated to a
point of the real
spectrum of a regular local ring $A$. We give two versions of the construction:
the first, much simpler, in a special
case (roughly speaking, that of rank 1 valuations), the second --- in the case of complete regular local rings and valuations of arbitrary rank.

We then describe certain subsets $C\subset\sper\ A$ by explicit
formulae in terms of approximate roots; we conjecture that these sets
satisfy the Connectedness (respectively, Definable Connectedness)
conjecture. Establishing this for a certain regular ring $A$ would imply that $A$ is a Pierce--Birkhoff ring (this means that the Pierce--Birkhoff conjecture holds in $A$).

Finally, we use these constructions and results to prove the
Definable Connectedness conjecture (and hence \textit{a fortiori} the
Pierce--Birkhoff conjecture) in the special case when $\dim\ A=2$.}

\section*{Introduction}
\label{Intro}
\setcounter{section}{0}
\setcounter{subsection}{1}
All the rings in this paper will be commutative with 1. Let $R$ be a
real closed field. Let $B=R[x_1,\dots,x_n]$. If $A$ is a ring and
$\mathfrak p$ a prime ideal of $A$,
$\kappa(\mathfrak p)$ will denote the residue field of $\mathfrak p$.

The Pierce--Birkhoff conjecture asserts that any
piecewise-polynomial function
$$
f:R^n\rightarrow R
$$
can be expressed as a maximum of minima of a finite family of polynomials in $n$
variables. We start by giving the precise statement of the conjecture as it was
first stated by M. Henriksen and J. Isbell in the early nineteen
sixties.
\begin{deft}\label{pw} A function $f:R^n\to R$ is said to be \textbf{piecewise
  polynomial} if $R^n$ can be covered by a finite collection of closed
  semi-algebraic sets $P_i$ such that for each $i$ there exists a
  polynomial $f_i\in B$ satisfying
  $\left.f\right|_{P_i}=\left.f_i\right|_{P_i}$.
\end{deft}
Clearly, any piecewise polynomial function is continuous. Piecewise
polynomial functions form a ring, containing $B$, which is denoted by
$PW(B)$.\medskip

On the other hand, one can consider the (lattice-ordered) ring of all
the functions obtained from $B$ by iterating the operations of $\sup$
and $\inf$. Since applying the operations of sup and inf to
polynomials produces functions which are piecewise polynomial, this
ring is contained in $PW(B)$ (the latter ring is closed under $\sup$
and $\inf$). It is natural to ask whether the two rings coincide. The
precise statement of the conjecture is:
\begin{conj}\textnormal{\textbf{(Pierce-Birkhoff)}}\label{PB} If
  $f:R^n\to R$ is in $PW(B)$, then there exists a finite family of
  polynomials $g_{ij}\in B$ such that
  $f=\sup\limits_i\inf\limits_j\{g_{ij}\}$ (in other words, for all
  $x\in R^n$, $f(x)=\sup\limits_i\inf\limits_j\{g_{ij}(x)\}$).
\end{conj}
This paper represents the second step of our program for proving the
Pierce--Birkhoff conjecture in its full generality. The
starting point of this program is the abstract formulation of the
conjecture in terms of the real spectrum of $B$ and separating ideals
proposed by J. Madden in 1989 \cite{Mad1}, which we
now recall, together with the relevant definitions. For a general
introduction to real spectrum, we refer the reader to
\cite{BCR}, Chapter 7, \cite{And}, Chapter II or \cite{PD}, 4.1, page
81 and thereafter; see also ``Bibliographical and historical
comments'' on p. 109 at the end of that chapter.

Let $A$ be a ring. A point $\alpha$ in the real spectrum of $A$ is, by
definition, the data of a prime ideal $\mathfrak{p}$ of $A$, and a
total ordering $\leq$ of the quotient ring $A/\mathfrak{p}$, or,
equivalently, of the field of fractions of $A/\mathfrak{p}$.  Another
way of defining the point $\alpha$ is as a homomorphism from $A$ to a
real closed field, where two homomorphisms are identified if they have
the same kernel $\mathfrak{p}$ and induce the same total ordering on
$A/\mathfrak{p}$.

The ideal $\mathfrak{p}$ is called the support of $\alpha$ and denoted
by $\mathfrak{p}_\alpha$, the quotient ring $A/\mathfrak{p}_{\alpha}$
by $A[\alpha]$, its field of fractions by $A(\alpha)$ and the real
closure of $A(\alpha)$ by $k(\alpha)$. The total ordering of
$A(\alpha)$ is denoted by $\leq_\alpha$. Sometimes we write
$\alpha=(\mathfrak{p}_\alpha,\leq_\alpha)$.

\begin{deft} The real spectrum of $A$, denoted by Sper $A$, is the
  collection of all pairs $\alpha=(\mathfrak{p}_\alpha,\leq_\alpha)$,
  where $\mathfrak{p}_\alpha$ is a prime ideal of $A$ and
  $\leq_\alpha$ is a total ordering of $A/\mathfrak{p}_{\alpha}$.
\end{deft}

We use the following notation: for an element $f\in A$,
$f(\alpha)$ stands for the natural image of $f$ in $A[\alpha]$ and the
inequality $f(\alpha)>0$ really means $f(\alpha)>_\alpha0$.

\noi The real spectrum $\sper\ A$ is endowed with two natural topologies.
The first one, called the {\bf spectral (or Harrison) topology}, has
basic open sets of the form
$$
U(f_1,\ldots,f_k)=\left\{\alpha\ \left|\ f_1(\alpha)>0,\ldots,f_k(\alpha)>0\right.\right\}
$$
with $f_1,...,f_k\in A$.

The second is the \textbf{constructible topology} whose basic open sets
are of the form
$$
V(f_1,\ldots,f_k,g)=\left\{\alpha\ \left|\
f_1(\alpha)>0,\ldots,f_k(\alpha)>0,g(\alpha)=0\right.\right\},
$$
where $f_1,...,f_n,g\in A$. Boolean combinations of sets of the form
$V(f_1,\ldots,f_n,g)$ are
called \textbf{constructible sets} of $\sper\ A$.

For more information about the real spectrum, see \cite{BCR}; there is
also a brief introduction to the real spectrum and its relevance to
the Pierce--Birkhoff conjecture in the Introduction to \cite{LMSS}.

\begin{deft}
Let
$$
f :\sper\ A \to \coprod_{\alpha \in \sper\ A} A(\alpha)
$$
be a map such that, for each $\alpha \in \sper\ A$, $f(\alpha) \in A(\alpha)$. We
say that $f$ is \textbf{piecewise polynomial}
(denoted by $f\in PW(A)$) if there exists a
covering of $\sper\ A$ by a finite family $(S_i)_{i \in I}$ of
constructible sets, closed in the spectral topology, and a family
$(f_i)_{i \in I}$, $f_i \in A$ such that, for each $\alpha \in S_i$,
$f(\alpha)=f_i(\alpha)$.

We call $f_i$ a \textbf{local representative} of $f$ at $\alpha$ and denote it
by $f_\alpha$ ($f_\alpha$ is, in general, not uniquely determined by
$f$ and $\alpha$; this notation means that one such local
representative has been chosen once and for all).
\end{deft}

Note that $PW(A)$ is naturally a lattice ring: it is equipped
with the operations of maximum and minimum. Each element of $A$
defines a piecewise polynomial function. In this way we get a
natural injection $A \subset PW(A)$.

\begin{deft}
A ring $A$ is a Pierce-Birkhoff ring if, for each $f \in PW(A)$, there
exist a finite collection of $f_{ij} \in A$ such that
$f=\sup\limits_i\inf\limits_j f_{ij}$.
\end{deft}

In \cite{Mad1} Madden reduced the Pierce--Birkhoff
conjecture to a purely local statement about separating ideals and the
real spectrum. Namely, he gave the following definition:

\begin{deft} Let $A$ be a ring. For $\alpha,\beta\in\mbox{Sper}\ A$,
  the \textbf{separating ideal} of $\alpha$ and $\beta$, denoted by
  $<\alpha,\beta>$, is the ideal of $A$ generated by all the elements
  $f\in A$ which change sign between $\alpha$ and $\beta$, that is,
  all the $f$ such that $f(\alpha)\geq0$ and $f(\beta)\leq0$.
\end{deft}

\begin{deft} \label{def:localPB}
A ring $A$ is \textbf{locally Pierce-Birkhoff} at $\alpha, \beta$ if the
following condition holds. Let $f$ be a piecewise  polynomial
  function, let $f_\alpha \in A$ be a local representative of $f$ at
  $\alpha$ and $f_\beta \in A$ a local representative of $f$ at
  $\beta$. Then $f_\alpha-f_\beta \in <\alpha,\beta>$.
\end{deft}

\begin{thm} (Madden)
A ring $A$ is Pierce-Birkhoff if and only if it is locally
  Pierce-Birkhoff for all $\alpha, \beta \in \sper(A)$.
\end{thm}

 Let $\alpha,\beta$ be points in $\sper\ A$.

\begin{conj}\textnormal{\textbf{(local Pierce-Birkhoff conjecture at
      $\alpha$, $\beta$)}}\label{PBS} Let $A$ be a regular ring
  and $f$ a piecewise  polynomial function. Let $f_\alpha\in A$ be a
  local representative of $f$ at
  $\alpha$ and $f_\beta\in A$ a local representative of $f$ at
  $\beta$. Then $f_\alpha-f_\beta\in<\alpha,\beta>$.
\end{conj}

There are known counterexamples in the case $A$ is not regular
(eg. $A=R[x,y]/(y^2-x^3)$) and even with $A$ normal.
\medskip

\begin{rek} \label{rek:spezial}
Assume that $\beta$ is a specialization of $\alpha$. Then

(1) $<\alpha,\beta> = \mathfrak{p}_\beta$.

(2) $f_\alpha -f_\beta \in \mathfrak{p}_\beta$. Indeed, we may assume
that $f_\alpha \neq f_\beta$, otherwise there is nothing to
prove. Since $\beta \in \sur{\{\alpha\}}$, $f_\alpha$ is also a local
representative of $f$ at $\beta$. Hence $f_\alpha(\beta)-f_\beta(\beta)=0$, so $f_\alpha-f_\beta \in \mathfrak{p}_\beta$.

Therefore, to prove that a ring $A$ is Pierce-Birkhoff, it is
sufficient to verify Definition \ref{def:localPB} for all
$\alpha,\beta$ such that neither of $\alpha,\beta$ is a specialization
of the other.
\end{rek}

In \cite{LMSS}, we introduced
\begin{conj}\textnormal{\textbf{(the Connectedness
conjecture)}}\label{conn} Let $A$ be a regular ring. Let
$$
\alpha,\beta\in\sper\ A
$$
and let $g_1,\dots,g_s$ be a finite collection of elements of
  $A\setminus<\alpha,\beta>$. Then there exists a connected set
  $C\subset\sper\ A$ such that $\alpha,\beta\in C$ and
  $C\cap\{g_i=0\}=\emptyset$ for $i\in\{1,\dots,s\}$ (in other words,
  $\alpha$ and $\beta$ belong to the same connected component of the
  set $\sper\ A\setminus \{g_1\dots g_s=0\}$).
\end{conj}

\begin{deft}\label{cd} A subset $C$ of $Sper(A)$ is said to be
  \textbf{definably connected} if it is not a union of two non-empty disjoint
  constructible subsets, relatively closed for the spectral topology.
\end{deft}

\begin{conj}\textnormal{\textbf{(Definable connectedness
      conjecture)}}\label{DCC} Let $A$ be a regular ring.
Let $\alpha,\beta\in\sper\ A $ and let
  $g_1,\dots,g_s$ be a finite collection of elements of $A$, not
  belonging to $<\alpha,\beta>$. Then there exists a
  definably connected set $C\subset\sper\ A$ such that
  $\alpha,\beta\in C$ and $C\cap\{g_i=0\}=\emptyset$ for
  $i\in\{1,\dots,s\}$ (in other words, $\alpha$ and $\beta$ belong to
  the same definably connected component of the set $\sper\ A\setminus
  \{g_1\dots g_s=0\}$).
\end{conj}

In the earlier paper \cite{LMSS} we stated the Connectedness
conjecture (in the special case $A=B$) and proved that it implies the
Pierce--Birkhoff conjecture. Exactly the same proof applies verbatim to show that the
Definable Connectedness conjecture implies the Pierce-Birkhoff
conjecture for any ring $A$.
\medskip

One advantage of the Connectedness conjecture is that it is a
statement about $A$ (respectively, about the polynomial ring if $A=B$) which makes no
mention of piecewise polynomial functions.

Our problem is therefore one of constructing connected subsets of
$\sper\ A$ having certain properties.
\medskip

\noi\textbf{Terminology}: If $A$ is an integral domain, the phrase
``valuation of $A$'' will mean ``a valuation of the field of fractions
of $A$, non-negative on $A$''. Also, we will sometimes commit the
following abuse of notation. Given a ring $A$, a prime ideal
$\mathfrak p\subset A$, a valuation $\nu$ of $\frac A{\mathfrak p}$
and an element $x\in A$, we will write $\nu(x)$ instead of
$\nu(x\mod\mathfrak p)$, with the usual convention that
$\nu(0)=\infty$, which is taken to be greater than any element of the
value group.
\medskip

Given any ordered domain $D$, let $\bar D$ denote the convex hull of
$D$ in its field of fractions $D_{(0)}$:
$$
\bar D:=\left\{\left. f\in D_{(0)}\ \right|\ d>|f|\text{ for some
  }d\in D\right\}.
$$
The ring $\bar D$ is a valuation ring, since for any element $f\in
D_{(0)}$, either $f\in\bar D$ or $f^{-1}\in\bar D$.
For a point $\alpha\in\sper\ A$, we define $R_\alpha:=\overline{A[\alpha]}$. In
this way, to every point $\alpha\in\sper\ A$ we can canonically associate a
valuation $\nu_\alpha$ of $A(\alpha)$, determined by the valuation ring
$R_\alpha$. The maximal ideal
of $R_\alpha$ is $M_\alpha=\left\{x\ \in A(\alpha) \left|\
    |x|<\frac{1}{|z|}, \ \forall z\in A[\alpha] \setminus \{0\}
  \right.\right\}$; its residue field $k_\alpha$ comes equipped
with a total ordering, induced by $\le_\alpha$.
\medskip

\noi Let $U(R_\alpha)$ denote the multiplicative group of
units of $R_\alpha$ and $\Gamma_\alpha$ the value group of
$\nu_\alpha$. Recall that
$$
\dsp{\Gamma_\alpha \cong\frac{A(\alpha)\setminus
    \{0\}}{U(R_\alpha)}}
$$
and that the valuation $\nu_\alpha$ can be identified with the natural homomorphism
$$
\dsp{A(\alpha)\setminus \{0\} \to \frac{A(\alpha)\setminus
  \{0\}}{U(R_\alpha)}}.
$$
\medskip

By definition, we have a natural ring homomorphism
\begin{equation}
A\rightarrow R_\alpha\label{eq:hom}
\end{equation}
whose kernel is $\mathfrak{p}_\alpha$.  \medskip

Conversely, the point $\alpha$ can be reconstructed
from the ring $R_\alpha$ by specifying a certain number of sign
conditions (finitely many conditions when $A$ is noetherian)
(\cite{Baer}, \cite{Krull}, \cite{BCR} 10.1.10, p. 217).
\medskip

The valuation $\nu_\alpha$ has the following properties:

(1) $\nu_\alpha(A[\alpha]) \geq 0$

(2) If $A$ is an $R$-algebra then for any positive elements $y,z \in
A(\alpha)$,
\begin{equation}
\nu_\alpha(y) < \nu_\alpha(z) \implies y > Nz,\ \forall N \in
R. \label{val2}
\end{equation}
A $\nu_\alpha$-ideal of $A$ is the preimage in $A$ of an ideal of
$R_\alpha$. See \cite{Pre} or \cite{And}, \S II.3 for more information
on this subject.
\medskip

As pointed out in \cite{LMSS}, the points of $\sper\ A$ admit the
following geometric interpretation (see also \cite{Fuc}, \cite{Kap},
\cite{Pre}, p. 89 and \cite{PC} for the construction and properties of
generalized power series rings and fields).

\begin{deft} Let $k$ be a field and $\Gamma$ an ordered abelian
  group. The generalized formal power series field
  $k\left(\left(t^\Gamma\right)\right)$ is the field formed by
  elements of the form $\sum\limits_{\gamma\in\Gamma}a_\gamma  t^\gamma$,
  $a_\gamma\in k$ such that the set $\left\{\gamma\ \left|\
      a_\gamma\ne0\right.\right\}$ is well ordered.
\end{deft}
The field $k\left(\left(t^\Gamma\right)\right)$ is equipped with the
natural $t$-adic valuation $v$ with values in $\Gamma$, defined by
$v(f)=\inf\{\gamma\ |\ a_\gamma\neq0\}$ for $f=\sum\limits_\gamma
a_\gamma t^\gamma\in k\left(\left(t^\Gamma\right)\right)$. The valuation ring of
this valuation is the ring $k\left[\left[t^\Gamma\right]\right]$ formed by all
the
elements of $k\left(\left(t^\Gamma\right)\right)$ of the form
$\sum\limits_{\gamma\in\Gamma_+}a_\gamma  t^\gamma$. Specifying
a total ordering on $k$ and $\dim_{\F_2}(\Gamma/2\Gamma)$ sign
conditions defines a total ordering on
$k\left(\left(t^\Gamma\right)\right)$. In this ordering $|t|$ is
smaller than any positive element of $k$. For example, if $t^\gamma>0$ for all
$\gamma\in\Gamma$ then $f>0$ if and only if $a_{v(f)}>0$.

For an ordered field $k$, let $\bar k$ denote the real closure of $k$.
The following result is a variation on a theorem of Kaplansky
(\cite{Kap}, \cite{Kap2}) for valued fields equipped with a total ordering.

\begin{thm} \textnormal{\textbf{(\cite{PC}, p. 62, Satz 21)}} Let $K$ be a
  real valued field, with residue field $k$ and value group $\Gamma$.
  There exists an injection $K\hookrightarrow\bar
  k\left(\left(t^\Gamma\right)\right)$ of real valued fields.
\end{thm}
\medskip

Let $\alpha\in\sper\ A$. In view of (\ref{eq:hom}) and the Remark above, specifying a
point $\alpha\in\sper\ A$ is equivalent to specifying a total order of
$k_\alpha$, a morphism
\begin{equation}
A[\alpha]\to\bar k_\alpha\left[\left[t^{\Gamma_\alpha}\right]\right]\label{eq:curvette}
\end{equation}
and $\dim_{\F_2}(\Gamma_\alpha/2\Gamma_\alpha)$ sign conditions.\medskip

\noi We may pass to Zariski spectra to obtain morphisms
$$
\mbox{Spec}\ \left(\bar
  k_\alpha\left[\left[t^{\Gamma_\alpha}\right]\right]\right)\to\mbox{Spec}\
A[\alpha]\to\mbox{Spec}\ A,
$$
induced by the ring homomorphism (\ref{eq:curvette}) and the natural
surjective homomorphism $A \twoheadrightarrow A[\alpha]$, respectively.

In particular, if $\Gamma_\alpha=\Z$, we obtain a
\textbf{formal curve} in Spec $A$ (an analytic curve if the series are
convergent). This motivates the following definition:

\begin{deft} Let $k$ be an ordered field. A $k$\textbf{-curvette} on
  $\sper(A)$ is a morphism of the form
$$
\alpha:A\to k\left[\left[t^\Gamma\right]\right],
$$
where $\Gamma$ is an ordered group. A $k$\textbf{-semi-curvette} is a
$k$-curvette $\alpha$ together with
  a choice of the sign data $\sgn\ x_1$,..., $\sgn\ x_r$, where
  $x_1,...,x_r$ are elements of $A$ whose $t$-adic values induce an
  $\F_2$-basis of $\Gamma/2\Gamma$.
\end{deft}

We have thus explained how to associate to a point $\alpha$ of Sper\ $A$ a
$\bar k_\alpha$-semi-curvette. Conversely, given an
ordered field $k$, a $k$-semi-curvette $\alpha$ determines a prime
ideal $\mathfrak{p}_\alpha$ (the ideal of all the elements of $A$
which vanish identically on $\alpha$) and a total ordering on
$A/\mathfrak{p}_\alpha$ induced by the ordering of the ring
$k\left[\left[t^\Gamma\right]\right]$ of formal power series.\medskip

Below, we will often describe points in the real spectrum by
specifying the corresponding semi-curvettes.\medskip
\medskip

Let $\nu$ be a  valuation centered in a
regular local ring $A$ (see \S \ref{valua}), let $\Phi=\nu(A \setminus
\{0\})$; $\Phi$ is a well-ordered set. For an ordinal $\lambda<\Phi$, let
$\gamma_\lambda$ be the element of $\Phi$ corresponding to $\lambda$.

\begin{deft}
A system of approximate roots of $\nu$ is a well-ordered set of elements
$$
\mathbf{Q}=\{Q_i\}_{i \in \Lambda}\subset A,
$$
satisfying the following condition: \label{def-sar} for every $\nu$-ideal $I$
in $A$, we have
\begin{equation} \label{eq:valideal0S}
I=\left\lbrace\left. \prod_j Q_j^{\gamma_j}\ \right|\ \sum_j \gamma_j \nu(Q_j)
  \geq \nu(I) \right\rbrace A;
\end{equation}
furthermore, we require the set $\mathbf Q$ to be minimal in the sense of
inclusion among those satisfying (\ref{eq:valideal0S}).

A system of approximate roots of $\nu$ up to $\gamma_\lambda$ is a
well-ordered set of elements of $A$ satisfying (\ref{eq:valideal0S})
for all the $\nu$-ideals $I$ such that $\nu(I)\le\gamma_\lambda$.
\end{deft}

The main results of this paper are:
\begin{enumerate}
\item Given a regular local ring $(A,\mathfrak m,k)$, a valuation $\nu$ centered
at $A$, as above, and an element $\gamma_\lambda\in\Phi$ such that the
$\nu$-ideal determined by $\gamma_\lambda$ is $\mathfrak m$-primary, we
construct a system of approximate roots up to $\gamma_\lambda$.
\item We construct a system of approximate roots for $A$ and $\nu$ under the
assumption that $A$ is $\mathfrak m$ -adically complete.
\item In the situation of the Connectedness (or Definable Connectedness)
conjecure we describe certain subsets $C\subset\sper\ A$ by explicit formulae
in terms of approximate roots; we conjecture that these sets satisfy the
Connectedness (respectively, Definable Connectedness) conjecture.
\item In the special case $\dim\ A=2$, we use the above results and
constructions to prove the Definable Connectedness conjecture (and hence
\textit{a fortiori} the Pierce--Birkhoff conjecture). We also prove the
Connectedness conjecture in dimension 2, provided the ring $A$ is excellent.
\end{enumerate}
The paper is organized as follows. Sections \ref{valua} to \ref{standardform} are purely
valuation-theoretic; sections \ref{AppS} and \ref{App} are devoted to the construction
of a system of approximate roots.

The approximate roots $Q_i$ are constructed
recursively in $i$. Roughly speaking, $Q_{i+1}$ is the lifting to $A$
of the minimal polynomial equation satisfied by $\init_\nu Q_i$ over
$k\left[\left\{\init_\nu Q_j\right\}_{j<i}\right]$ in $\gr_\nu A$. In sections
\ref{valua} to \ref{standardform}, we
prove that such systems of approximate roots exist in two situations: first, for any
$\mathfrak{m}$-primary $\nu$-ideal $J$ there exists a system of
approximate roots up to $\nu(J)$; secondly, there exists a system of
approximate roots whenever $A$ is $\mathfrak{m}$-adically complete.

Once these valuation-theoretic tools are developed, we continue with
the program announced in \cite{LMSS} for proving the Pierce-Birkhoff
conjecture. We place ourselves in the situation of Conjectures \ref{conn} and
\ref{DCC}. In \S\ref{DSI} we describe the separating ideal $<\alpha,\beta>$ by
describing monomials in the approximate  roots (common to the valuations
$\nu_\alpha$ and $\nu_\beta$) which generate it. In section
\ref{section:connexe}, we give an explicit description of a set $C \subset
\sper A \setminus \{g_1\ldots g_s=0\}$, containing $\alpha$ and
$\beta$, which we conjecture to be connected. The set $C$ is described
in terms of a finite family of approximate roots, common to the valuations
$\nu_\alpha$ and $\nu_\beta$.

Finally, we prove the Definable connectedness conjecture and hence the
Pierce-Birkhoff conjecture for an arbitrary regular
2-dimensional local ring $A$; we also prove
Conjecture \ref{conn} assuming
that $A$ is excellent which provides a second proof of the
Pierce-Birkhoff conjecture in the case of excellent rings. The outline of the
proof of the two conjectures is as follows. First, we use a sequence of point
blowings up and Zariski's theory of complete ideals (recalled and refined in
\S\ref{Zariski}) to transform the set $C$ into a set $U$ of a very simple form,
which informally we call a quadrant. Namely, $U$ is the set of all the points
$\delta$ of $\sper\ A'$ (where $A'$ is a regular two-dimensional local ring
obtained after a sequence of blowings up with regular system of parameters
$x',y'$), centered at the origin, which induce a specified total order on $k$
and which satify the sign conditions $x'(\delta)>0$, $y'(\delta)>0$ (resp.
$x'(\delta)>0$). This is accomplished in \S\ref{Ctoquadrant}.

In the special case when $A'$ is essentially of finite type over a real closed
field the connectedness of $U$ is well  known and follows easily from the
results of \cite{BCR} (which allow to reduce connectedness of $U$ to that of a
quadrant in the usual Euclidean plane). However, for more general regular rings
this result seems to us to be new and non-trivial.

In \S\ref{connectedexcellent}, we use results from \cite{And} to reduce the
connectedness of $U$ to that of a quadrant in the usual Euclidean space,
assuming the ring $A$ is excellent. This completes the proof of the
connectedness conjecture for excellent regular 2-dimensional rings. In
\S\ref{proof} we prove the definable connectedness of $U$, without any
excellence assumptions, by using a new notion of a graph, associated to a
sequence of point blowings-up of a real surface.

Our proof is based on Madden's unpublished preprint \cite{Mad2}. As
well, we would like to acknowledge a recent paper by S. Wagner
\cite{W} which gives a proof of the Definable Connectedness and the
Pierce--Birkhoff conjecture in the case of smooth 2-dimensional
algebras of finite type over real closed fields.
\medskip

The overall structure of our proof is similar to that of \cite{Mad2}
and \cite{W}, with the following differences:
\begin{enumerate}
\item Here, we have tried to present a proof which should provide a
pattern for a general proof of the conjecture, that is, have a hope of
generalizing to higher dimensions. In particular, we went to great lengths to
phrase everything in terms of approximate roots rather than work directly with
connected sets as in \cite{Mad2} and \cite{W}.
\item We make no assumptions on the real closedness of the residue field of $A$
which introduces certain extra  complications.
\item Because we work with arbitrary regular two-dimensional rings, we have to
overcome a serious difficulty: proving that the ``quadrant'' $U$, defined
above, is connected. This is well known for algebras of finite type over a real
closed field (see, for example, \cite{BCR}) but as far as we can tell, for
general rings this result is new and non-trivial. Its proof occupies most of
section \ref{proof}.
\end{enumerate}

We thank the referee for his very careful reading of the manuscript and for
many useful suggestions which helped improve the paper.

\section*{\Large Part 1. Valuations and approximate roots.}
\setcounter{section}{1}
\setcounter{subsection}{0}
\subsection{Generalities on valuations.} \label{valua}
In this section we review some basic facts of valuation theory.

Let $A$ be a noetherian ring and $\nu:A \to \Gamma \cup \{\infty\}$ a
valuation centered at a prime ideal of $A$. Let $\Phi = \nu(A \setminus
\{0\}) \subset \Gamma$.

For each $\gamma \in \Phi$, consider the ideals
\begin{equation} \label{valua:eq1}
\begin{array}{ccl} P_\gamma &=& \left\{x \in A\ \left|\ \nu(x) \geq \gamma\right.\right\}
\\ P_{\gamma +}&=& \left\{x \in A\ \left|\ \nu(x) > \gamma\right.\right\}.
\end{array}
\end{equation}
$P_\gamma$ is called the $\nu$-\textbf{ideal} of $A$ of value
$\gamma$.

\begin{rek} It is easy to see that, as $A$ is noetherian, $\nu(A)$ is
  well-ordered.
\end{rek}

\noi \textbf{Notation.} If $I$ is an ideal of $A$ and $\nu$ a valuation
of $A$, $\nu(I)$ will denote $\min\{\nu(x) \ |\ x \in I\}$.
\medskip

We now define certain natural graded algebras associated to a
valuation. Let $A$, $\nu$ and $\Phi$ be as above. For $\gamma \in \Phi$,
let $P_\gamma$ and $P_{\gamma +}$ be as in (\ref{valua:eq1}). We
define $$\gr_\nu A = \bigoplus_{\gamma \in
    \Phi}\frac{P_\gamma}{P_{\gamma +}}.$$  The algebra $\gr_\nu(A)$ is
  an integral domain. For any element $f \in A$ with $\nu(f)=\gamma$,
  we may consider the natural image of $f$ in
  $\dsp{\frac{P_\gamma}{P_{\gamma +}}} \subset \gr_\nu(A)$. This image
  is a homogeneous element of $\gr_\nu(A)$ of degree $\gamma$, which
  we denote by $\init_\nu f$. The grading induces an obvious valuation
  on $\gr_\nu(A)$ with values in $\Phi$; this valuation will be
  denoted by ord.
\medskip

We end this section with the notion of a \textit{monomial}
valuation. Let $(A,\mathfrak{m},k)$ be a regular local ring,
and $\mathbf{u}=(u_1,\ldots,u_n)$ a regular system of parameters of
$A$. Let $\Phi$ be an ordered semigroup and let
$\beta_1,\ldots,\beta_n$ be strictly positive elements of $\Phi$. Let $\Phi_*$
denote the ordered semigroup, contained in $\Phi$, consisting of all
the  $\N_0$-linear combinations of $\beta_1,\ldots,\beta_n$. For
$\gamma \in \Phi_*$, let $I_\gamma$ denote the ideal of $A$, generated
by all the monomials $u^\alpha$ such that
$\sum\limits_{j=1}^n\alpha_j\beta_j \geq \gamma$ (we take $I_0 = A$). Let $x$ be
a non-zero element of $A$. Let $\Phi_x =\left\{\gamma \in \Phi_* \ \left|\ x \in
I_\gamma \right.\right\}$. Then it is not
difficult to prove that the set $\Phi_x$ contains a maximal element and there
exists a unique valuation $\nu$, centered at $\mathfrak m$, such that
\begin{equation} \label{valua:eq3}
\nu(u_j)=\beta_j,\ 1 \leq j \leq n
\end{equation} and
\begin{equation} \label{valua:eq4}
\nu(x) = \max \{\gamma
\in \Phi_x\},\ x \in A \setminus \{0\}.
\end{equation}
This valuation is called the
\textbf{monomial valuation} of $A$, associated to $\mathbf{u}$ and the
$n$-tuple $(\beta_1,\ldots,\beta_n)$. A valuation $\nu$, with values
in a group $\Gamma$, centered in $\mathfrak m$, is said to be
\textbf{monomial with respect to} $\mathbf{u}$ if there exist
$\beta_1,\ldots,\beta_n \in \Gamma_+$ such that (\ref{valua:eq4})
holds for all $x \in A \setminus \{0\}$.
\medskip

For further results on valuations, see also \cite{V5} or \cite{Zar}.
\medskip

\noi The following result is an immediate consequence of definitions:
\begin{prop} Let $G_\nu$ be the graded algebra associated to a valuation
$\nu:K\rightarrow\Gamma$, as above. Consider a sum of the form
$y=\sum\limits_{i=1}^sy_i$, with $y_i\in K$. Let $\beta=\min\limits_{1\le i\le
s}\nu(y_i)$ and
$$
S=\left\{\left.i\in\{1,\dots,n\}\ \right|\ \nu(y_i)=\beta\right\}.
$$
The following two conditions are equivalent:

(1) $\nu(y)=\beta$

(2) $\sum\limits_{i\in S}\init_\nu y_i\ne0$.
\end{prop}

\subsection{Approximate roots up to $\nu(J)$ for an
 $\mathfrak{m}$-primary ideal $J$}
\label{AppS}

Let $A$ be a regular local ring of dimension $n$, $\mathfrak{m}$ its maximal
ideal, $\dsp{k= \frac{A}{\mathfrak{m}}}$,
$\mathbf{u}=(u_1,\ldots,u_n)$ a regular
system of parameters and
$$
\nu:A\setminus\{0\}\to\Gamma
$$
a valuation, centered in $\mathfrak{m}$ (this means $\nu(\mathfrak{m})>
0$).

Let $\mathbf{1}=\nu(\mathfrak{m})=\min\{\gamma\in\Phi\ |\ \gamma>0\}$ and
$\Phi_1=\{\gamma\in\Phi\ |\ \exists
a\in\N;\gamma<a\cdot\mathbf{1}\}$. For the sake of simplicity, we will
write $a$ instead of $a\cdot\mathbf{1}$. We shall study the structure
of $\nu$-ideals $P_\gamma$ where $\gamma\in\Phi$.

If $\nu$ were monomial with respect to $\mathbf{u}$ then $\init_\nu
u_1,\ldots,\init_\nu u_n$ would generate $\gr_\nu A$ as a
$k$-algebra. We are interested in analyzing valuations
which are not necessarily monomial. We fix an $\mathfrak{m}$-primary valuation
ideal $J$. The purpose of sections \ref{AppS} and \ref{Sstandardform} is
to construct a system of approximate roots up to $\nu(J)$, that is, a finite
sequence of elements
$\mathbf{Q}=\{Q_i\}_{i \in \Lambda}$ of $A$ such that for every
$\nu$-ideal $I$ in $A$ containing $J$ we have
\begin{equation}
I = \left\lbrace\left. \prod\limits_j Q_j^{\gamma_j} \ \right|\
\sum\limits_j\gamma_j\nu(Q_j)\ge\nu(I)\right\rbrace A \label{eq:validealS}
\end{equation}
(in particular, the images $\init_\nu Q_i$ of the $Q_i$ in $\gr_\nu A$
generate $\gr_\nu A$ as a $k$-algebra up to degree $\nu(J)$). In this construction, each $Q_{i+1}$ will be
described by an explicit formula (given later in this section) in
terms of $Q_1,...,Q_i$.
\medskip

The earliest precursor of approximate roots appears in a series of
papers by Saunders MacLane and O.F.G. Schilling \cite{Mac1},
\cite{Mac2} and \cite{Mac3}. In dimension 2, they were defined
globally in $k[x,y]$ by S. Abhyankar and T. T. Moh (\cite{AM1}, \cite{AM2}) and locally by
M. Lejeune-Jalabert \cite{LJ}. See also the papers \cite{Kuo1} and
\cite{Kuo2} by T. C. Kuo, \cite{GT} by R. Goldin and B. Teissier and
\cite{Spi1} by M. Spivakovsky, \cite{GAS} by F.J. Herrera Govantes,
M.A. Olalla Acosta, M. Spivakovsky, \cite{V1}-\cite{V4}  by Michel Vaquié. We
also refer the reader to the paper \cite{Te} by B. Teissier for a
different approach to the theory of approximate roots in higher
dimensions.
\medskip

Let $\dsp{k=\frac{A}{\mathfrak{m}}=\frac{A}{\mathfrak{m}_\nu\cap A}}$
be the residue field
of $A$. Fix an isomorphism $\displaystyle{\frac{A}{J} \cong
  \frac{k[u_1,\ldots,u_n]}{J_0}}$, where $J_0$ is an ideal of
$k[u_1,\ldots,u_n]$. In this way, we will
view $k$ as a subring of $A/J$.

We fix, once and for all, a section $k\to A$ of the natural map
$A\rightarrow k$ which composed with the natural map $\dsp{A \to
\frac{A}{J}}$ maps $k$ isomorphically onto its image in
$\dsp{\frac{A}{J}}$. The image of $k$ in $A$ will be denoted by $\mathbf{k}$.
\medskip

According to Definition \ref{def-sar}, we are
looking for a finite set of
elements $\mathbf{Q}=\{Q_i\}_{i \in \Lambda}$, $Q_i \in A$ satisfying
(\ref{eq:validealS}).
\begin{rek} This means, in particular, that the initial forms
$\init_\nu(Q_1),\init_\nu(Q_2),\ldots$ generate $\gr_\nu(A)$, up to
degree $\nu(J)$. In other words, we want to build $\mathbf{Q}$ such that, for
$f\in A$, we have $\init_\nu(f)\in k[\init_\nu \mathbf{Q}]$
provided $\nu(f)\le\nu(J)$.
\end{rek}

Since $J$ is an $\mathfrak{m}$-primary ideal, there are only finitely many
elements of $\Phi$ less than or equal to $\nu(J)$. We proceed by induction on
the finite set $\{\gamma \in \Phi\ |\ \gamma\le\nu(J)\}$.

\begin{deft}\label{genmon}
Let
$E$ be an ordered set of elements of $A$. A generalized monomial
$\mathbf{Q}^{\alpha}$ in $E$ is a formal expression $$\mathbf{Q}^{\alpha} = \prod_{Q \in E}
Q^{\alpha_Q}$$ where $\alpha_Q \in \N$ and $\alpha_Q = 0$ for
all $Q$ outside of a finite subset of $E$.
\end{deft}

We view the set $\N^E$ as being ordrered lexicographically and order the
set of generalized monomials by the lexicographical
order of the pairs $(\nu(\mathbf{Q}^\alpha),\alpha)$.
\medskip

The semigroup $\Phi$ is well ordered. For a natural number $\lambda$,
$\gamma_\lambda$ will denote the $\lambda$-th element of $\Phi$.
\medskip

\noi We start by choosing a coordinate system adapted to the situation.

\begin{deft}
Take $j \in \{2,\ldots,n\}$. We say that $u_j$ is $(\nu,J)$-prepared
if either $u_j \in J$ or there does not exist $f \in A$ such that
\begin{eqnarray}
\init_\nu u_j&=& \init_\nu f\quad\text{ and}\label{eq:equalityin}\\
f \mod J &\in&
\frac{k[u_1,\ldots,u_{j-1}]}{k[u_1,\ldots,u_{j-1}] \cap J_0}.\label{eq:j-1variables}
\end{eqnarray}
The coordinate system $\mathbf{u} = \{u_1,\ldots,u_n\}$ is
$(\nu,J)$-prepared if $u_j$ is $(\nu,J)$-prepared for all $j \in
\{2,\ldots,n\}$.
\end{deft}

\begin{prop}
There exists a $(\nu,J)$-prepared coordinate system.
\end{prop}

\noi Proof: \ We construct a $(\nu,J)$-prepared coordinate system
recursively in $j$. \ Assume that
$u_2,\ldots,u_{j-1}$ are already $(\nu,J)$-prepared, but $u_j$ is not. Take $f \in A$
satisfying (\ref{eq:equalityin}) and (\ref{eq:j-1variables}).

Let $\tilde{u}_j= u_j-f$; then $\nu(\tilde{u}_j) >
\nu(u_j)$.

Since there are only finitely many elements of $\Phi$ less than $\nu(J)$, after finitely many
repetitions of the above procedure, we may assume that $u_j$ is
$(\nu,J)$-prepared. This completes the proof by induction on $j$. $\Box$
\medskip

\noi We construct, recursively in $\lambda$, two finite
ordered sets $\Lambda(\gamma_\lambda)$ and
$\Theta(\gamma_\lambda)$
with
$$
\Lambda(\gamma_\lambda) \subset \bigcup_{\lambda' < \lambda}
\Theta(\gamma_{\lambda'}),
$$
and a total ordering of the set
$\Lambda(\gamma_\lambda) \cup \Theta(\gamma_{\lambda-1})$, compatible with the
orders on $\Lambda(\gamma_\lambda)$ and $\Theta(\gamma_{\lambda-1})$.
We do not impose a total order on the union $\bigcup_{\lambda' < \lambda}
\Theta(\gamma_{\lambda'})$.
At each step we define additional finite ordered sets
\begin{equation}
\mathcal{V}(\gamma_\lambda) \subset \Psi(\gamma_\lambda) \subset
\Lambda(\gamma_\lambda),\label{eq:ordercompatible}
\end{equation}
where the inclusions in (\ref{eq:ordercompatible}) are inclusions of ordered
sets. Both
collections of sets $\Lambda(\gamma_\lambda)$ and $\Nu(\gamma_\lambda)$ will
be increasing with $\lambda$.
A typical element of each of
those sets will have the form $(Q,\mbox{Ex}(Q))$ where $Q \in A$ and
$\mbox{Ex}(Q)$ is a sum of monomials in
$\Lambda(\gamma_\lambda) \cup \Theta(\gamma_{\lambda-1})$,
written in the increasing order according to the on monomials, defined above.


Given an element
$(Q,\mbox{Ex}(Q)) \in \Lambda(\gamma_\lambda) \cup
\Theta(\gamma_\lambda)$, $Q$ is called an \textit{approximate root} and
$\mbox{Ex}(Q)$ is called the \textit{expression} of $Q$. In what follows, we adopt the convention
$$
\Theta(\gamma_\lambda)=\mathcal{V}(\gamma_\lambda)=\Psi(\gamma_\lambda)=\Lambda(\gamma_\lambda)=\emptyset
$$
whenever $\lambda<0$.
\medskip

For a natural number $\ell$, $\gamma_\ell \leq \nu(J)$, and for
$(Q,\mbox{Ex}(Q)) \in \Lambda(\gamma_\ell) \cup \Theta(\gamma_\ell)$, let
 $\mbox{In}\ Q$ denote the smallest monomial of $\mbox{Ex}(Q)$. Let
$$
\mbox{In}(\ell) = \left\{\left.\alpha \in \N^{\Nu(\gamma_\ell)} \ \right|\
\exists (Q,\mbox{Ex}(Q))
 \in \Lambda(\gamma_\ell) \mbox{ such that } \mathbf{Q}^\alpha = \mbox{In}\ Q\right\}.
$$
\medskip

\begin{thm}\label{thmS-aproot}
For a natural number $\lambda$, $\gamma_\lambda\le \nu(J)$, there exist
finite ordered sets
$$
\mathcal{V}(\gamma_\lambda) \subset
  \Psi(\gamma_\lambda) \subset \Lambda(\gamma_\lambda)
  $$
  and $\Theta(\gamma_ \lambda)$ (and a total ordering of
$\Lambda(\gamma_\lambda) \cup \Theta(\gamma_{\lambda-1})$)
consisting of elements $(Q,\mbox{Ex}(Q))$, with $Q \in A$ and
$\mbox{Ex}(Q)$ a sum of monomials in
$\mathcal{V}(\gamma_\lambda) \cup \Theta(\gamma_{\lambda-1})$, increasing with
respect to the given order on
monomials, and having the following properties:
\begin{eqnarray}
 \nu(Q) < \gamma_\lambda \mbox{ whenever } (Q,\mbox{Ex}(Q)) \in
 \Lambda(\gamma_\lambda) \label{thmS1}\\
  \nu(Q) \geq \gamma_\lambda \mbox{ whenever } (Q,\mbox{Ex}(Q)) \in
 \Theta(\gamma_\lambda). \label{thmS2}
\end{eqnarray}

Moreover, for any $(Q,\mbox{Ex}(Q)) \in
\Lambda(\gamma_\lambda)$, any monomial
$\mathbf{Q}^\alpha$ appearing in $\mbox{Ex}(Q)$ is a monomial in
$\Nu(\gamma_{\lambda-1})$ provided $Q \notin \{u_1,\ldots,u_n\}$. For
any $(Q,\mbox{Ex}(Q)) \in \Theta(\gamma_\lambda)$, any monomial
$\mathbf{Q}^\alpha$ appearing in $\mbox{Ex}(Q)$ is a monomial in
$(\Nu(\gamma_{\lambda+1}) \cap \Theta(\gamma_{\lambda-1})) \cup
\Nu(\gamma_\lambda)$ provided $Q \notin \{u_1,\ldots,u_n\}$.
An element
$$
(Q,\mbox{Ex}(Q)) \in\Psi(\gamma_\lambda) \cup  \Theta(\gamma_\lambda)
$$
is completely determined by $\mbox{In}\ Q$.
\end{thm}

\noi Proof: We proceed by induction on $\lambda$.

First define $\Psi(\mathbf{1}) = \Lambda(\mathbf{1}) =
\emptyset$ and $\Theta(\mathbf{1})  = \{(u_1,u_1),\ldots,(u_n,u_n)\}$
where we assume
$$
\nu(u_1) \leq \nu(u_2) \leq \cdots \leq \nu(u_n).
$$
We define the total ordering on $\Theta(\mathbf{1})$ by $(u_1,u_1) < (u_2,u_2) <
\cdots < (u_n,u_n)$.
\medskip

Let $\lambda>0$ be a natural number such that $\gamma_\lambda \leq
\nu(J)$. Assume that for each $\ell <
\lambda$ we have constructed sets $\mathcal{V}(\gamma_\ell) \subset
  \Psi(\gamma_\ell) \subset
\Lambda(\gamma_\ell)$
and $\Theta(\gamma_\ell)$ having the
properties required in the theorem.

Let \begin{equation} \Lambda(\gamma_\lambda)= \Lambda(\gamma_{\lambda-1}) \cup
  \{(Q,\mbox{Ex}(Q)) \in \Theta(\gamma_{\lambda-1}) \ |\ \nu(Q) < \gamma_\lambda\}.
\end{equation}

\begin{deft}\label{Spred-ines} An element $(Q,\mbox{Ex}(Q)) \in
  \Lambda(\gamma_\lambda)$ is an \textbf{inessential predecessor} of an approximate root
  $(Q',\mbox{Ex}(Q')) \in \Lambda(\gamma_\lambda)$ if $\mbox{Ex}(Q') =
  \mbox{Ex}(Q) + \sum\limits_\alpha c_\alpha \mathbf{Q}^\alpha$, where
  $c_\alpha\in k$ and the $\mathbf{Q}^\alpha$ are monomials in
  $\Nu(\gamma_\lambda)$.
\smallskip

An element $(Q,\mbox{Ex}(Q)) \in \Lambda(\gamma_\lambda)$ is said to be
\textbf{essential at the level $\gamma_\lambda$} if $Q$ is not an inessential
predecessor of an element of $\Lambda(\gamma_\lambda)$.
\end{deft}


Let $\Psi(\gamma_\lambda)$ be the subset of $\Lambda(\gamma_\lambda)$
consisting of all the essential roots at the level
$\gamma_\lambda$. Let $\mathcal{V}(\gamma_\lambda)$ be the
subset of $\Psi(\gamma_\lambda)$ consisting of all $(Q,\mbox{Ex}(Q))$ such
that $\init_\nu(Q)$ does not belong to the $k$-vector space of
$G=\gr_\nu(A)$ generated by the set $\{\init_\nu \mathbf{Q}^\gamma \}$
where $\mathbf{Q}^\gamma$ runs over the set of all the generalized
monomials on roots preceding $Q$ in the above ordering.

\noi We extend the total ordering from $\Lambda(\gamma_{\lambda-1})$ to
$\Lambda(\gamma_\lambda)$ by postulating that
$\Lambda(\gamma_{\lambda-1})$ is the initial segment of
$\Lambda(\gamma_\lambda)$. Moreover, we extend this order to
$\Lambda(\gamma_\lambda) \cup \Theta(\gamma_{\lambda-1})$ by postulating that
$\Lambda(\gamma_\lambda)$ is the initial segment of $\Lambda(\gamma_\lambda)
\cup \Theta(\gamma_{\lambda-1})$.
\medskip

For a natural number $\ell$, let $E(\ell) =
\mbox{In}(\ell)+ \N^{\mathcal{V}(\gamma_\ell)} \subset
\N^{\mathcal{V}(\gamma_\ell)}$.
\medskip

Now consider the ordered set
$\{\mathbf{Q}^{\alpha_1},\ldots,\mathbf{Q}^{\alpha_s} \}$ of
monomials
\begin{equation}\label{eq:listmonomgammalambda}
\mathbf{Q}^\alpha = \prod Q^{\alpha_Q},\ (Q,\mbox{Ex}(Q)) \in
\mathcal{V}(\gamma_\lambda)
\cup \left\{\left.(Q,\mbox{Ex}(Q)) \in \Theta(\gamma_{\lambda-1}) \ \right|\
\nu(Q) =
\gamma_\lambda \right\}
\end{equation}
of value $\gamma_\lambda$, such that the natural projection of $\alpha$ to
$\N^{\mathcal{V}(\gamma_\lambda)}$ does not belong to $E(\lambda)$.
\smallskip

Let $i_1 = \max \left\{ i \in \{ 1,\ldots,s \}\ \left|\
\init_\nu(\mathbf{Q}^{\alpha_i})  \in
\sum\limits_{j=i+1}^s k \ \init_\nu(\mathbf{Q}^{\alpha_j})\right.\right \}$ and
consider the unique
relation  $\init_\nu(\mathbf{Q}^{\alpha_{i_1}}) - \sum\limits_{j={i_1}+1}^s
\sur{c_{1j}} \
\init_\nu(\mathbf{Q}^{\alpha_j})=0$.  Let $P_1 =
  \mathbf{Q}^{\alpha_{i_1}} - \sum\limits_{j={i_1}+1}^s
c_{1j} \mathbf{Q}^{\alpha_j}$ where $c_{1j} \in \mathbf{k}$ is the image of
$\sur{c_{1j}}$ under the chosen section $k \to A$.

Let $i_2 = \max\left\{i \in \{1,\ldots,i_1-1\}\ \left|\
\init_\nu(\mathbf{Q}^{\alpha_i}) \in \sum\limits_{j=i+1}^s k\
\init_\nu(\mathbf{Q}^{\alpha_j})\right.\right\}$ and, as before, consider
the unique $\displaystyle{P_2 = \mathbf{Q}^{\alpha_{i_2}} -
  \sum_{\stackrel{j={i_2}+1}{j \neq i_1}}^s c_{2j}
\mathbf{Q}^{\alpha_j}}$ such that the vector
$(\alpha_j)_{j=i_1+1,\ldots,s},\ \sur{c_{2j}} \neq
  0$, is minimal in the lexicographical order. We continue in this way and
define $P_3,\ldots, P_t$.
\smallskip

Let
\begin{equation}\Theta(\gamma_\lambda) =
\left\{\left.(Q,\mbox{Ex}(Q)) \in \Theta(\gamma_{\lambda-1}) \ \right|\ \nu(Q) \geq
\gamma_\lambda\right\}
\cup \left\{(P_1,\mbox{Ex}(P_1)),\ldots,(P_t,\mbox{Ex}(P_t))\right\}
\end{equation}
where
\begin{equation} \label{Sformul2}
\mbox{Ex}(P_j) = \mbox{Ex}(Q) - \sum_k c_{jk} \mathbf{Q}^{\alpha_k}
\end{equation}
if $\mathbf{Q}^{\alpha_{i_j}}=Q$ with
$(Q,\mbox{Ex}(Q)) \in \left\{\left.(Q,\mbox{Ex}(Q)) \in \Theta(\gamma_{\lambda-1}) \ \right|\ \nu(Q) =
\gamma_\lambda \right\}$ and
\begin{equation}\label{Sformul1} \mbox{Ex}(P_j) =
\mathbf{Q}^{\alpha_{i_j}} - \sum_k c_{jk}
\mathbf{Q}^{\alpha_k}
\end{equation}
 otherwise.

We define the order on $\Theta(\gamma_\lambda)$ by
$\Theta(\gamma_{\lambda-1}) <
\left\{(P_1,\mbox{Ex}(P_1)),\ldots,(P_t,\mbox{Ex}(P_t))\right\}$ and
$(P_1,\mbox{Ex}(P_1)) < \cdots < (P_t,\mbox{Ex}(P_t))$.
\medskip
\begin{rek}
Note that, because the coordinate system is prepared, $u_1,\ldots,u_n$
are always essential.
\end{rek}

\begin{rek} \label{rk-continu}
Suppose given two approximate roots $Q_1$ and $Q_2$ such that
$$
In(Q_1)=In(Q_2)=\mathbf{Q}^\alpha
$$
and suppose that $Q_1$ appears before $Q_2$ in the process of construction of
the approximate roots decribed above. Because of the uniqueness of the
construction of the $P_i$'s above, we have
$$
\nu(Q_2) > \nu(Q_1).
$$
Now, if $\nu(\mathbf{Q}^\alpha) = \gamma_\ell$, then $\alpha \in E(\ell)$, so
the only way the monomial $\mathbf{Q}^\alpha$ can appear as an
initial form of $Q_2$ is when $P_k=Q'+\sum c_j\mathbf{Q}^{\alpha_j}$ where
$In(Q') = \mathbf{Q}^\alpha$ and then $\nu(Q') < \nu(Q_2)$. Then,
either $\nu(Q')= \nu(Q_1)$ and so $Q'=Q_1$ because of the uniqueness in
the construction process, or $\nu(Q') > \nu(Q_1)$, but we conclude by
descending induction that $Q_2=Q_1+\sum c_j\mathbf{Q}^{\alpha_j}$ and
$\mbox{Ex}(Q_2) = \mbox{Ex}(Q_1) + \sum c_j\mathbf{Q}^{\alpha_j}$.
\end{rek}

So finally, the expression of an approximate root has the form
\begin{equation}\label{Seq:apprroot}
\mbox{Ex}(Q) = \mathbf{Q}^{\alpha_{i_j}} + \sum_k a_k
\mathbf{Q}^{\alpha_k}
\end{equation}
the sum being  written in the increasing order of the monomials.
\medskip

\begin{rek} \label{Sunic} This construction is very
  similar to finding a basis of the space of relations by row reduction.
\end{rek}

\begin{rek} We just showed that there is a one to one
   correspondence $Q \leftrightarrow \mbox{In}(Q)$ between the
   approximate roots $Q \in \Psi(\gamma_\ell)$ and the set of
   monomials which are the first term of the expression $\mbox{Ex}(Q)$
   of such an approximate root $Q$. Let us denote by
   $\mathbf{M}(\ell)$ \label{Sem-de-ell} the set of those monomials.
\end{rek}

The last part of the theorem holds by construction. $\Box$

\subsection{Standard form up to $\nu(J)$}\label{Sstandardform}

Consider the integer $\lambda$ such that $\gamma_\lambda =
\nu(J)$. Assume that the system of coordinates $u$ of $A$ is
$(\nu,J)$-prepared.
\begin{deft}
A monomial in $\Psi(\gamma_\lambda)\cup \Theta(\gamma_\lambda)$  is called
\textbf{standard with respect to} $\lambda$ if all the approximate roots
appearing in it belong to $\mathcal{V}(\gamma_\lambda)$ and it is not divisible
by any $In(Q)$ where $Q$ is
an approximate root in
$\left(\Psi(\gamma_\lambda)\cup \Theta(\gamma_\lambda)\right)
\setminus\{(u_1,u_1),\dots,(u_n,u_n)\}$.
\end{deft}

\begin{deft} Let $f \in A$ and let $\ell$ be a positive integer,
  $\ell\le\lambda$. An expression of the form
  $$
  f = \sum c_\alpha \mathbf{Q}^\alpha,
  $$
  where the $\mathbf{Q}^\alpha$ are monomials in
  $\Psi(\gamma_\lambda)\cup \Theta(\gamma_\lambda)$, written in the
  increasing order, is \textbf{a standard form of level} $\gamma_\ell$ with
  respect to $\lambda$ if for all $\gamma' < \gamma_\ell$ and
for all $\alpha$ such that $\nu(\mathbf{Q}^\alpha)= \gamma'$ and $c_\alpha\ne0$,
$\mathbf{Q}^\alpha$ is a standard monomial with respect to $\lambda$.
\end{deft}

We now construct, by induction on $\ell$, a standard form of $f$ of
level $\gamma_\ell$. We will write this standard form as
$$
f = f_\ell + \sum c_\alpha \mathbf{Q}^\alpha
$$
where, for all $\alpha$, $\mathbf{Q}^\alpha$ is a generalized monomial
in $\Psi(\gamma_\lambda)\cup \Theta(\gamma_\lambda)$,
$\nu(\mathbf{Q}^\alpha) \geq
\gamma_\ell$ and $f_\ell$ is a sum of standard monomials in
$\mathcal{V}(\gamma_\lambda)$ of value strictly less than
$\gamma_\ell$.

To start the induction, let $f_0=0$. The \textit{standard form of $f$
  of level 0 with respect to $\lambda$} will be its expansion $f =
f_0 + \sum c_\alpha \mathbf{u}^\alpha$ as a formal
power series in the $u_i$, with the monomials written in the increasing order
according to the monomial order defined above,.
\medskip

Let $\ell$ be a natural number, $\ell<\lambda$. Let us define $f_{\ell +1}$ and the
standard form of $f$ of level $\gamma_{\ell +1}$ as follows. Assume we already
have an expression $f = f_\ell + \sum c_\alpha \mathbf{Q}^\alpha$ with
$\nu(\mathbf{Q}^\alpha) \geq
\gamma_\ell$, for all $\alpha$, and the value of any monomial of
$f_\ell$ is strictly less than $\gamma_\ell$.

Take the homogeneous part of $\sum c_\alpha \mathbf{Q}^\alpha$ of
value $\gamma_\ell$, with the monomials arranged in the increasing order, and
consider the first monomial $\mathbf{Q}^\alpha$ which
is not standard. Since $\mathbf{Q}^\alpha$ is not standard, one of the
following two conditions holds:
\begin{enumerate}
\item There exists an approximate
root $Q\in\left(\Psi(\gamma_\lambda)\cup \Theta(\gamma_\lambda)\right)
\setminus\{(u_1,u_1),\dots,(u_n,u_n)\}$ such that $In(Q)$ divides
$\mathbf{Q}^\alpha$. Write $Q=In(Q)+\sum c_\beta  \mathbf{Q}^\beta$ and replace
$In(Q)$ by $Q-\sum c_\beta \mathbf{Q}^\beta$ in
$\mathbf{Q}^\alpha$.
\item There exists $Q \in \Psi(\gamma_\lambda) \setminus
\mathcal{V}(\gamma_\lambda)$ which divides $\mathbf{Q}^\alpha$. Since $Q
\notin \mathcal{V}(\gamma_\lambda)$, there exists
$$
Q' \in \Psi(\gamma_\lambda)\cup\Theta(\gamma_\lambda)
$$
of the form $Q'=Q+\sum\limits_\delta
d_\delta \mathbf{Q}^\delta$ where $\mathbf{Q}^\delta$ are monomials in
$\mathcal{V}(\gamma_\lambda)$ of value greater than or equal to $\gamma_\ell$. Replace $Q$ by $Q'-\sum\limits_\delta d_\delta \mathbf{Q}^\delta$.
\end{enumerate}
In both cases, those changes introduce new monomials, but either
they are of value strictly greater than $\gamma_\ell$ or they are
of value exactly $\gamma_\ell$ but greater than
$\mathbf{Q}^\alpha$ in the monomial ordering. We repeat this
procedure as many times as we can. After a finite number of steps, no
more changes are available at level $\gamma_{\ell+1}$. Then, let
$f_{\ell+1} = f_\ell + \sum d_\beta \mathbf{Q}^\beta$ with
$\nu(\mathbf{Q}^\beta) = \gamma_\ell$, so that $f = f_{\ell+1} +
\sum c_\alpha \mathbf{Q}^\alpha$ where $\nu(\mathbf{Q}^\alpha) >
\gamma_\ell$.
\medskip

The expression thus constructed satisfies the definition of standard
form of level $\gamma_{\ell+1}$ because all the non-standard monomials
$\mathbf{Q}^\alpha$ of value less than or equal to $\gamma_\ell$ have been eliminated.
\begin{prop}\label{standard} Let
$$
f = f_\ell + \sum c_\alpha \mathbf{Q}^\alpha
$$
be a standard form of $f$ of level $\gamma_\ell$ and
$\gamma<\gamma_\ell$ an element of $\Phi$. Then
$\sum\limits_{\nu(\mathbf{Q}^\beta)=\gamma} c_\beta
\mathbf{Q}^\beta \notin P_{\gamma +}$.
\end{prop}

\noi Proof : We give a proof by contradiction. Suppose there exists a
relation of the form
\begin{equation} \label{eq:saut}
\sum\limits_{\nu(\mathbf{Q}^\beta)=\gamma} c_\beta
\mathbf{Q}^\beta \in P_{\gamma+}.
\end{equation}
Let $\mathbf{Q}^\alpha$ be the smallest monomial on the left hand side
of (\ref{eq:saut}). By construction of approximate roots, there exists
a finite collection $Q_1,\ldots,Q_s \in \Lambda(\gamma+)\cup
\theta(\gamma+)$ and generalized monomials
$\mathbf{Q}^{\omega_1},\ldots,\mathbf{Q}^{\omega_s}$ such that
$\sum\limits_{i=1}^s Q_i\mathbf{Q}^{\omega_i} =
\sum\limits_{\nu(\mathbf{Q}^\beta)=\gamma}c_\beta\mathbf{Q}^\beta$.

\noi There exists $i\in \{1,\ldots,s\}$ such that one of the two
conditions holds : either
$$
\mathbf{Q}^\alpha =\mathbf{Q}^{\omega_i}\cdot In(Q_i)
$$
or
$$
Q_i= Q'_i + \sum_\epsilon b_\epsilon \mathbf{Q}^\epsilon,\ Q'_i \in
\Lambda(\gamma+) \setminus \Psi(\gamma+).
$$
In either case, the monomial $\mathbf{Q}^\alpha$ is
not standard, which gives the desired contradiction.
$\Box$
\medskip

For each $\ell$, the part $f_\ell$ of a standard form of $f$ of
level $\gamma_\ell$ is uniquely determined. This is a straightforward
consequence of the Proposition.

As a consequence of Proposition \ref{standard}, note that if
$\gamma_\ell > \nu(f)$ then $\nu(f)$
equals the smallest value of a monomial appearing in
the standard form of $f$ of level $\gamma_\ell$.

\begin{thm} \label{Sgraded-gen}
(1) Take $\gamma \in \Phi$, $\gamma < \gamma_\lambda$. Then
$\dsp{\frac{P_\gamma}{P_{\gamma +}}}$ is generated as a $k$-vector
space by $\{\init_\nu \mathbf{Q}^\beta \}$ where $\mathbf{Q}^\beta$
runs over the set of all the standard monomials with respect to $\lambda$,
satisfying $\nu(\mathbf{Q}^\beta)=\gamma$.

(2) The part of the graded $k$-algebra $\gr_\nu(A)$
  of degree strictly less than $\gamma_\lambda$ is generated by the
  initial forms of the approximate roots of $\mathcal{V}(\gamma_\lambda)$.
\end{thm}

\noi Proof : Take an element $\gamma\in\Phi$, $\gamma<\gamma_\lambda$. Let $h
\in P_\gamma/P_{\gamma^+}$ be a  homogeneous
element of degree $\gamma$ of $\gr_\nu(A)$ and let $f \in P_\gamma$ be
such that $\init_\nu(f) = h$. Let $\sum c_\beta
\mathbf{Q}^\beta$ denote the homogeneous part of least value of a
standard form of $f$ of level $\gamma_\lambda$. Then the initial form of
$f$ is $\sum  \init_\nu(c_\beta\mathbf{Q}^\beta)$. $\Box$
\medskip

\noi\textbf{The Alvis--Johnston--Madden example.} Let $\alpha$ be
the point of $\sper(R[x,y,z])$ given by the curvette $x(t)=t^6,\
y(t)=t^{10}+ut^{11},z(t)=t^{14}+t^{15}$ where $u$ is some fixed
element of $R$ with $u>2$. Let $J$ be a $\nu_\alpha$-ideal of value greater than
or equal to 37.

\noi The calculation of the first few approximate roots gives
\begin{eqnarray}
Q_1&=&x,\\
Q_2&=&y,\\
Q_3&=&z,\\
Q_4&=&y^2-xz=(2u-1)t^{21}+u^2t^{22},\ \nu(Q_4)= 21\\
Q_5&=&yz-x^4=(u+1)t^{25}+ut^{26},\ \nu(Q_5) = 25\\
Q_6&=&z^2-x^3y=(2-u)t^{29}+t^{30},\ \nu(Q_6) = 29\\
Q_7^{(31)}&=&yQ_4-\alpha(u)xQ_5,\ \alpha(u)= (2u-1)/(u+1), \
\nu\left(Q_7^{(31)}\right) = 32\\
Q_7^{(32)}&=&yQ_4-\alpha(u)xQ_5-\beta(u)x^3z,\
\nu\left(Q_7^{(32)}\right)=33
\end{eqnarray}
\begin{eqnarray}
Q_7^{(33)}&=&yQ_4-\alpha(u)xQ_5-\beta(u)x^3z-\gamma(u)x^2Q_4 \\
Q_7^{(34)}&=&yQ_4-\alpha(u)xQ_5-\beta(u)x^3z-\gamma(u)x^2Q_4-\delta(u)x^4y
\\
Q_7^{(35)}&=&yQ_4-\alpha(u)xQ_5-\beta(u)x^3z-\gamma(u)x^2Q_4-\delta(u)x^4y-
\epsilon(u)xQ_6\\
Q_8^{(35)}&=&zQ_4+\zeta(u)xQ_6\\
Q_9^{(35)}&=&yQ_5+\eta(u)xQ_6,
\end{eqnarray}
where $\beta(u)$, $\gamma(u)$, $\delta(u)$, $\epsilon(u)$, $\zeta(u)$, $\eta(u)$
are functions of $u$ which can  be calculated explicitly.

The elements listed above belong to $\Lambda(37)$; we chose
to index them as $Q_i^{(j)}$. In this notation, the approximate root
$Q_i^{(j)}$ is an inessential predecessor of $Q_i^{(j+1)}$ whenever
$Q_i^{(j+1)}$ is defined.
\smallskip

\noi We also note the relation $xQ_6-yQ_5+zQ_4 = 0$, which is the simplest
example of a syzygy, an important  phenomenon, responsible for much of the
difficulty of the Pierce--Birkhoff conjecture.
\medskip

In the same vein, we can describe the standard form of different levels of an
element of $A$, say for instance,
\begin{equation}\label{fermat}
f=x^3+y^3+z^3
\end{equation}
(which is a standard form of level 0). For $\gamma\le30$, the
standard form of $f$ of level $\gamma$ is given by (\ref{fermat}). Then, as
$y^2 \in E(8)$ (this is so because $21$ is the eighth positive element of the
value semigroup $\Phi$), we
replace $y^3$ by $y(Q_4+xz)$ to obtain
\begin{equation}\label{fermat1}
f=x^3+yQ_4+xyz+z^3.
\end{equation}
Since
$yz\in E(11)$ (note that $25$ is the eleventh positive element of the
value semigroup $\Phi$), we replace $xyz$ in (\ref{fermat1}) by $xQ_5+x^5$,
to obtain the standard form of level 31:
\begin{equation}
f= x^3+x^5+yQ_4+xQ_5+z^3\label{fermat2}
\end{equation}
(the monomials being written in the order of increasing values 18, 30, 31, 31,
42). Next, we replace $yQ_4$ by  $\alpha(u)xQ_5$ in (\ref{fermat2}), so the
standard form of levels 32, 33, 34 and 35  is given by
$$
f = x^3+x^5+(1+\alpha(u))xQ_5+\beta(u)x^3z+\gamma(u)x^2Q_4+\delta(u)x^4y+Q_7^{(34)}+z^3,
$$
and so on ...
\bigskip

Let $\ell$ be an integer such that $\gamma_\ell \leq \gamma_\lambda$. Let
$X=X_{\mathcal{V}(\gamma_\ell)}$ be a set of independent variables, indexed by
$\mathcal{V}(\gamma_\ell)$, and consider the graded $k$-algebra
$k\left[X_{\mathcal{V}(\gamma_\ell)}\right]$, where we define
$$
\deg\ X_j=\nu(Q_j).
$$
Let $P$ denote the homogeneous monomial ideal of
$k\left[X_{\mathcal{V}(\gamma_\ell)}\right]$ generated by all  the monomials
in $X_{\mathcal{V}(\gamma_\ell)}$ of degree greater than or equal to
$\gamma_{\ell}$. We have a natural map
$$
 \begin{array}{rcl}
\phi_\ell:\frac{k\left[X_{\mathcal{V}(\gamma_\ell)}\right]}P & \rightarrow
  & \frac{\gr_\nu A}{P_{\gamma_{\ell}}} \\ X_j & \mapsto &\init_\nu
  Q_j. \end{array}
  $$

Now, for $\ell=0$, let $I_0=(0)$. For $\ell>0$, let $I_\ell$ denote the
ideal of $\frac{k\left[X_{\mathcal{V}(\gamma_\ell)}\right]}P$ generated by
all the homogeneous polynomials of the form

\begin{equation}
X^{\alpha_0}+\lambda_1X^{\alpha_1}+\lambda_2X^{\alpha_2}+\cdots+
\lambda_{b_0}X^{\alpha_{b_0}} 
\end{equation}

\noi where
$\mathbf{Q}^{\alpha_0}+\lambda_1\mathbf{Q}^{\alpha_1} +
\lambda_2\mathbf{Q}^{\alpha_2}+\cdots+
\lambda_{b_0}\mathbf{Q}^{\alpha_{b_0}}$ is the homogeneous part of
least degree of $Ex(Q)$ for an approximate root $Q \in
\mathcal{V}(\gamma_\ell)\cup\Theta(\gamma_\ell)$.

\begin{cor} \label{cor:ker}
We have $\ker\ \phi_\ell=I_\ell$.
\end{cor}

\noi Proof : The inclusion $I_\ell\subset\ker\ \phi_\ell$ is
immediate. To prove the opposite
inclusion, we argue by contradiction. Take a homogeneous
element
\begin{equation}\label{eq:hinker}
h=a_{\lambda_1}X^{\lambda_1}+a_{\lambda_2}X^{\lambda_2}+
\dots+a_{\lambda_s}X^{\lambda_s} \in \ker(\phi_\ell) \setminus I_\ell
\end{equation}
of degree $b$, $b< \gamma_\ell$, such that $\lambda_1$ is
lexicographically smallest among all the elements $h \in
\ker(\phi_\ell) \setminus I_\ell$ of degree $b$.

The inclusion (\ref{eq:hinker}) implies that
\begin{equation}
a_{\lambda_1}\init_\nu\mathbf{Q}^{\lambda_1} +
  a_{\lambda_2}\init_\nu\mathbf{Q}^{\lambda_2}+
  \dots+a_{\lambda_s}\init_\nu\mathbf{Q}^{\lambda_s}=0.
\end{equation}
in $P_b/P_{b+}$.

By definition of $I_\ell$, there exists an element $g\in I_\ell
$ of the form $X^\epsilon+\sum_pc_pX^{\epsilon_p}$ and a monomial $X^\delta$
with $\epsilon_p>\epsilon$ for all $p$ and
$\lambda_1=\epsilon+\delta$. Then, as $g \in I_\ell \subset
\ker(\phi_\ell)$, we have $h-a_{\lambda_1}X^\delta g \in
\ker(\phi_\ell)$ and the greatest monomial of
$h-a_{\lambda_1}X^\delta g$ is strictly bigger than
$X^{\lambda_1}$. This contradicts the choice of $h$.
$\Box$

\begin{cor} \label{cor:val-ideal}
Take an element $\gamma\in\Phi$, $\gamma
  < \gamma_\lambda$. The valuation ideal $P_{\gamma}$ is generated by all the
  generalized monomials of value greater than or equal to $\gamma$ in
    $\{ Q \ |\ (Q,Ex(Q)) \in \Psi(\gamma_\lambda)\}$. The ideal
    $P_{\gamma_\lambda}$ is generated by all the
  generalized monomials of value greater than or equal to $\gamma_\lambda$ in
    $\{ Q \ |\ (Q,Ex(Q)) \in \Psi(\gamma_\lambda)\cup
    \Theta(\gamma_\lambda)\}$.
\end{cor}

\noi Proof: Let $f \in P_\gamma$ (resp. $f \in
P_{\gamma_\lambda}$). By the very definition of the standard form of
level $\ell$ such that $\gamma_\ell = \gamma$, $f$ can be written as
an $A$-linear combination of generalized monomials
of value greater than or equal to $\gamma$ in $\{Q \
  |\ (Q,Ex(Q)) \in \Psi(\gamma_\lambda)\}$ (resp. $\in
  \Psi(\gamma_\lambda)\cup \Theta(\gamma_\lambda)$). Thus $P_\gamma$
  (resp. $P_{\gamma_\lambda}$) is  generated by the generalized
  monomials of value at least $\gamma$, as desired. $\Box$

\subsection{Approximate roots in a complete regular local ring}
\label{App}

We now generalize the notion of approximate root to a
\textbf{complete} regular local ring $A$ of dimension $n$, with maximal
ideal $\mathfrak{m}$, and residue field $\dsp{k=\frac{A}{\mathfrak{m}}}$. Let
  $\mathbf{u}=(u_1,\ldots,u_n)$ be a regular
system of parameters and
$$
\nu:A\setminus\{0\}\to\Gamma
$$
a valuation, centered in $\mathfrak{m}$. Denote by
$\nu_{\mathfrak{m}}$ the $\mathfrak{m}$-adic valuation.

We keep the same notation as in \S2.

The purpose of this section is to construct, for a general $\nu$, a system of
approximate roots of $\nu$, that is, a well-ordered collection of elements
$\mathbf{Q}=\{Q_i\}_{i \in \Lambda}$ of $A$ such that for every
$\nu$-ideal $I$ in $A$, we have
\begin{equation}
I = \left\lbrace\left. \prod\limits_j Q_j^{\gamma_j} \ \right|\
\sum\limits_j\gamma_j\nu(Q_j)\ge\nu(I)\right\rbrace A \label{eq:valideal}
\end{equation}
(in particular, the images $\init_\nu Q_i$ of the $Q_i$ in $\gr_\nu A$
generate $\gr_\nu A$ as a $k$-algebra). Each $Q_{i+1}$ will be described by an
explicit formula (given later in this section) in terms of the $Q_j$, $j<i$.
\medskip

In this general setting, we have to proceed by transfinite induction
on the well-ordered semigroup $\Phi$. Since we
are not assuming that $rk\ \Gamma=1$ or that $\Phi$ is Archimedean, we
have to work with ordinals other than the natural numbers.
\medskip

\noi\textit{Remark on the use of transfinite induction.} Since the
ring $A$ is noetherian, the group $\Gamma$ of values of $\nu$ has
finite rank. Therefore all the ordinals $\ell$ we will encounter in
this paper will be of type $\ell \leq \omega^n$ (cf. \cite{V5} and
\cite{CT}). Thus we will be using a very special form of transfinite induction,
which amounts to usual induction, applied finitely many times. We will,
however, stick to the language of
transfinite induction to simplify the exposition.
\medskip

Recall the definition of generalized monomial with respect to a
totally ordered set $E\subset A$ (Definition \ref{genmon}). Assume in addition
that $E$ is well-ordered. We well-order the set $\N^E$ by the lexicographical
ordering and the set of generalized monomials by the lexicographical ordering on
the set of triples
$\left(\nu\left(\mathbf{Q}^\alpha\right),\nu_{\mathfrak{m}}\left(\mathbf{Q}
^\alpha\right),\alpha\right)$.
\medskip

The semigroup $\Phi$ is well ordered. By abuse of notation, we will
sometimes write $\Phi$ for the ordinal given by the order type of
$\Phi$. Let $\lambda < \Phi$ be an ordinal and $\gamma_\lambda$ the
element of $\Phi$ corresponding to $\lambda$.
\medskip

\noi We start by choosing a coordinate system adapted to the situation.
Fix an isomorphism
\begin{equation} \label{ident2}
A \cong k[[u_1,\ldots,u_n]].
\end{equation}
\begin{deft}
Take $j \in \{2,\ldots,n\}$. We say that $u_j$ is $\nu$-prepared
if there does not exist $f \in A$ such
that $\init_\nu u_j= \init_\nu f$ and  $f \in
k[[u_1,\ldots,u_{j-1}]]$. The coordinate system $\mathbf{u} =
\{u_1,\ldots,u_n\}$ is $\nu$-prepared if $u_j$ is $\nu$-prepared for
all $j \in \{2,\ldots,n\}$.
\end{deft}

\begin{prop}
There exists a $\nu$-prepared coordinate system.
\end{prop}

\noi Proof: We construct a $\nu$-prepared coordinate system
recursively in $j$. Assume that
$u_2$, \ldots, $u_{j-1}$ are
$\nu$-prepared, but $u_j$ is not.

We will construct the prepared coordinate $\tilde{u}_j$ recursively by
transfinite induction on $\Phi$. More precisely, we will construct a
well ordered set $\{u_{ji}\}$ of successive approximation to
$\tilde{u}_j$ in the $\mathfrak{m}$-adic topology. We will show that this set
satisfies the hypothesis of Zorn's lemma and let $\tilde{u}_j$ be its
maximal element.

The details go as follows. Let $u_{j0}=u_j$. Suppose that $u_{ji}$ is
constructed and that it is not prepared. Let $f_{ji}$ be the element
$f$ of $k[[u_1,\ldots,u_{j-1}]]$ appearing in the definition of ``not
prepared''. Put $u_{j,i+1}=u_{ji}-f_{ji}$. Then
$\nu(u_{ji})=\nu(f_{ji})<\nu(u_{j,i+1})$. Next, suppose given a
sequence $u_{ji},u_{j,i+1},\ldots$ of elements of $k[[u_1,\ldots,u_
j]]$ such that $(u_1,\ldots,u_{j-1},u_{jq})$ is a regular system of
parameters of $k[[u_1,\ldots,u_j]]$ for each $q$ and
$$
\nu(u_{ji}) <\nu(u_{j,i+1}) < \nu(u_{j,i+2})< \cdots.
$$
Let $\beta_q =\nu(u_{jq})$. Since the ring $A$ is noetherian, the semi-group
$\Phi$ is well-ordered. Let
$\bar{\beta}= \min\left\{\beta \in \Phi\ \left|\ \beta >
\beta_q, \forall q \in \N\right.\right\}$. By Chevalley's lemma, applied to the
nested sequence of ideals $\displaystyle{\frac{P_{\beta_q} \cap
    k[[u_1, \ldots,u_{j-1}]]}{P_{\bar{\beta}}\cap k[[u_1,
    \ldots,u_{j-1}]] }}$ in the complete local ring
$\displaystyle{\frac{k[[u_1, \ldots,u_{j-1}]]}{P_{\bar{\beta}}\cap k[[u_1,
    \ldots,u_{j-1}]] }}$, we see that $\lim\limits_{q \to
  \infty}(f_{jq} \mod P_{\bar{\beta}}) =0$
in the $(u_1,\ldots,u_{j-1})$-adic topology.

Hence, modifying each $f_{jq}$ by an element of $P_{\bar{\beta}}$ if
necessary, we may assume that
$$
\lim\limits_{q \to \infty} f_{jq}=0.
$$
We define $u_{j,i+\omega}$ to be the formal power series
$u_{ji}-f_{ji}-f_{j,i+1}- \cdots$. By construction,
$$
\nu(u_{j,i+\omega}) \geq \bar{\beta}.
$$
To complete our construction , we need to consider countable well ordered sets
$\{u_{jt}\}$ of order type greater than $\omega$. This presents no problem: by
countability, we can always choose a cofinal subsequence in each such set. Then
the above construction of $u_{j,i+\omega}$ applies verbatim.
$\Box$
\medskip

\noi We construct, inductively in $\lambda$, two well-ordered
sets $\Lambda(\gamma_\lambda)$ and $\Theta(\gamma_\lambda)$ and, in the case
$\lambda$ is not a limit ordinal, a well ordering of the set
$\Lambda(\gamma_\lambda) \cup \Theta(\gamma_{\lambda-1})$, compatible with the
orders on $\Lambda(\gamma_\lambda)$ and $\Theta(\gamma_{\lambda-1})$.  At each
step we define two additional well-ordered sets
$\mathcal{V}(\gamma_\lambda) \subset \Psi(\gamma_\lambda) \subset
\Lambda(\gamma_\lambda)$ where the inclusions are inclusions of ordered
sets. Both collections of sets $\Lambda(\gamma_\lambda)$ and
$\Nu(\gamma_\lambda)$ will be increasing with $\lambda$.

A typical element of each of
those sets will have the form $(Q,\mbox{Ex}(Q))$ where $Q \in A$ and
$\mbox{Ex}(Q)$ is an increasing sum of monomials in
$\mathcal{V}(\gamma_\lambda) \cup \Theta(\gamma_{\lambda-1})$ if $\lambda$ is
not a limit ordinal, resp. monomials in $\Nu(\gamma_\lambda)$ if $\lambda$ is a
limit ordinal. The sum in Ex$(Q)$
may be finite or infinite, but it is always convergent in the
$\mathfrak{m}$-adic topology. Given an element
$(Q,\mbox{Ex}(Q)) \in \Lambda(\gamma_\Phi) \cup \Theta(\gamma_\Phi)$,
$Q$ is called an \textit{approximate root} and $\mbox{Ex}(Q)$ is
called the \textit{expression} of $Q$.
\medskip

For an ordinal $\ell < \Phi$ and for $(Q,\mbox{Ex}(Q)) \in
\Lambda(\gamma_\ell) \cup \Theta(\gamma_\ell)$, let
 $\mbox{In}\ Q$ denote the smallest monomial of $\mbox{Ex}(Q)$. Let
 $\mbox{In}(\ell) = \left\{\left.\alpha \in \N^{\Nu(\gamma_\ell)} \ \right|\
\exists (Q,\mbox{Ex}(Q))
 \in \Lambda(\gamma_\ell) \mbox{ such that } \mathbf{Q}^\alpha =
 \mbox{In}\ Q\right\}$.

\begin{thm}\label{thm-aproot}
For $\lambda < \Phi$, there exist well ordered sets
$\Nu(\gamma_\lambda) \subset \Psi(\gamma_\lambda) \subset
\Lambda(\gamma_\lambda)$ and
$\Theta(\gamma_\lambda)$, and a well ordering of $\Lambda(\gamma_\lambda) \cup
\Theta(\gamma_{\lambda-1})$ when $\lambda$ is not a limit ordinal, having
the following properties. Let
\begin{eqnarray}
\Psi(<\gamma_\lambda) &=&
\Psi(\gamma_{\lambda-1})\quad\text{ if
$\lambda$ is not a limit ordinal and }\\
\Psi(<\gamma_\lambda) &=&\Psi(\gamma_\lambda)\qquad\text{otherwise}\label{eq:Psi<lambda}
\end{eqnarray}
and similarly for $\Nu(<\gamma_\lambda)$. Then each set
$\Nu(\gamma_\lambda),\Psi(\gamma_\lambda),\ \Lambda(\gamma_\lambda),\
\Theta(\gamma_\lambda)$ consists of elements of the form
$(Q,\mbox{Ex}(Q))$, with $Q \in A$ and $\mbox{Ex}(Q)$ is an
increasing (with respect to the monomial order defined above) sum of monomials in
$\Nu(<\gamma_\lambda) \cup \Theta(\gamma_{\lambda-1})$ when $\lambda$ is not a
limit ordinal, resp. $\Nu(\gamma_\lambda)$ when $\lambda$ is a limit ordinal, of
value $< \nu(Q)$, provided $Q \notin \{u_1,\ldots,u_n\}$,
such that
\begin{eqnarray}
 \nu(Q) < \gamma_\lambda \mbox{ whenever } (Q,\mbox{Ex}(Q)) \in
 \Lambda(\gamma_\lambda) \label{thm1}\\
  \nu(Q) \geq \gamma_\lambda \mbox{ whenever } (Q,\mbox{Ex}(Q)) \in
 \Theta(\gamma_\lambda) \label{thm2}
\end{eqnarray}
and the sets

\begin{equation}\label{thm3}
\{(Q,\mbox{Ex}(Q)) \in \Theta(\gamma_\lambda) \cup \Lambda(\gamma_\lambda)
\ |\ \nu(Q) = \gamma\}, \  \gamma \in \Phi
\end{equation} and
\begin{equation} \label{thm4}
\left\lbrace \left. (Q,\mbox{Ex}(Q))
\in \Psi(\gamma_\lambda) \cup \Theta(\gamma_\lambda) \ \right|\ Q
\notin \mathfrak{m}^s
\right\rbrace, \ s \in \N
\end{equation} are finite.
An element $(Q,\mbox{Ex}(Q)) \in
\Psi(\gamma_\lambda) \cup  \Theta(\gamma_\lambda)$ is completely determined
by $\mbox{In}\ Q$; moreover
$\nu_{\mathfrak{m}}(\mbox{In}\ Q)=\nu_{\mathfrak{m}}(Q)$.
\end{thm}

\noi In what follows, $\Lambda(<\gamma_\lambda)$
will stand for $\bigcup\limits_{\ell <  \lambda}\Lambda(\gamma_\ell)$.

\noi Proof : We proceed by transfinite induction.

First define $\Psi(\mathbf{1}) = \Lambda(\mathbf{1}) =
\emptyset$ and $\Theta(\mathbf{1})  = \{(u_1,u_1),\ldots,(u_n,u_n)\}$
where we
assume $$\nu(u_1) \leq \nu(u_2) \leq \cdots \leq \nu(u_n).$$

\noi We define the well ordering on $\Theta(\mathbf{1})$ by $(u_1,u_1) <
(u_2,u_2) < \cdots < (u_n,u_n)$.
\medskip

Let $\lambda < \Phi$ be an ordinal. Assume that for each $\ell <
\lambda$ we have constructed sets $\Psi(\gamma_\ell) \subset
\Lambda(\gamma_\ell)$ and $\Theta(\gamma_\ell)$ and a well ordering of
$\Lambda(\gamma_\ell) \cup \Theta(\gamma_{\ell-1})$, having the
properties required in the theorem.

Let
\begin{eqnarray} \Lambda(\gamma_\lambda)&=&
  \Lambda(<\gamma_\lambda) \text{ if } \lambda \text{ is a limit
    ordinal }\label{eq:Lambda<lambda} \\
\Lambda(\gamma_\lambda)&=& \Lambda(\gamma_{\lambda-1}) \cup
\{(Q,\mbox{Ex}(Q)) \in \Theta(\gamma_{\lambda-1}) \ |\ \nu(Q) <
\gamma_\lambda\} \text{ otherwise.}\label{eq:Lambda<lambdanotlimit}
\end{eqnarray}

\begin{deft}\label{pred-ines} An element $(Q,\mbox{Ex}(Q)) \in
  \Lambda(\gamma_\lambda)$ is an \textbf{ inessential predecessor} of
  a root
  $(Q',\mbox{Ex}(Q')) \in \Lambda(\gamma_\lambda)$ if $\mbox{Ex}(Q') =
  \mbox{Ex}(Q) + \sum_\alpha c_\alpha \mathbf{Q}^\alpha$, where each
  $c_\alpha$ is a unit in $A$ and $\mathbf{Q}^\alpha$ a monomial in
  $\Nu(\gamma_\lambda)$.
\smallskip

An element $(Q,\mbox{Ex}(Q)) \in \Lambda(\gamma_\lambda)$ is said to be
\textbf{essential at the level $\gamma_\lambda$} if $Q$ is not an inessential
predecessor of an element of $\Lambda(\gamma_\lambda)$.
\end{deft}


Let $\Psi(\gamma_\lambda)$ be the subset of $\Lambda(\gamma_\lambda)$
consisting of all the essential roots at the level $\gamma_\lambda$.
Let $\mathcal{V}(\gamma_\lambda)$ be the
subset of $\Psi(\gamma_\lambda)$ consisting of all $(Q,\mbox{Ex}(Q))$ such
that $\init_\nu(Q)$ does not belong to the $k$-vector space of
$\gr_\nu(A)$ generated by the set $\{\init_\nu \mathbf{Q}^\gamma \}$
where $\mathbf{Q}^\gamma$ runs over the set of all the generalized
monomials on roots preceding $Q$ in the above ordering.

\noi We extend the well ordering from $\Lambda(<\gamma_\lambda)$ to
$\Lambda(\gamma_\lambda)$ by postulating that $\Lambda(<\gamma_\lambda)$ is the
initial segment of $\Lambda(\gamma_\lambda)$.   Moreover, we extend this well
ordering from $\Lambda(\gamma_\lambda)$ to $\Lambda(\gamma_\lambda) \cup
\Theta(\gamma_{\lambda-1})$.
\medskip

If $\ell$ is not a limit ordinal, let
$E(\ell) =
\mbox{In}(\ell)+ \N^{\Nu(\gamma_\ell)} \subset
\N^{\Nu(\gamma_\ell)}$. Now, if $\ell' < \ell''$, we have
$\Nu(\gamma_{\ell'}) \subset \Nu(\lambda_{\ell''})$, which induces an
inclusion
$\N^{\Nu(\gamma_{\ell'})} \subset \N^{\Nu(\gamma_{\ell''})}$. If
$\ell$ is a limit ordinal, define $E(\ell) = \bigcup\limits_{\ell' <
  \ell}E(\ell')$.
\medskip

\noi \textbf{Notation.} Denote by $\Theta(<\gamma_\lambda)$ the set
$\dsp{\bigcup_{\ell < \lambda}\Theta(\gamma_\ell) \setminus
\Lambda(<\gamma_\lambda)}$.

\begin{rek}\label{rek-Psi<gammalambda} We have
\begin{equation}
\Psi(\gamma_\lambda)\cup\Theta(<\gamma_\lambda)=\Psi(<\gamma_\lambda)\cup\Theta(
<\gamma_\lambda) .\label{eq:Psi<gammalambda}
\end{equation}
Indeed, consider an element
$(Q,\mbox{Ex}(Q))\in\Psi(\gamma_\lambda)\cup\Theta(<\gamma_\lambda)$. If
$\lambda$  is a limit ordinal, then
\begin{equation} \label{eq:indhypapplicable}
(Q,\mbox{Ex}(Q))\in\Psi(<\gamma_\lambda)\cup\Theta(<\gamma_\lambda)
\end{equation}
by (\ref{eq:Psi<lambda}). If $\lambda$ is not a limit ordinal and $(Q,\mbox{Ex}(Q))\in\Psi(\gamma_\lambda)\setminus\Psi(\gamma_{\lambda-1})$ then
$$
(Q,\mbox{Ex}(Q))\in\Theta(\gamma_{\lambda-1})
$$
by (\ref{eq:Lambda<lambdanotlimit}). Thus (\ref{eq:indhypapplicable}) holds in
all the cases and  (\ref{eq:Psi<gammalambda}) is proved.
\end{rek}

\begin{lem} \label{fini1} The set $\mathbf{Q}(h) = \left\{\left.(Q,\mbox{Ex}(Q))
  \in \Psi(\gamma_\lambda) \cup \Theta(<\gamma_\lambda)\ \right|\ Q
  \notin \mathfrak{m}^h\right\}$ is finite for every $h \in \N$.
\end{lem}

\noi Proof: Consider an element $(Q,\mbox{Ex}(Q))\in
\mathbf{Q}(h)$. If $(Q,Ex(Q)) \in \Theta(<\gamma_\lambda)$, then there
exists
$\ell < \lambda$ such that  $(Q,Ex(Q)) \in
\Theta(<\gamma_\ell)$. If $(Q,Ex(Q)) \in \Psi(<\gamma_\lambda) \subset
\Lambda(\gamma_\lambda)= \bigcup_{\ell <
  \lambda}\Lambda(\gamma_\ell)$, then there exists
$\ell < \lambda$ such that $(Q,Ex(Q)) \in \Lambda(\gamma_\ell)$. Since
$Q$ is essential at level $\gamma_\lambda$, it is also essential at
level $\gamma_\ell$, so $(Q,Ex(Q)) \in \Psi(\gamma_\ell)$. Thus by the
induction hypothesis on $\lambda$, for any $Q \in\mathbf{Q}(h)$, we
have $\nu_{\mathfrak{m}}(Q) = \nu_{\mathfrak{m}}(\mbox{In}\ Q)$.

Write $\mbox{Ex}(Q) = \mathbf{Q}^{\alpha_0} + \cdots $ where, by construction,
$\mathbf{Q}^{\alpha_0}$ is either a $u_r$ or a product of at least 2 terms,
$\mathbf{Q}^{\alpha_0} = \prod Q_s^{\beta_s}$.

\noi In the first case, the number of such $\mathbf{Q}^{\alpha_0}$ is
finite, because the number of $u_k$ is finite.

\noi In the second case,
$\nu_{\mathfrak{m}}(Q_s) < \nu_{\mathfrak{m}}(\mathbf{Q}^{\alpha_0}) \leq
\nu_{\mathfrak{m}}(Q) < h$. So $\nu_{\mathfrak{m}}(Q_s) < h-1$ and, by
induction on $h$, the number of such $Q_s$ is finite. If
$$
m = \min \left\{\nu_{\mathfrak{m}}(Q_s)\ \left|\
Q_s \text{ divides } \mathbf{Q}^{\alpha_0}\right.\right\},
$$
then $|\alpha_0|m \leq\nu_{\mathfrak{m}}(\mathbf{Q}^{\alpha_0})\le h-1$, so
there is a finite number of such $\alpha_0$ possible which means that the
number of such $\mathbf{Q}^{\alpha_0}$ is finite. By the induction hypothesis,
$Q$ is completely determined by $\mbox{In}\ Q$ whenever $(Q,\mbox{Ex}(Q)) \in
\Psi(\gamma_\lambda) \cup \Theta(<\gamma_\lambda)$. Therefore $\mathbf{Q}(h)$
is finite. $\Box$
\medskip

\begin{cor} The set of monomials $\left\{\mathbf{Q}^\alpha\ \left|\
\mathbf{Q}^\alpha \notin \mathfrak{m}^s\right.\right\}$ in
$\Psi(\gamma_\lambda) \cup \Theta(<\gamma_\lambda)$ is finite for every $s \in
\N$.
\end{cor}

\begin{cor} \label{cor:infinite}
(1) Any infinite sequence of generalized monomials in $\Psi(\gamma_\lambda)
\cup \Theta(<\gamma_\lambda)$, all of whose members are distinct,
converges to 0 in the $\mathbf{m}$-adic topology.

(2) Any infinite series, all of whose terms are distinct generalized
monomials in $\Psi(\gamma_\lambda) \cup \Theta(<\gamma_\lambda)$
converges in the $\mathbf{m}$-adic topology.
\end{cor}

\begin{lem} \label{AR}
The set
$$
\mathbf{Q}^\alpha = \prod Q^{\alpha_Q} \mbox{ such that }  (Q,\mbox{Ex}(Q))\in
\Psi(\gamma_\lambda) \cup \{(Q,,Ex(Q)) \in \Theta(<\gamma_\lambda) \ |\
\nu(Q) = \gamma_\lambda \}$$ and $\nu(\mathbf{Q}^\alpha) =
\gamma_\lambda$ is finite.
\end{lem}

\noi Proof: By the Artin-Rees lemma, there exists $p_0$ such that,
for $p \geq p_0$,
$$
\mathfrak{m}^p \cap P_{\gamma_\lambda} =
\mathfrak{m}^{p-p_0}(\mathfrak{m}^{p_0}\cap P_{\gamma_\lambda}).$$
Take $p > p_0$, then
\begin{equation}
\mathfrak{m}^p \cap P_{\gamma_\lambda} \subset \mathfrak{m}
P_{\gamma_\lambda} \subset P_{\gamma_\lambda +}.
\end{equation}
This equation shows that the set of the lemma is disjoint from
$\mathfrak{m}^p$. So by the above corollary, the set of the lemma is
finite.
$\Box$
\medskip

Consider now the ordered set
$\{\mathbf{Q}^{\alpha_1},\ldots,\mathbf{Q}^{\alpha_s} \}$ of
monomials \begin{equation}
    \label{ens-mon}
\mathbf{Q}^\alpha = \prod Q^{\alpha_Q},\ (Q,\mbox{Ex}(Q)) \in
\Nu(\gamma_\lambda)
\cup \{(Q,\mbox{Ex}(Q)) \in \Theta(<\gamma_\lambda) \ |\ \nu(Q) =
\gamma_\lambda  \}
          \end{equation}
of value $\gamma_\lambda$ such that the natural projection of $\alpha$
to $\N^{\mathcal{V}(\gamma_\lambda)}$ does not belong to
$E(\lambda)$. The fact that this set is \textbf{finite} follows from
the above Lemma and the fact that $\Nu(\gamma_\lambda) \subset
\Psi(\gamma_\lambda)$.
\smallskip

Let $i_1 = \max \left\{ i \in \{ 1,\ldots,s \}\ \left|\
\init_\nu\left(\mathbf{Q}^{\alpha_i}\right)  \in
\sum\limits_{j=i+1}^s k \
\init_\nu\left(\mathbf{Q}^{\alpha_j}\right)\right.\right\}$ and
consider the unique relation  $\init_\nu\left(\mathbf{Q}^{\alpha_{i_1}}\right) -
\sum\limits_{j={i_1}+1}^sc_{1j} \
\init_\nu\left(\mathbf{Q}^{\alpha_j}\right)=0$.  Let $P_1 =
  \mathbf{Q}^{\alpha_{i_1}} - \sum\limits_{j={i_1}+1}^s
c_{1j} \mathbf{Q}^{\alpha_j}$ where we view $k$ as a subring of $A$
via the identification (\ref{ident2}).

Let $i_2 = \max\left\{i \in \{1,\ldots,i_1-1\}\ \left|\
\init_\nu\left(\mathbf{Q}^{\alpha_i}\right) \in \sum\limits_{j=i+1}^s k\
\init_\nu\left(\mathbf{Q}^{\alpha_j}\right)\right.\right\}$ and, as
before, consider the unique $\displaystyle{P_2 = \mathbf{Q}^{\alpha_{i_2}} -
  \sum_{\stackrel{j={i_2}+1}{j \neq i_1}}^s c_{2j}
\mathbf{Q}^{\alpha_j}}$ such that the vector
$(\alpha_j)_{j=i_1+1,\ldots,s},\ c_{2j} \neq
  0$, is minimal in the lexicographical order and define so on
  uniquely $P_3,\ldots, P_t$.
\smallskip

Now, if $\lambda$ has a predecessor, we let
\begin{equation}\Theta(\gamma_\lambda) =
\left\{\left.(Q,\mbox{Ex}(Q)) \in \Theta(<\gamma_\lambda) \ \right|\
  \nu(Q) \geq \gamma_\lambda\right\}
\cup \{(P_1,\mbox{Ex}(P_1)),\ldots,(P_t,\mbox{Ex}(P_t))\}
\end{equation} where
\begin{equation}\label{formul1} \mbox{Ex}(P_j) =
\mathbf{Q}^{\alpha_{i_j}} - \sum_k c_{jk}
\mathbf{Q}^{\alpha_k}
\end{equation}
if $\mathbf{Q}^{\alpha_{i_j}}$ is not a
preceding root $Q$ and
\begin{equation} \label{formul2}
\mbox{Ex}(P_j) = \mbox{Ex}(Q) - \sum_k c_{jk} \mathbf{Q}^{\alpha_k}
\end{equation} in the other case. We define the order on
$\Theta(\gamma_\lambda)$ by
$\Theta(\gamma_{\lambda-1}) <
\left\{(P_1,\mbox{Ex}(P_1)),\ldots,(P_t,\mbox{Ex}(P_t))\right\}$ and
$(P_1,\mbox{Ex}(P_1)) < \cdots < (P_t,\mbox{Ex}(P_t))$.
\begin{rek}
Note that, because the system of coordinates is prepared,
$u_1,\ldots,u_n$ are always essential.
\end{rek}

\begin{rek}\label{rk-continubis} Note that Remark \ref{rk-continu}
  remains valid in this context, with the obvious modification that
  the expressions of approximate roots are now allowed to be infinite,
  but convergent in the $\mathfrak m$-adic topology.
\end{rek}

Suppose now $\lambda$ is a limit ordinal. Let $(Q_0,Ex(Q_0)) \in
  \Lambda(\gamma_{\ell_0})$ for some $\ell_0< \lambda$ and $\mathbf{Q}^\alpha =
  In(Q_0)$. Let $L(Q_0)$ be the following infinite well ordered set
  of approximate roots, indexed by ordinals $\ell$,
$\ell_0 \leq \ell < \lambda$
$$L(Q_0)=\{\ (Q^{(\ell)},Ex(Q^{(\ell)})) \in
  \Psi(\gamma_\ell)\ \}_{\ell_0 \leq \ell < \lambda}$$ such that
  $\mbox{In} Q^{(\ell)} = \mathbf{Q}^\alpha$.

By Remarks \ref{rk-continu} and \ref{rk-continubis}, for
$\ell_0 \leq \ell < \ell' < \lambda$, we have
\begin{equation}\label{exquell}
\mbox{Ex}(Q^{(\ell')})= \mbox{Ex}(Q^{(\ell)}) + \sum_{j \in W}c_j
\mathbf{Q}^{\alpha_j}
\end{equation}
where $\nu(\mathbf{Q}^{\alpha_j}) \geq \nu(Q^{(\ell)})$.

Let $p$ be a positive integer. By induction assumption, all the approximate
roots $Q$ appearing in any of the monomials $\mathbf{Q}^{\alpha_j}$
belong to $\Nu(\gamma_\lambda)$ and, by lemma \ref{fini1}, the number
of such roots outside $\mathfrak{m}^p$
is finite. Thus, all but finitely many $\mathbf{Q}^{\alpha_j}$ belong
to $\mathfrak{m}^p$. This proves that $L(Q_0)$ has a
limit in $A$ with respect to the $\mathfrak{m}$-adic topology :
  $(\lim\limits_\to Q,\lim\limits_\to\mbox{Ex}(Q))$.

Let
\begin{equation}\Theta(\gamma_\lambda) = \{(Q,\mbox{Ex}(Q)) \in
  \Theta(<\gamma_\lambda)
  \ |\ \nu(Q) \geq \gamma_\lambda\} \cup \hat{L}
\end{equation} where $\hat{L}$ consists of  all
couples of the form $(\lim\limits_\to Q,\lim\limits_\to\mbox{Ex}(Q))$.

So finally, the expression of an approximate root has the form
\begin{equation}\label{eq:apprroot}
\mbox{Ex}(Q) = \mathbf{Q}^\alpha + \sum_k a_k
\mathbf{Q}^{\alpha_k}
\end{equation}
the sum, written in the increasing order of the monomials, being
finite or infinite.
\medskip

We now prove the finiteness of sets (\ref{thm3}) and
(\ref{thm4}). First, note that  the set
\begin{equation}
\{(Q,\mbox{Ex}(Q)) \in \Theta(<\gamma_\lambda) \cup \Lambda(\gamma_\lambda)
\ |\ \nu(Q) = \gamma\}, \  \gamma \in \Phi
\end{equation}
is finite by the induction hypothesis and the set
\begin{equation}
\left\lbrace \left. (Q,\mbox{Ex}(Q))
\in \Psi(\gamma_\lambda) \cup \Theta(<\gamma_\lambda) \ \right|\ Q
\notin \mathfrak{m}^p \right\rbrace, \ p \in \N
\end{equation}
is finite by the induction hypothesis and lemma
(\ref{AR}). If $\lambda$ is not a limit ordinal,
the finiteness of (\ref{thm3}) and (\ref{thm4}) follows from the fact
that the set $\Theta(\gamma_\lambda)\setminus\Theta(<\gamma_\lambda)$
is finite by construction. If $\lambda$ is a limit ordinal, to prove
finiteness of (\ref{thm3}) and (\ref{thm4}), it remains to prove that
the set
\begin{equation}
    \{ (Q,\mbox{Ex}(Q))
\in \hat{L} \ |\ \nu(Q) = \gamma \}
\end{equation} is finite. This is proved in exactly the same way as
lemma (\ref{AR}). This completes the proof of the finiteness of
(\ref{thm3}) and (\ref{thm4}).
\medskip

The property that the monomials appearing in $\mbox{Ex}(Q)$ are
arranged in increasing order with respect to the $\nu$-adic value
holds for all the newly constructed approximate roots. Next we show
that $\nu_{\mathfrak{m}}(\mbox{In}Q) = \nu_{\mathfrak{m}}(Q)$ for all
those new approximate
roots. Indeed, if $\lambda$ is not a limit ordinal and $\mbox{Ex}(Q)$
is given by formula (\ref{formul1}), all the monomials appearing in
$\mbox{Ex}(Q)$ have the same $\nu$-adic value and their
$\nu_{\mathfrak{m}}$-adic
values are increasing because of the order we imposed on monomials
which proves that $\nu_{\mathfrak{m}}(\mbox{In}Q) =
\nu_{\mathfrak{m}}(Q)$. If $(Q',\mbox{Ex}(Q'))$ is an approximate root
whose expression is given by formula (\ref{formul2}), with $P_j$
playing the role of $Q'$, let $\mathbf{Q}^{\alpha_0}= \mbox{In
}Q$. We have $Q' = Q + \sum c_\alpha \mathbf{Q}^\alpha$, where
$\nu(\mathbf{Q}^\alpha) = \nu(Q)$. Then $\nu_{\mathfrak{m}}(Q) \leq
\nu_{\mathfrak{m}}(\mathbf{Q}^\alpha)$ for all $\alpha$, because of
the order on monomials. So that finally,
$\nu_{\mathfrak{m}}(\mathbf{Q}^{\alpha_0}) \leq \nu_{\mathfrak{m}}(Q) \leq
\nu_{\mathfrak{m}}(\mathbf{Q}^\alpha)$, which proves that
$\nu_{\mathfrak{m}}(\mbox{In}Q') = \nu_{\mathfrak{m}}(Q')$. The property
that $\nu_{\mathfrak{m}}(\mbox{In}Q) = \nu_{\mathfrak{m}}(Q)$ is
clearly preserved by passing to the limit, so it also
holds in the case when $\lambda$ is a limit ordinal.

\begin{rek} We just showed that there is a one to one
   correspondence between the approximate roots $Q \in
   \Psi(\gamma_\ell)$ and the set of
   monomials which are the first term of the expression $\mbox{Ex}(Q)$
   of such an approximate root $Q$. Let us denote by
   $\mathbf{M}(\ell)$ \label{em-de-ell} the set of those monomials.
\end{rek}

We well order $\hat{L}$ by the lexicographical order of the triples
$(\nu(Q),\nu_{\mathfrak{m}}(Q),\text{In}(Q))$,$Q \in \hat{L}$. We extend this
ordering to $\Theta(\gamma_\lambda)$ by postulating that $\hat{L}$ is the final
segment in $\Theta(\gamma_\lambda)$.

The rest of Theorem \ref{thm-aproot} holds by construction.

\subsection{Standard form in the case of complete regular local
rings}\label{standardform}

Let $\Psi(\gamma_\Phi) = \bigcup_{\ell <\Phi} \bigcap_{\ell \leq \ell' <
  \Phi} \Psi(\gamma_{\ell'})$ and let $\mathcal{V}(\gamma_\Phi)$ be the set of
approximate roots, essential at the level $\gamma_\Phi$.

In this section, we fix an ordinal $\lambda \leq \Phi$.

\begin{deft}
A monomial in $\Psi(\gamma_\lambda)\cup \Theta(\gamma_\lambda)$ is
called standard with
respect to $\lambda$ if all the approximate roots appearing in it belong to
$\mathcal{V}(\gamma_\lambda)$ and it is not
divisible by any $In Q$ where $Q$ is an approximate root in
$(\Psi(\gamma_\lambda)\cup \Theta(\gamma_\lambda)) \setminus
\{(u_1,u_1),\ldots,(u_n,u_n)\}$.
\end{deft}
Take an ordinal $\ell\le\lambda$.
\begin{deft} Let $f \in A$. An expansion of $f$ of the form $f = \sum
  c_\alpha \mathbf{Q}^\alpha$ where the $\mathbf{Q}^\alpha$ are
  monomials in $\Psi(\gamma_\lambda)\cup \Theta(\gamma_\lambda)$,
  written in increasing order, is a standard form of level
  $\gamma_\ell$ if $\forall \gamma' < \gamma_\ell$ and for all
  $\alpha$ such that $\nu(\mathbf{Q}^\alpha)= \gamma'$,
$\mathbf{Q}^\alpha$ is a standard monomial.
\end{deft}

We now construct by induction on $\ell$ a standard form of $f$ of
level $\gamma_\ell$. We will write this standard form as $$f = f_\ell + \sum
c_\alpha \mathbf{Q}^\alpha$$ where, for all $\alpha$,
$\mathbf{Q}^\alpha$ is a generalized monomial in
$\Psi(\gamma_\lambda)\cup \Theta(\gamma_\lambda)$,
$\nu(\mathbf{Q}^\alpha)\ge\gamma_\ell$ and $f_\ell$ is a sum of standard
monomials in $\Nu(\gamma_\lambda)$ of value strictly less than $ \gamma_\ell$.
\medskip

To start the induction, let $f_0=0$. The \textit{standard form
of $f$ of level 0} will be its expansion, $f = f_0 + \sum
c_\alpha \mathbf{u}^\alpha$, written in increasing order according to the
monomial order defined above, as a formal power series in the $u_i$.
\medskip

Let $\ell < \lambda$ be an ordinal. Let us define $f_{\ell +1}$ and the
standard form of $f$ of level $\gamma_{\ell +1}$ as follows. Assume,
inductively, that a standard form of level  $\gamma_\ell$ is already defined: $f
= f_\ell + \sum c_\alpha \mathbf{Q}^\alpha$ with
$\nu(\mathbf{Q}^\alpha) \geq
\gamma_\ell$, for all $\alpha$, and the value of any monomial of
$f_\ell$ is strictly less than $\gamma_\ell$.

Take the homogeneous part of $\sum c_\alpha \mathbf{Q}^\alpha$ of
value $\gamma_\ell$, the monomials being
written in increasing order. Assume that not all the $\mathbf{Q}^\alpha$ are
standard with respect to $\lambda$, and take the smallest non standard
$\mathbf{Q}^\alpha$. Since $\mathbf{Q}^\alpha$ is not standard, one of the two
following conditions holds:
\begin{enumerate}
\item  There exists an approximate
root $Q\in\left(\Psi(\gamma_\lambda)\cup \Theta(\gamma_\lambda)\right)
\setminus\{(u_1,u_1),\dots,(u_n,u_n)\}$ such that $In(Q)$ divides
$\mathbf{Q}^\alpha$. Write $Q=In(Q)+\sum  c_\beta \mathbf{Q}^\beta$ and replace
$In(Q)$ by $Q-\sum c_\beta \mathbf{Q}^\beta$ in
$\mathbf{Q}^\alpha$.
\item An approximate root $Q \in \Psi(\gamma_\lambda)\setminus
\Nu(\gamma_\lambda)$ divides $\mathbf{Q}^\alpha$. Since $Q \notin
\Nu(\gamma_\lambda)$, there exists
$$
Q' \in\Psi(\gamma_\lambda)\cup\Theta(\gamma_\lambda)
$$
of the form $Q'=Q+ \sum_\beta d_\beta \mathbf{Q}^\beta$,
where the $\mathbf{Q}^\beta$ are monomials in
$\Nu(\gamma_\lambda)$ of value greater than or equal to
$\gamma_\ell$. Replace $Q$ by $Q'-\sum_\beta d_\beta \mathbf{Q}^\beta$.
\end{enumerate}
In both cases, those changes introduce new monomials, with increasing
$\nu_{\mathfrak{m}}$ value, but either they are of
value strictly greater than $\gamma_\ell$ or they are of value
exactly $\gamma_\ell$ but greater than $\mathbf{Q}^\alpha$ in the
monomial ordering. We repeat this procedure as many times as we can.
After a finite number of steps, no more changes are available involving
monomials of value exactly $\gamma_\ell$. Then, let $f_{\ell+1} = f_\ell +
\sum d_\rho
\mathbf{Q}^\rho$ with $\nu(\mathbf{Q}^\rho) = \gamma_\ell$, so
that $f = f_{\ell+1} + \sum c_\alpha \mathbf{Q}^\alpha$ where
$\nu(\mathbf{Q}^\alpha) > \gamma_\ell$.
\medskip

\noi Suppose now that $\mu$ is a limit ordinal. For each $\ell < \mu$,
write $f=f_\ell+ \delta_\ell$ where $f_\ell$ is a sum of standard
monomials, with respect to $\lambda$, of value strictly less than
$\gamma_\ell$ and
$\delta_\ell$ is a sum of monomials in $\Psi(\gamma_\lambda) \cup
\Theta(\gamma_\lambda)$, of value greater than or equal to $ \gamma_\ell$. We
assume inductively that, for each $\ell < \mu$ and for each generalized
monomial $\mathbf{Q}^\tau$ in $\Psi(\gamma_\lambda) \cup
\Theta(\gamma_\lambda)$, there exist
$c_\tau, b_\tau \in k$ and an ordinal $\ell_0 < \ell$ such that, for
all $\ell'$, $\ell_0 < \ell' < \ell$, the monomial $\mathbf{Q}^\tau$
appears in $f_{\ell'}$ with coefficient $c_\tau$ and in $\delta_{\ell'}$
with coefficient $b_\tau$. Moreover, assume that
$f_\ell = \lim\limits_{\stackrel{\to}{\ell' < \ell}}f_{\ell'}
= \sum_\tau c_\tau \mathbf{Q}^\tau$ and
$\delta_{\ell} = \lim\limits_{\stackrel{\to}{\ell' <
    \ell}}\delta_{\ell'} = \sum_\tau b_\tau \mathbf{Q}^\tau$.

\begin{lem} \label{lem:conv}
Consider a generalized monomial $\mathbf{Q}^\tau$ in
$\Psi(\gamma_\lambda) \cup \Theta(\gamma_\lambda)$. There exist
$c_\tau, b_\tau \in k$ and an ordinal $\ell_0 < \mu$ such that, for
all $\ell$, $\ell_0 < \ell < \mu$, the monomial $\mathbf{Q}^\tau$
appears in $f_\ell$ with coefficient $c_\tau$ and in $\delta_\ell$
with coefficient $b_\tau$.
\end{lem}

\begin{cor} \label{cor-limitexist}
The limits $\lim\limits_{\stackrel{\to}{\ell < \mu}} f_\ell$ and
$\lim\limits_{\stackrel{\to}{\ell < \mu}} \delta_\ell$ exist in the
$\mathfrak{m}$-adic topology.
\end{cor}

\noi Proof of Corollary \ref{cor-limitexist} : This is an immediate
consequence of the Lemma and Corollary \ref{cor:infinite}. $\Box$
\medskip

\noi Proof of Lemma \ref{lem:conv}: The existence of $c_\tau$ in the
lemma follows immediately from the construction and the induction
hypothesis.

If $\nu(\mathbf{Q}^\tau) < \gamma_\mu$, put $b_\tau = 0$. Assume
$\nu(\mathbf{Q}^\tau) \geq \gamma_\mu$.
For $\ell < \mu$, let $b_\tau(\ell)$ denote the coefficient of
$\mathbf{Q}^\tau$ in $\delta_\ell$.
Take an ordinal $\ell < \mu$. Suppose
\begin{equation} \label{ineq:btau}
b_\tau(\ell) \neq b_\tau(\ell+1).
\end{equation}
This means that in the above construction of
$f_{\ell+1} + \delta_{\ell+1}$ from $f_\ell+\delta_\ell$,
$\mathbf{Q}^\tau$ appears in one of the expressions
$\dsp{\frac{\mathbf{Q}^\alpha}{In
    Q}Q}$,$\dsp{\frac{\mathbf{Q}^\alpha}{In Q}\sum_\beta d_\beta
  \mathbf{Q}^\beta}$ (case 1 of the construction) or
$\dsp{\frac{\mathbf{Q}^\alpha}{Q}Q'}$,
$\dsp{\frac{\mathbf{Q}^\alpha}{Q} \sum_\beta d_\beta\mathbf{Q}^\beta}$
(case 2 of the construction). Then
\begin{equation} \label{ineq:num}
\nu_{\mathfrak{m}}(\mathbf{Q}^\alpha) \leq
\nu_{\mathfrak{m}}(\mathbf{Q}^\tau).
\end{equation}
Suppose that there were infinitely many $\ell$ for which (\ref{ineq:btau})
holds. This would mean that there are
infinitely many monomials
$\mathbf{Q}^\alpha$ (all distinct because
$\nu(\mathbf{Q}^\alpha)=\gamma_\ell$), satisfying
(\ref{ineq:num}). This contradicts Lemma \ref{fini1}; hence there are
finitely many such $\ell$. Together with the induction hypothesis,
this proves that $b_\tau(\ell)$ stabilizes for $\ell$ sufficiently
large. This completes the proof of the lemma.
$\Box$
\medskip

For each $\mathbf{Q}^\tau$ as
above, let $c_\tau, b_\tau$ be as in Lemma \ref{lem:conv}. Let
$f_\mu = \lim\limits_{\stackrel{\to}{\ell < \mu}}f_\ell
= \sum_\tau c_\tau \mathbf{Q}^\tau$ and
$\delta_\mu = \lim\limits_{\stackrel{\to}{\ell < \mu}} \delta_\ell =
\sum_\tau b_\tau \mathbf{Q}^\tau$. We define the standard form
of $f$ of level $\gamma_\mu$ as $f=f_\mu + \delta_\mu$.
\medskip

This completes the construction of standard form of level
$\gamma_\ell$ for
$\ell \leq \lambda$.

\begin{prop}\label{standard2} Let
$$
f = f_\ell + \sum c_\alpha \mathbf{Q}^\alpha
$$
be a standard form of $f$ of level $\gamma_\ell$ and
$\gamma<\gamma_\ell$ an element of $\Phi$. Then
$\sum\limits_{\nu(\mathbf{Q}^\beta)=\gamma} c_\beta
\mathbf{Q}^\beta \notin P_{\gamma +}$.
\end{prop}

\noi The proof is entirely the same as the proof of the analogous
Proposition \ref{standard}.
\medskip

For each $\ell$, the part $f_\ell$ of a standard form of $f$ of
level $\gamma_\ell$ is uniquely determined. This is a straightforward
consequence of the proposition.

By Proposition \ref{standard2}, if $\gamma_\ell > \nu(f)$ then
$\nu(f)$ equals the smallest value of a monomial appearing in
the standard form of $f$ of level $\gamma_\ell$.
\medskip

\begin{thm} \label{graded-gen}
(1) Take $\gamma \in \Phi$, $\gamma < \gamma_\lambda$. Then
$\dsp{\frac{P_\gamma}{P_{\gamma +}}}$ is generated as a $k$-vector
space by $\{\init_\nu \mathbf{Q}^\beta \}$ where $\mathbf{Q}^\beta$
runs over the set of all standard monomials with respect to $\lambda$, satisfying $\nu(\mathbf{Q}^\beta)=\gamma$.

(2) The part of the graded $k$-algebra $\gr_\nu(A)$
  of degree strictly less than $\gamma_\lambda$ is generated by the
  initial forms of the approximate roots of $\mathcal{V}(\gamma_\lambda)$.
\end{thm}

\noi The same proof as that of Theorem \ref{Sgraded-gen} works here.
\medskip

Now, for each ordinal $\ell$, let $X=X_{\mathcal{V}(\gamma_\ell)}$
be a set of independent variables, indexed by
$\mathcal{V}(\gamma_\ell)$ and consider
the graded $k$-algebra $k\left[X_{\mathcal{V}(\gamma_\ell)}\right]$,
where we define $\deg\ X_j=\nu(Q_j)$. Let $P$ denote the homogeneous
monomial ideal of
$k\left[X_{\mathcal{V}(\gamma_\ell)}\right]$ generated by all the
monomials in $X_{\mathcal{V}(\gamma_\ell)}$ of
degree greater than or equal to $\gamma_{\ell+1}$. We have the natural
map
$$ \phi_\ell : \begin{array}{rcl}
  \frac{k\left[X_{\mathcal{V}(\gamma_\ell)}\right]}P &
  \rightarrow
  & \frac{\gr_\nu A}{P_{\gamma_{\ell+1}}} \\ X_j & \mapsto &\init_\nu
  Q_j \end{array}.$$

Now, for $\ell=0$, let $I_0=(0)$. For $\ell>0$, let $I_\ell$ denote the
ideal of $\frac{k\left[X_{\mathcal{V}(\gamma_\ell)}\right]}P$ generated by
$I_{<\ell}$ and all the homogeneous polynomials of the form

\begin{equation}
X^{\alpha_0}+\lambda_1X^{\alpha_1}+\lambda_2X^{\alpha_2}+\cdots+
\lambda_{b_0}X^{\alpha_{b_0}}\label{eq:distgen}
\end{equation}

\noi where
$\mathbf{Q}^{\alpha_0}+\lambda_1\mathbf{Q}^{\alpha_1} +
\lambda_2\mathbf{Q}^{\alpha_2}+\cdots+
\lambda_{b_0}\mathbf{Q}^{\alpha_{b_0}}$ is the homogeneous part of
least degree of $Ex(Q)$, $Q \in \Nu(\gamma_\ell) \cup
\Theta(\gamma_\ell)$.
\medskip

\noi Once again the proofs of Corollary \ref{cor:ker} and Corollary
\ref{cor:val-ideal} give the analogous corollaries :

\begin{cor}\label{ker1} We have $Ker\ \phi_\ell=I_\ell$.
\end{cor}

\begin{cor} \label{val-ideal}
Take an element $\gamma\in\Phi$, $\gamma
  < \gamma_\lambda$. The valuation ideal $P_{\gamma}$ is generated by all the
  generalized monomials of value $\geq \gamma$ in
    $\{ Q \ |\ (Q,Ex(Q)) \in \Psi(\gamma_\lambda)\}$. The ideal
    $P_{\gamma_\lambda}$ is generated by all the
  generalized monomials of value $\geq \gamma_\lambda$ in
    $\{ Q \ |\ (Q,Ex(Q)) \in \Psi(\gamma_\lambda)\cup
    \Theta(\gamma_\lambda)\}$.
\end{cor}

\section*{Part 2.  Separating ideal and connectedness}
\setcounter{section}{2}
\setcounter{subsection}{0}
\subsection{A description of the separating ideal.}
\label{DSI}

Let $A$ be a noetherian ring and $\alpha$ and $\beta$ points in
$\sper\ A$. The purpose of this section is twofold. First we prove a
general result on the behaviour of $<\alpha,\beta>$ under
localization. Secondly, we restrict attention to the case when $A$ is
regular and is either complete or $<\alpha,\beta>$ is primary to a
maximal ideal of $A$. In this case, we describe generators of the
separating ideal $<\alpha,\beta>$ as generalized monomials in those
approximate roots $Q_j$ which are common to $\nu_\alpha$ and $\nu_\beta$.
\medskip

We will need the following basic properties of the separating ideal,
proved in \cite{Mad1}:
\begin{prop}\label{sepideal} Let the notation be as above. We
have:

(1) $<\alpha,\beta>$ is both a $\nu_\alpha$-ideal and a $\nu_\beta$-ideal.

(2) $\alpha$ and $\beta$ induce the same ordering on $\frac
A{<\alpha,\beta>}$ (in particular, the set of
$\nu_\alpha$-ideals containing $<\alpha,\beta>$ coincides with the set
of $\nu_\beta$-ideals containing $<\alpha,\beta>$).

(3) $<\alpha,\beta>$ is the smallest ideal (in the sense of
inclusion), satisfying (1) and (2).

(4) If $\alpha$ and $\beta$ have no common specialization then
$<\alpha,\beta>=A$.
\end{prop}

\noi \textbf{Notation.} If $\mathfrak{p} \in \sper\ A$,
$\mathfrak{p}_\alpha \subset \mathfrak{p}$, the notation $\alpha
 A_{\mathfrak{p}}$ will stand for the point of $\sper\ A_{\mathfrak{p}}$
  with support $\mathfrak{p}_\alpha A_{\mathfrak{p}}$ and the total
  order on $\dsp{\frac{A_{\mathfrak{p}}}{\mathfrak{p}_\alpha
       A_{\mathfrak{p}}}}$ given by $\leq_\alpha$.

\begin{prop} \label{loca1} Let $A$ be a ring. Consider points
  $\alpha,\beta\in\mbox{Sper}A$ whose respective supports are
  $\mathfrak{p}_\alpha,\mathfrak{p}_\beta$ and let $\epsilon$ be a
  common specialization of $\alpha$ and
  $\beta$ with support $\mathfrak{p}$.

(1) We have
   $<\alpha,\beta>A_{\mathfrak{p}}=<\alpha A_{\mathfrak{p}},\beta
   A_{\mathfrak{p}}>$.

(2) Let $\mathfrak p$ be a prime ideal of $A$, containing $<\alpha,\beta>$. Then
\begin{equation} \label{loca}
<\alpha,\beta> \subset <\alpha,\beta>A_{\mathfrak{p}}\cap A.
\end{equation}
with equality if $<\alpha,\beta>$ is $\mathfrak p$-primary.

(3) If $\mathfrak{p}=\mathfrak{p}_\epsilon$ with $\epsilon$ the unique
common specialization of $\alpha$ and $\beta$ (in particular, whenever
$$
\mathfrak{p}= \sqrt{<\alpha,\beta>}
$$
and $\mathfrak{p}$ is maximal), we have equality in (\ref{loca}).
 \end{prop}
\begin{rek} In (2) of the Proposition, the special case of interest for
  applications is $\mathfrak{p}=\mathfrak{p}_\epsilon$, with
  $\epsilon\in\mbox{Sper }A$ a common specialization of $\alpha$ and
  $\beta$.
\end{rek}

\noi Proof: Let $f$ be a generator of $<\alpha,\beta>$ such that $f$
changes sign between $\alpha$ and $\beta$. Say, $f(\alpha)\geq0$ and
$f(\beta)\leq0$. As the orders on $A/\mathfrak{p}_\alpha$ and
$A_{\mathfrak{p}}/\mathfrak{p}_\alpha A_{\mathfrak{p}}$ are the same
(the quotient field is the same) --- and similarly for
$\mathfrak{p}_\beta$ --- $f$ changes sign between $\alpha
A_{\mathfrak{p}}$ and $\beta A_{\mathfrak{p}}$. Thus $f \in<\alpha
A_{\mathfrak{p}},\beta A_{\mathfrak{p}}>$.

Conversely, a generator of $<\alpha A_{\mathfrak{p}},\beta
A_{\mathfrak{p}}>$ is of the form $g/s$, $s\notin\mathfrak{p}$, such
that $\displaystyle{\frac{g}{s}(\alpha A_{\mathfrak{p}})\geq0}$ and
$\displaystyle{\frac{g}{s}(\beta A_{\mathfrak{p}})\leq0}$, for
instance. But, as $\mathfrak{p}$ is a specialisation of $\alpha$ and
$\beta$ and $s\notin\mathfrak{p}$, $s$ has the same sign on $\alpha$
and $\beta$ (and is non-zero at both points), so $g$ keeps different signs
on $\alpha$ and $\beta$ which means that $g\in<\alpha,\beta>$, and,
consequently,
$\displaystyle{\frac{g}{s}\in<\alpha,\beta>A_{\mathfrak{p}}}$. This
proves (1) of the Proposition.

(2) of the Proposition is a standard general statement about
localization of ideals at a prime ideal.

(3) of the Proposition follows immediately from the
fact that $\mathfrak{p}$ is the center of the valuation
$\nu_\alpha$ and $<\alpha,\beta>$ is a $\nu_\alpha$-ideal. $\Box$
\bigskip

Let $(A,\mathfrak{m},k)$ be a regular local ring and $\alpha$ and
$\beta$ two points of Sper($A$) having a common
specialization $\epsilon$ whose center is the maximal ideal
$\mathfrak{m}$ of $A$. Then $\nu_\alpha$ and $\nu_\beta$ are both
centered at $\mathfrak{m}$.

Let $\Phi_\alpha=\nu_\alpha(A\setminus \{0\})$ and
$\Phi_\beta=\nu_\beta(A\setminus \{0\})$.
Let $\gamma_{\alpha s}$ be
the $s$-th element of $\Phi_\alpha$ and similarly for $\beta$. Let
$P_{\gamma_{\alpha s}}$ denote the $\nu_\alpha$-ideal of value
$\gamma_{\alpha s}$
and similarly for
$P_{\gamma_{\beta s}}$. Let $r$
be the ordinal such that $\gamma_{\alpha r}=
\nu_\alpha(<\alpha,\beta>)$. Then $\gamma_{\beta r}
=\nu_\beta(<\alpha,\beta>)$ by Proposition
 \ref{sepideal}. We have $P_{\gamma_{\alpha s}}=P_{\gamma_{\beta s}}$
 for $s=1,\ldots,r$ by Proposition \ref{sepideal}.
\medskip

Let $Q_j(\alpha)$ denote
the $j$-th approximate root for $\nu_\alpha$ (in the case when $A$ is
complete $j$ is an ordinal rather than a natural number); we will
denote the monomials in these approximate roots by
$\mathbf{Q}(\alpha)^\gamma$; similarly for $Q_j(\beta)$ and
$\mathbf{Q}(\beta)^\gamma$. Let us consider the sequences of vectors
$\mathbf{m}_i=(m_{i1},m_{i2},\ldots,m_{it_{i\alpha}})$, $m_{ij}\in
P_{\gamma_{\alpha i}}/P_{\gamma_{\alpha,i+1}}$ which are the initial
forms of the monomials $\mathbf{Q}(\alpha)^{\alpha_{ij}}$ of value
$\gamma_{\alpha i}$ (see section \ref{AppS} and (\ref{ens-mon})). We do the
same with $\nu_\beta$ and write
$\mathbf{n}_1,\mathbf{n}_2,\ldots$ the
corresponding sequences of initial forms.

Let $M_{\alpha h}$ be the set of all the
generalized monomials in $\mathbf{Q}(\alpha)$, of value $\gamma_{\alpha h}$ with
respect to $\nu_\alpha$. Let $M_{\beta h}$ be the same kind of set with respect
to $\nu_\beta$. Now, let $s_{\alpha h}$ denote the cardinality of $M_{\alpha
h}$; similarly for $s_{\beta h}$.
\medskip

\noi For a given $\ell$, consider the following three conditions
(1)$_\ell$, (2)$_\ell$, (3)$_\ell$:

(1)$_\ell$ $s_{\alpha i}=s_{\beta i}$, $1\leq i\leq\ell$

(2)$_\ell$ $M_{\alpha i} = M_{\beta i}$ for $i\leq\ell$

(3)$_\ell$ For any $i\leq\ell$ and
$\bar{\lambda}_1,\ldots,\bar{\lambda}_{s_{\alpha i}}\in k$, the
sign on $\alpha$ of the linear combination
$\sum\limits_{j=1}^{s_{\alpha i}}\bar{\lambda}_jm_{ij}$ is the same
as the sign on $\beta$ of
$\sum\limits_{j=1}^{s_{\alpha i}}\bar{\lambda}_jn_{ij}$ (here we
adopt the convention that the sign can be strictly positive,
strictly negative or zero) where $m_{ij}, n_{ij}$ are the initial
forms of the monomials
$\mathbf{Q}(\alpha)^{\alpha_{ij}},\mathbf{Q}(\beta)^{\alpha_{ij}}$ in
the graded rings $\gr_{\nu_\alpha}(A),\ \gr_{\nu_\beta}(A)$.
\medskip
Note that if conditions (1)$_\ell$--(3)$_\ell$ hold then the set of
$k$-linear relations among the $m_{ij}$, $i\leq\ell$, is the same as
the set of $k$-linear relations among the $n_{ij}$.
\medskip

\begin{prop}\label{valsepid} The ordinal $r$ is the smallest ordinal
  $r'$
  such that at least one of the conditions (1)$_{r'}$--(3)$_{r'}$ does not
  hold.
\end{prop}

\noi Proof:
Let $r'$ be the smallest ordinal such that at least one of the conditions
(1)$_{r'}$--(3)$_{r'}$ does not hold.
 By definitions, we have $M_{\alpha r'} \ne \emptyset$ and
$M_{\beta r'} \ne \emptyset$. We have the following 2 possibilities:

First, suppose $M_{\alpha r'} \neq M_{\beta r'}$ (which includes the case
$s_{\alpha r'} \neq s_{\beta r'}$). Say, $M_{\alpha r'} \not\subset
M_{\beta r'}$. Take generalized monomials $\mathbf{Q}^\gamma\in
M_{\alpha r'} \setminus M_{\beta r'}$, and $\mathbf{Q}^\delta\in M_{\beta r'}$.
Then $\nu_\alpha(\mathbf{Q}^\gamma)\leq\nu_\alpha(\mathbf{Q}^\delta)$,
but $\nu_\beta(\mathbf{Q}^\gamma)>\nu_\beta(\mathbf{Q}^\delta)$.

Then there exists a linear combination, with coefficients in
$(A\setminus \mathfrak{m})$, of
$\mathbf{Q}^\gamma$ and $\mathbf{Q}^\delta$, of value
$\gamma_{\alpha r'}$ with respect to $\nu_\alpha$, which changes sign between
$\alpha$ and $\beta$. This shows that
$$
\nu_\alpha(<\alpha,\beta>)\leq\gamma_{\alpha r'}
$$
in this case.

The second case is $M_{\alpha r'}=M_{\beta r'}$ and there exist
$\bar{\lambda}_1,\ldots,\bar{\lambda}_{s_{\alpha r'}}$ such that the
sign on $\alpha$ of
$\sum\limits_{j=1}^{s_{\alpha r'}}\bar{\lambda}_jm_{r'j}$ differs from
the sign on $\beta$ of
$\sum\limits_{j=1}^{s_{\alpha r'}}\bar{\lambda}_jn_{r'j}$ (by assumption,
we are in the case $s_{\alpha r'}=s_{\beta r'}$). By a small perturbation
of the $\bar{\lambda}_j$ (for instance, by adding or subtracting a
``small'' element of $k$ to $\bar{\lambda}_1$), we can ensure both
that $\sum\limits_{j=1}^{s_{\alpha r'}}\bar{\lambda}_jm_{r'j}\ne0$ in
$gr_{\nu_\alpha}A$ and
$\sum\limits_{j=1}^{s_{\alpha r'}}\bar{\lambda}_jn_{r'j}\ne0$ in
$gr_{\nu_\beta}A$. But this gives an
$f=\sum\limits_{j=1}^{s_{\alpha r'}}\lambda_j\mathbf{Q}^{\alpha_{r'j}}\in
A$ which changes signs between $\alpha$ and $\beta$. We have
$\nu_\alpha(f)=\gamma_{\alpha r'}$ (and
$\nu_\beta(f)=\gamma_{\beta r'}$), so
$\nu_\alpha(<\alpha,\beta>)\leq\gamma_{\alpha r'}$ also in this case.

Now take an $f \in A$ with $\nu_\alpha(f)<\gamma_{\alpha r'}$. Then
$f\in P_{\gamma_{\alpha s}}$,
\begin{equation}
\gamma_{\alpha s}<\gamma_{\alpha r'},\label{eq:lessvalue}
\end{equation}
so $\init_{\nu_\alpha}(f)\in
P_{\gamma_{\alpha s}}/P_{\gamma_{\alpha s+}}$. By theorem
\ref{graded-gen},
$\init_{\nu_\alpha}(f)$ is a $k$-linear combination of
$m_{s1}$, \dots, $m_{st_{s\alpha}}$. By (\ref{eq:lessvalue}) and the
definition of $r'$, this linear combination has the same sign for
$\alpha$ and for $\beta$ (in other words,
$P_{\gamma_{\alpha s}}/P_{\gamma_{\alpha s+}}=
P_{\gamma_{\beta s}}/P_{\gamma_{\beta s +}}$ with same order induced
by $\alpha$ and by $\beta$. This means that $\init_{\nu_\alpha}(f)$
has the same sign on $\alpha$ and $\beta$, so
$\nu_\alpha(<\alpha,\beta>)\geq\gamma_{\alpha r'}$). This completes the
proof. $\Box$

\begin{cor} Let $\alpha,\beta\in\sper(A)$, both centered in the
   maximal ideal. Let $r$ be as above. Denote by $\gamma=
\gamma_{\alpha r}$ the $\nu_\alpha$-value of $<\alpha,\beta>$. Let
   $Q_1,\ldots,Q_q$ be the common approximate roots of the
   valuations $\nu_\alpha$ and
   $\nu_\beta$. Then $<\alpha,\beta>$ is
generated by the generalized monomials in $Q_1,\ldots,Q_q$ of
$\nu_\alpha$-value $\geq\gamma$ (and the same with $\nu_\beta$ instead of
$\nu_\alpha$).
\end{cor}

\noi Proof: As $<\alpha,\beta>$ is a $\nu_\alpha$-ideal (and a
$\nu_\beta$-ideal), this is a consequence of Corollary
\ref{val-ideal}.
\medskip

\begin{deft} For a graded algebra $G$, we define
$$
\displaystyle{G^*=\left\{\left.\frac{f}{g}\ \right|\ f,g\in
     G,g\neq0\mbox{ and homogeneous}\right\}/\sim}.
$$ where $ \dsp{\frac{f}{g} \sim \frac{f'}{g'}}$ whenever $fg'=f'g$.
\end{deft}

\noi\textbf{The Alvis--Johnston--Madden example.} Let us consider
$\alpha$ and $\beta$ in $\sper(R[x,y,z])$ given by curvettes
\begin{eqnarray}
x(t)&=&t^6,\\
y(t)&=&t^{10}+ut^{11},\\
z(t)&=&t^{14}+t^{15}
\end{eqnarray}
where $u$ takes 2 distinct values $u_\alpha>2$ and
$u_\beta>2$. Applying the above procedure, we show that
$\nu_\alpha(<\alpha,\beta>)=31$.

\noi Indeed, we have $Q_1=x,Q_2=y,Q_3=z$ for $\alpha$ and $\beta$. The
first level approximate roots are
\begin{eqnarray}
Q_4&=&y^2-xz=(2u-1)t^{21}+u^2t^{22},\\
Q_5&=&yz-x^4=(u+1)t^{25}+ut^{26},\\
Q_6&=&z^2-x^3y=(2-u)t^{29}+t^{30}
\end{eqnarray}
for both $\alpha$ and $\beta$. Let $T$ denote the preimage of $\init_vt$
under the natural map $$(gr_{\nu_\alpha} R[x,y,z])^* \hookrightarrow (\gr_v
R[[t]])^*,$$ so that
$$
(gr_{\nu_\alpha} R[x,y,z])^*\cong(R[T])^*.
$$
Then $\init_{\nu_\alpha}(yQ_4)=(2u_\alpha-1)T^{31}$ and
$\init_{\nu_\alpha}(xQ_5)=(u_\alpha+1)T^{31}$, and similarly for
$\beta$. Since $u_\alpha\ne u_\beta$, the matrix
$$
\left(\begin{array}{cc}
(2u_\alpha-1)&(u_\alpha+1)\\
(2u_\beta-1)&(u_\beta-1)
\end{array}\right)
$$
is non-singular, so there exists an $R$-linear combination of
$\init_{\nu_\alpha}(yQ_4)$ and $\init_{\nu_\alpha}(xQ_5)$ which is
strictly positive on $\alpha$ and strictly negative on
$\beta$. According to Proposition \ref{valsepid},
$$
\nu_\alpha(<\alpha,\beta>)\le31.
$$
One can check that 31 is the lowest value for which either there is a
linear combination of generalized monomials with this property or the
set of monomials of that value for $\alpha$ does not equal the
corresponding set for $\beta$, so that in fact
$\nu_\alpha(<\alpha,\beta>)=31$.

For the next approximate root
\begin{equation}
Q_7=yQ_4+\frac{2u-1}{u+1}Q_5\label{eq:Q8},
\end{equation}
we have $Q_7(\alpha) \neq Q_7(\beta)$.

\subsection{Some sets which are conjecturally connected}
\label{section:connexe}
\noi Let $(A,\mathfrak{m},k)$ be a regular local ring.
Take $\alpha, \beta \in \sper A$, both centered at $\mathfrak{m}$,
and elements $f_1,\ldots,f_r \in A \setminus <\alpha,\beta>$. The
Connectedness Conjecture \ref{conn} asserts that
there exists a connected set $C$, containing $\alpha,\beta$, such that
$C$ is disjoint from the zero set of $f_1\cdots f_r$.
\medskip

Assume that either $A$ is complete or $\sqrt{<\alpha,\beta>}=
\mathfrak{m}$.
\medskip

In this section, we describe a set $C$, which contains $\alpha, \beta$,
disjoint from the set $f_1\cdots f_r =0$, and which we conjecture to
be connected. Under the above assumptions, this reduces the
Connectedness Conjecture for $\alpha$ and $\beta$ to proving the
connectedness of $C$.

Let $\mathbf{Q}_\Lambda = \{Q_\lambda, \lambda \in \Lambda\}$ be the
approximate roots common to $\alpha$ and $\beta$. Let
$\mathbf{Q}^{\gamma_1},\mathbf{Q}^{\gamma_2},\ldots$ be the list of
monomials in $\mathbf{Q}_\Lambda$, arranged in the increasing order of
the $\nu_\alpha$ values. There exists an ordinal $s$ such that
$<\alpha,\beta>$ is generated by the set $\{\mathbf{Q}^{\gamma_j}; \ j
\leq s,\ \mathbf{Q}^{\gamma_j} \in <\alpha,\beta> \}$. Let $\sigma$ be the
unique ordinal such that $\mathbf{Q}^{\gamma_a} \notin <\alpha,\beta>$
for $a<\sigma$ and $\mathbf{Q}^{\gamma_\sigma},\mathbf{Q}^{\gamma_{\sigma+1}},\ldots \in <\alpha,\beta>$.
\medskip

Next, we study the standard form of $f_i$ of level
$\nu_\alpha(<\alpha,\beta>)$. In the case when $A$ is complete, this
standard form may contain infinitely many generalized monomials
$\mathbf{Q}^\gamma$. Since $A$ is noetherian, we can choose a finite
subset $\mathbf{Q}^{\epsilon_{ji}}$, $1 \leq j \leq n_i$, of these
monomials such that all of the others lie in the ideal
$(\mathbf{Q}^{\epsilon_{ji}}, 1\leq j \leq n_i)A$. For
$i \in \{1,\ldots,r\}$, let
\begin{equation}
f_i=\sum_{j=1}^{m_i} b_{ji}\mathbf{Q}^{\theta_{ji}} + \sum_{j'=1}^{n_i}
  c_{j'i} \mathbf{Q}^{\epsilon_{j'i}}
\end{equation}
be the standard expansion of $f_i$ of level
$\nu_\alpha(<\alpha,\beta>)$ where
$\nu_\alpha(\mathbf{Q}^{\theta_{ji}}) = \nu_\alpha(f_i) <
\nu_\alpha(\mathbf{Q}^{\epsilon_{j'i}})$ for all $j
\in\{1,\ldots,m_i\}$ and $j' \in \{1,\ldots,n_i\}$.

\begin{rek}
1. If $k=k_\alpha$ (in particular, if $k$ is real closed), then
$m_i=1$.

2. By Proposition \ref{standard2}, $\sum_{j=1}^{m_i} b_{ji}
\init_{\nu_\alpha}\mathbf{Q}^{\theta_{ji}} \neq 0$.
\end{rek}

\begin{conj} 1. Let \begin{equation} \label{ensC}
C = \left\lbrace \delta \in \sper A
    \ \left| \begin{array}{cl} \ \nu_\delta(\mathbf{Q}^{\theta_{ji}}) <
        \nu_\delta(\mathbf{Q}^{\epsilon_{j'i}}) & \mbox{ for all } j \in
        \{1,\ldots,m_i\},\ j' \in \{1,\ldots,n_i\}
    \\
     sgn_\delta(Q_q) = sgn_\alpha(Q_q) & \mbox{ for all } Q_q
     \mbox{ appearing in } \mathbf{Q}^{\theta_{ji}} 
\\
      sgn_\delta(\sum_{j=1}^{m_i} b_{ji} \mathbf{Q}^{\theta_{ji}}) =
&     sgn_\alpha(\sum_{j=1}^{m_i} b_{ji} \mathbf{Q}^{\theta_{ji}})
 \end{array} \right. \right\rbrace.
\end{equation}Then $C$ is connected.

2. Let $C'$ defined by the inequalities
\begin{equation} \label{eq:dominant1}
\left|\sum_{j=1}^{m_i} b_{ji}\mathbf{Q}^{\theta_{ji}}\right| >_\delta
n_i|\mathbf{Q}^{\epsilon_{j'i}}| \
\forall i \in \{1,\ldots,r\}, \ \forall j' \in \{1,\ldots,n_i\}
\end{equation}
and the two sign conditions appearing in (\ref{ensC}). Then $C'$ is
connected.
\end{conj}

\begin{rek}\label{Cworks} 1. We have $\alpha, \beta \in C$.

2. $C \cap \{f_1\cdots f_r = 0\} = \emptyset$. Indeed,
inequalities (\ref{ensC}) imply that, for every $\delta \in
C$, $f_i$ has the same sign as $\sum_{j =1}^{m_i} b_{ji}
\mathbf{Q}^{\theta_{ji}}$; in particular, none of the $f_i$ vanish on $C$.

3. Either of those conjectures implies the Connectedness Conjecture.
\end{rek}

\section*{\Large Part 3.  A proof of the conjecture for arbitrary
  regular 2-dimensional rings.}
\setcounter{section}{3}
\setcounter{subsection}{0}

We start with a general plan of the proof and an outline of different
sections of Part 3. In \S\ref{Zariski} we recall Zariski's theory of
complete ideals. We explain how the construction of approximate roots
in arbitrary dimension restricts to the special case of dimension 2
(and that the standard construction in dimension 2 is, indeed,
recovered from the general one as a special case) and prove some
general lemmas about approximate roots in regular two dimensional
local rings and their behaviour under sequences of point blowings
up. In \S\ref{realsurfaces} we define the notion of real geometric
surfaces which are glued from affine charts of the form $\sper\ A_j$,
where $A_j$ is a regular two-dimensional ring, in order to be able to
talk about point blowings up of $\sper\ A$. We also define the notion
of a segment on the exceptional divisor of a blowing up and prove that
such a segment is connected; another notion useful later in the proof
is that of a maximal segment. One slightly delicate point here is that
since the residue field $k$ of $A$ is not assumed real closed we need
to fix an order on $k$ and always restrict attention to points of the
real spectra of various $A_j$ which induce the given order on $k$. The
bulk of the proof \textit{per se} is contained in
\S\S\ref{proofPB}--\ref{proof}. As explained above,
our problem is one of proving connectedness (resp. definable
connectedness) of the set $C$.
\medskip

In \S\ref{proofPB} we use Zariski's theory and other results from
\S\ref{Zariski} to construct a sequence of point blowings up which
transform $C$ into a quadrant, that is, a set $\tilde U$ of all points
$\delta$ of a suitable affine chart $\sper\ A_j$ centered at the
origin satisfying either $x'(\delta)>0$, $y'(\delta)>0$, or just
$x'(\delta)>0$. In \S\ref{connectedexcellent} we use results from
\cite{And} to prove connectedness of $\tilde U$ by reducing it to that
of a quadrant in the usual Euclidean space, assuming that $A$ is
excellent. In \S\ref{proof} we prove the definable connectedness of
$\tilde U$ (without any excellence assumptions) after introducing a
new object called the graph associated
to $\tilde U$ and a finite sequence of point blowings up of $\sper\
A$.

\subsection{Approximate roots  in dimension 2 and Zariski's
theory.}\label{Zariski} \label{subs-proof}

In the special case of regular 2-dimensional local rings, the theory
of approximate roots is well known: see, for instance \cite{Zar},
Appendix 5 or \cite{Spi1}. We briefly recall the construction here
since it is much simpler than in the general case.
\medskip

We start with two purely combinatorial lemmas about semigroups. Take
an integer $g \geq 2$.
\begin{lem} \label{lem:semigp}
Let $\beta_1,\beta_2\dots,\beta_g$ be positive elements in some
ordered group. Let $\alpha_j$, $j \in \{2,\ldots,g\}$ be positive
integers. Assume
\begin{equation} \label{eq:semigp}
\beta_i \geq \alpha_{i-1}\beta_{i-1}, \ i \in \{3,\ldots g\}.
\end{equation}
 Let $\gamma_1,\ldots,\gamma_g$ be integers such that $0 \leq \gamma_j
 < \alpha_j$ for $2 \leq j \leq g$ and $\sum_{j=1}^g \gamma_j\beta_j
 \geq \alpha_g\beta_g$. Then $\gamma_1 > 0$.
\end{lem}

\noi Proof : We prove by descending induction that $\sum_{j=1}^i
\gamma_j\beta_j \geq \alpha_i\beta_i$ for $i\geq 2$. The case $i=g$ is
given by hypothesis. Assume then that $\sum_{j=1}^{i+1}
\gamma_j\beta_j \geq \alpha_{i+1}\beta_{i+1}$. Subtracting
$\gamma_{i+1}\beta_{i+1}$ and using the fact that $\gamma_{i+1} <
\alpha_{i+1}$, we obtain $\sum_{j=1}^i \gamma_j\beta_j \geq
(\alpha_{i+1}-\gamma_{i+1})\beta_{i+1} \geq \alpha_i\beta_i$. This
completes the induction. So for $i=2$, we obtain $\gamma_1\beta_1 +
\gamma_2 \beta_2 \geq \alpha_2\beta_2$; subtracting $\gamma_2\beta_2$
and using the fact that $\gamma_2 < \alpha_2$, we get $\gamma_1\beta_1
\geq (\alpha_2-\gamma_2)\beta_2 > 0$, hence $\gamma_1 > 0$.
$\Box$
\medskip

\noi\textbf{Notation.} Let $\beta_1,\beta_2\dots,\beta_g$ be positive
elements in some ordered group.
We will denote by $(\beta_1,\ldots,\beta_{i-1})$ the group generated by
$\beta_1,\ldots,\beta_{i-1}$ and by $sg(\beta_1,\ldots,\beta_{i-1})$
the semigroup generated by $\beta_1,\ldots,\beta_{i-1}$, that is, the
semigroup formed by all the $\N$-linear combinations of
$\beta_1,\ldots,\beta_{i-1}$. For $i\in\{2,\dots,g\}$, $\alpha'_i$ will denote
the smallest positive integer such that
$\alpha'_i\beta_i \in
(\beta_1,\ldots,\beta_{i-1})$. If there is no such
  integer, we put $\alpha'_i = \infty$.
Write
\begin{equation}\label{eq:relation}
\alpha'_i\beta_i  = \sum_{j=1}^{i-1}\alpha_{ji}\beta_j\quad\text{
  where
}\alpha_{ji}
  \in \Z.
\end{equation}

\begin{lem}\label{combin} Let $\beta_1,\beta_2\dots,\beta_g$ be
  positive rational numbers such that
%
$\beta_g \geq \alpha'_{g-1}\beta_{g-1}$. If $g\ge3$, assume that
\begin{equation}\label{gp=sg}
\{a\in(\beta_1,\ldots,\beta_{g-1})\ |\
a\ge\alpha'_{g-1}\beta_{g-1}\}=\{a\in sg(\beta_1,\ldots,\beta_{g-1})\
|\ a\ge\alpha'_{g-1}\beta_{g-1}\};
\end{equation}
in particular, we can choose
$\alpha_{jg}\ge0$ for all $j\in\{1,\dots,g-1\}$ in
(\ref{eq:relation}) when $i=g$. Then
\begin{equation}\label{gp=sg1}
\{a\in(\beta_1,\ldots,\beta_{g})\ |\
a\ge\alpha'_{g}\beta_{g}\}=\{a\in sg(\beta_1,\ldots,\beta_{g})\
|\ a\ge\alpha'_{g}\beta_{g}\}.
\end{equation}
\end{lem}

\noi\textbf{Proof.} Multiplying all the $\beta_i$ by the same rational number does
not change the problem, so we may assume that
$\beta_1,\beta_2,\ldots,\beta_g$ are positive integers, such that
$\gcd(\beta_1,\beta_2,\ldots,\beta_g)=1$.

For $g=2$, we have $\alpha'_2=\beta_1$. If $a \in (\beta_1)$ and $a
\geq \beta_1 \beta_2$, then $a >0$, hence $a \in sg(\beta_1)$; thus
$$\{a \in (\beta_1)\ |\ a \geq \beta_1\beta_2\} \subset \{a \in
sg(\beta_1)\ |\ a \geq \beta_1\beta_2\},$$ the opposite inclusion
being obvious.

Assume that $g\geq3$. Write
\begin{equation}\label{eq:relation1}
\alpha'_{g}\beta_{g}= \sum_{j=1}^{g-1}\alpha_{jg}\beta_j.
\end{equation}

To prove
(\ref{gp=sg1}), let $\beta = \gamma_1 \beta_1 + \gamma_2
\beta_2+\cdots+\gamma_{g} \beta_{g}$ be an element of
$$
\{a\in (\beta_1,\dots\beta_{g})\ |\ a \geq \alpha'_{g}\beta_{g}\}.
$$
Using the relation (\ref{eq:relation1}) we can write, for each $n \in \Z$,
$$
\beta = \sum\limits_{j=1}^{g-1}(\gamma_j -n\alpha_{jg})\beta_j+
(\gamma_{g} + n\alpha'_{g})\beta_{g}
= \sum\limits_{j=1}^{g-1}\gamma'_j\beta_j + (\gamma_{g} + n
\alpha'_{g})\beta_{g}.
$$
After replacing $\gamma_{g}$ by
$\gamma_{g} + n \alpha'_{g}$ for a suitable $n\in \Z$, we may assume that
$0 \leq \gamma_{g} < \alpha'_{g}$. Since $\beta \geq
\alpha'_{g}\beta_{g}$, this implies that
\begin{equation}\label{eq:semigp1}
\sum\limits_{j=1}^{g-1}\gamma_j\beta_j \geq \beta_{g} \geq\alpha'_{g-1}
\beta_{g-1}.
\end{equation}
By (\ref{gp=sg})
, we may take $\gamma_j\geq 0$ in
(\ref{eq:semigp1}). This completes the proof of the
Lemma. $\Box$

\begin{cor} \label{cor:combin}
Let $\beta_1,\ldots,\beta_g$ be positive rational numbers satisfying
\begin{equation} \label{eq:semigp2}
\beta_i \geq \alpha'_{i-1}\beta_{i-1}, \ i \in \{3,\ldots g\}.
\end{equation}
Then equalities (\ref{gp=sg}) and (\ref{gp=sg1}) hold.
\end{cor}

\noi Proof : For $i=3$, (\ref{gp=sg}) is immediate. Now the corollary
follows from Lemma \ref{combin} by induction on $i$.
$\Box$
\medskip

Let $\nu$ be a valuation centered at $A$ and let $(x,y)$ be a
$\nu$-prepared system of coordinates, such that $\nu(x) =
\nu(\mathfrak{m})$. In what follows, we will omit the description of
$\Nu(\gamma), \Lambda(\gamma), \Theta(\gamma)$, since in the
simplified situation of $n=2$, the sets $\Psi(\gamma)$ suffice to
carry out the entire construction.
\medskip

Put $Q_1=x$, $Ex(Q_1)= x$, $Q_2=y$, $Ex(Q_2)=y$  and
$\beta_i=\nu(Q_i)$, $i\in\{1,2\}$. If $\beta_1, \beta_2$ are
rationally independent, then $\alpha'_2= \infty$ and the construction
stops, there are no more approximate roots. In this case, all the
$\nu$-ideals are generated by monomials in ($x,y)$. Assume then
$\alpha'_2 < \infty$. This means that there is a relation
$\alpha'_2 \beta_2= \alpha_{12} \beta_1$ for a positive integer
$\alpha_{12}$.

Let $\alpha'_2$ and $\alpha_{12}$
be as above. Let $\Psi(\beta_1)=
\emptyset$. For $\gamma \in \Phi$, $\beta_1 < \gamma < \beta_2$,
$\Psi(\gamma) = \{x\}$ and $\Psi(\beta_2+)= \{x,y\}$. Let $k_1=k$,
$\dsp{k_2 = k\left(\frac{\init_\nu\left(
    Q_2^{\alpha'_2}\right)}{\init_\nu\left(Q_1^{\alpha_{12}}\right)}\right)}$.

 We prove that (\ref{gp=sg}) is satisfied for $i=3$.
Let $\beta = \gamma_1 \beta_1 + \gamma_2 \beta_2$ be an element of
$$
\{a\in (\beta_1,\beta_2)\ |\
 a\ge\alpha'_2\beta_2\}.
$$
As $\alpha'_2 \beta_2 = \alpha_{12} \beta_1$,
 we have, for each $n \in \Z$, $\beta = (\gamma_1-n\alpha_{12})\beta_1
 +(\gamma_2 + n \alpha'_2)\beta_2$. After replacing $\gamma_2$ by
$\gamma_2 + n \alpha'_2$ for a suitable $n$, we may assume that $0 \leq
 \gamma_2 < \alpha'_2$. Since $\beta\ge \alpha'_2 \beta_2$, this
 implies that $\gamma_1 \geq 0$.

\noi Then we have $\nu\left(Q_2^{\alpha'_2}\right)=
\nu\left(Q_1^{\alpha_{12}}\right)$, hence the image of
$\displaystyle{\frac{\init_\nu(Q_2^{\alpha'_2})}{\init_\nu(Q_1^{\alpha_{12}})}}$
in $k_\nu$ is not zero.
If
$\displaystyle{\frac{\init_\nu(Q_2^{\alpha'_2})}{\init_\nu(Q_1^{\alpha_{12}})}}$
is algebraic over $k$, this means that it satisfies
an algebraic equation of the form
\begin{equation}
X^d+ \sur{a}_1X^{d-1} + \cdots + \sur{a}_d =0, \ \sur{a}_i \in k.
\end{equation}
Let $a_i$ be a representative of $\sur{a}_i$ in $A$. Let
$\alpha_2=d\alpha'_2$ and
\begin{equation}
Q_{3}^{(1)} = Q_2^{\alpha_2} + \sum_{\ell=1}^d
a_\ell Q_2^{\alpha'_2(d-\ell)} Q_1^{\alpha_{12}\ell}.
\end{equation}
\noi The expression $Ex\left(Q_{3}^{(1)}\right)$ is just the right hand side of
this formula.

Let
$
\beta_{3}^{(1)} = \nu\left(Q_{3}^{(1)}\right)$, then $\beta_3^{(1)} >
\nu(Q_2^{\alpha_2}) =  \alpha_2\beta_2 \geq \alpha'_2\beta_2$
and the elements $\left(\beta_1,\beta_2,\beta_3^{(1)}\right)$ satisfy
the conclusion of Lemma \ref{combin}.

By construction, $\alpha_2 \beta_2$ is the smallest
element of $\Phi$ such that the monomials
$$
\left\{Q_1^{\gamma_1}Q_2^{\gamma_2}\ \left|\ \nu(Q_1^{\gamma_1}Q_2^{\gamma_2}) =
\gamma_1\beta_1 + \gamma_2\beta_2= \alpha_2\beta_2\right.\right\}
$$
are $k$-linearly dependent. The unique $k$-linear dependence relation is given by
$Q_3^{(1)}$. Hence, according to the general construction of \S 2, we
have $\Theta(\beta) =\{Q_3^{(1)}\} $ for $\alpha_2 \beta_2 \geq \beta
\geq \beta_3^{(1)}$ and $\Psi(\beta_3^{(1)}+) =
\{Q_1,Q_2,Q_3^{(1)}\}$.
\medskip

Assume that $i\geq 3$ and that elements $Q_1,\ldots,Q_{i-1},Q_i^{(j)}$
are already defined. Let
\begin{eqnarray}
\beta_q&=& \nu(Q_q),\\
\beta_i^{(j)}&=&\nu\left(Q_i^{(j)}\right).
\end{eqnarray}
Assume that the initial form $\init_\nu Q_q$ is
algebraic over $k[\init_\nu Q_1,\ldots,\init_\nu Q_{q-1}]$ for $q \in
\{2,\ldots,i-1\}$. Let $\alpha_q$ denote the degree of its minimal
polynomial. Note that, in particular, all of
$\beta_2,\dots,\beta_{i-1}$ are rational multiples of $\beta_1$.
Assume that $\beta_q > \alpha_{q-1} \beta_{q-1}$, $q \in
\{3,\ldots,i-1 \}$ and  $\beta_i^{(j)} >
\alpha_{i-1}\beta_{i-1}$. Assume that, in the notation of \S
\ref{AppS}, we have
$$
\Psi\left(\beta_{i\ +}^{(j)}\right)
=\left\{Q_1,\ldots,Q_{i-1},Q_i^{(j)}\right\}.
$$
A monomial $\dsp{\prod_{\ell=1}^{i-1} Q_\ell^{\epsilon_\ell}}$ is
  \textbf{standard} if
\begin{equation}\label{eq:condit}0 \leq
    \epsilon_\ell < \alpha_\ell \mbox{ for } \ell \in
    \{2,\ldots,i-1\}.
\end{equation}
This allows us to extend the notion of
standard to monomials with $\epsilon_1 <0$: such a monomial is called
standard if (\ref{eq:condit}) is satisfied. Similarly, we may talk
about standard monomials in $\init_\nu Q_1,\ldots, \init_\nu
Q_{i-1}$.

Assume, in addition, that we have defined elements
$z_2,\ldots,z_{i-1} \in k_\nu$, algebraic over $k$, where $z_\ell$ is
a $k$-linear combination of standard monomials in $\init_\nu Q_1, \ldots,
\init_\nu Q_\ell$ of degree 0. Let $k_\ell= k(z_2,\ldots,z_\ell)$. We
obtain a tower of finite field extensions $k=k_1 \subset k_2 \cdots
\subset k_{i-1} \subset k_\nu$.
\smallskip

If $\init_\nu Q_i^{(j)}$ is transcendental over $k[\init_\nu Q_1,
\ldots, \init_\nu Q_{i-1}]$, put $Q_i=Q_i^{(j)}$ and the construction
stops.

\noi Assume $\init_\nu Q_i^{(j)}$ is algebraic over $k[\init_\nu Q_1,
\ldots, \init_\nu Q_{i-1}]$. Then  $\beta_i^{(j)} \in
\sum\limits_{q=1}^{i-1}\Q \beta_q$. Let $\alpha^{, (j)}_i$ be the smallest
positive integer such that $\alpha_i^{,(j)} \beta_i^{(j)} \in
(\beta_1,\ldots,\beta_{i-1})$.

\noi Then $\nu\left(\left(Q_i^{(j)}\right)^{\alpha_i^{,(j)}}\right)=
\nu\left(\prod\limits_{r=1}^{i-1}Q_r^{\alpha_{ri}^{(j)}}\right)$, hence
the image of
$\dsp{\frac{\init_\nu\left(Q_i^{(j)}\right)^{\alpha_i^{,(j)}}}
{\init_\nu\prod_{r=1}^{i-1}Q_r^{\alpha_{ri}^{(j)}}}}$
in $k_\nu$ is not zero. By Corollary \ref{cor:combin}, we may take
$\alpha_{ri}^{(j)} \ge 0$ for $1\le r\le i-1$.

\noi The assumption on $\init_\nu Q_i^{(j)}$ implies that
$\dsp{\frac{\init_\nu\left(Q_i^{(j)}\right)^{\alpha_i^{,(j)}}}
{\init_\nu\prod_{r=1}^{i-1}Q_r^{\alpha_{ri}^{(j)}}}}$
satisfies an algebraic equation of the form
\begin{equation}
X^d+ \sur{a}_1X^{d-1} + \cdots + \sur{a}_d =0, \ \sur{a}_\ell \in k_{i-1}.
\end{equation}

For $\ell \in \{1,\ldots,d\}$, write
\begin{equation} \label{eq:expositif}
\sur{a}_\ell
\left(\prod_{r=1}^{i-1} \init_\nu Q_r^{\alpha_{ri}^{(j)}} \right)^\ell
= \sum_{\gamma= (\gamma_1,\ldots,\gamma_{i-1})} \sur{b}_{\ell \gamma}
\prod_{r=1}^{i-1}\init_\nu Q_r^{\gamma_r}
\end{equation}
as a $k$-linear combination of standard monomials. By Lemma
\ref{lem:semigp},
we have $\gamma_1 \geq 0$ whenever $\sur{b}_{\ell \gamma} \neq 0$.

Let $b_{\ell \gamma}$ be a representative of $\sur{b}_{\ell \gamma}$
in $A$. Let $\alpha_i^{(j)}=d\alpha_i^{,(j)}$ and
\begin{equation} \label{eqn:appr2-1}
Q= \left(Q_i^{(j)}\right)^{\alpha_i^{(j)}} + \sum_{\ell=1}^d
 \left(\sum_{\gamma = (\gamma_1,\ldots,\gamma_{i-1})} b_{\ell \gamma}
   \prod_{r=1}^{i-1} Q_r^{\gamma_r}\right)Q_i^{\alpha_i^{, (j) (d-\ell)}}.
\end{equation}
Then
\begin{equation}\label{eqn:appr2-2}
\nu(Q) > \nu\left(\left(Q_i^{(j)}\right)^{\alpha_i^{(j)}}\right) =
\alpha_i^{(j)}\beta_i^{(j)} \ge \alpha_i^{,(j)} \beta_i^{(j)} >
\alpha_{i-1}\beta_{i-1} \ge \alpha^,_{i-1} \beta_{i-1}.
\end{equation}
\smallskip

\noi If $\init_\nu Q_i^{(j)} \notin k[\init_\nu Q_1,\ldots, \init_\nu
Q_{i-1}]$ (which is equivalent to saying that $\alpha_i^{(j)} > 1$), put
$Q_i=Q_i^{(j)}$, $Ex(Q_i) =Ex(Q_i^{(j)})$, $\beta_i=\beta_i^{(j)}$,
$\alpha_i=\alpha_i^{(j)}$, $\alpha_i^,=\alpha_i^{,(j)}$,
$Q_{i+1}^{(1)}=Q$, $\beta_{i+1}^{(1)}=
\nu\left(Q_{i+1}^{(1)}\right)$.
\smallskip

\noi Formulae (\ref{eqn:appr2-1}) and (\ref{eqn:appr2-2}) become
\begin{equation} \label{eqn:appr2-3}
Q_{i+1}^{(1)} = Q_i^{\alpha_i} + \sum_{\ell=1}^d
 \left(\sum_{\gamma = (\gamma_1,\ldots,\gamma_{i-1})} b_{\ell \gamma}
   \prod_{r=1}^{i-1} Q_r^{\gamma_r}\right)Q_i^{\alpha_i^{,} (d-\ell)}.
\end{equation}
and
\begin{equation}\label{eqn:appr2-4}
\beta_{i+1}^{(1)} = \nu\left(Q_{i+1}^{(1)}\right) > \nu(Q_i^{\alpha_i}) =
\alpha_i\beta_i \ge \alpha_i^{,} \beta_i.
\end{equation}
\noi The expression $Ex\left(Q_{i+1}^{(1)}\right)$ is just the right hand side of
(\ref{eqn:appr2-3}).
\smallskip

\noi For $\beta_i^{(j)}<\gamma \le \beta_{i+1}^{(1)}$,
we have  $\Psi(\gamma) = \Psi\left(\beta_i^{(j)}+\right)$ and
$\Psi\left(\beta_{i+1}^{(1)}+\right) =
\left\{Q_1,\ldots,Q_i,Q_{i+1}^{(1)}\right\}$. Moreover,
the elements $\left(\beta_1,\dots,\beta_{i+1}^{(1)}\right)$ satisfy
the hypothesis of Corollary \ref{cor:combin}, hence also its
conclusion.

\noi If $\init_\nu Q_i^{(j)} \in k[\init_\nu Q_1,\ldots, \init_\nu
Q_{i-1}]$ (which is equivalent to saying that $\alpha_i^{(j)} = 1$), put
$Q_i^{(j+1)}=Q$ and $\beta_i^{(j+1)}= \nu\left(Q_i^{(j+1)}\right)$.

\noi Formulae (\ref{eqn:appr2-1}) and (\ref{eqn:appr2-2}) become
\begin{equation} \label{eqn:appr2-5}
Q_i^{(j+1)} = Q_i^{(j)} + \sum_{\gamma=(\gamma_1,\ldots,\gamma_{i-1})}
b_{1 \gamma} \prod_{r=1}^{i-1}Q_r^{\gamma_r}
\end{equation}
and
\begin{equation}\label{eqn:appr2-6}
\beta_i^{(j+1)} = \nu\left(Q_i^{(j+1)}\right) > \beta_i^{(j)} >
\alpha^,_{i-1} \beta_{i-1}.
\end{equation}
\noi The expression $Ex\left(Q_{i}^{(j+1)}\right)$ is just the right hand side of
(\ref{eqn:appr2-5}).
\smallskip

\noi For $\beta_i^{(j)} <\gamma \le \beta_i^{(j+1)}$
we have
\begin{eqnarray}
\Psi(\gamma) &=& \Psi\left(\beta_i^{j}+\right)\quad\text{ and}\\
\Psi\left(\beta_i^{(j+1)}+\right) &=&
\left\{Q_1,\ldots,Q_{i-1},Q_i^{(j+1)}\right\}.
\end{eqnarray}
Moreover, the elements
$\left(\beta_1,\dots,\beta_{i-1},\beta_i^{(j+1)}\right)$
satisfy the hypothesis of Corollary \ref{cor:combin}, hence also its
conclusion.
\medskip

\begin{rek} Either the process stops after a finite number of steps or we
obtain an infinite sequence
\begin{equation}
\mathbf{Q}=Q_1,Q_2, \ldots,Q_i,\ldots\label{eq:alphai>1}
\end{equation}
or a sequence
\begin{equation}
\mathbf{Q}=Q_1,Q_2,\ldots,Q_{i-1},Q_i^{(j)},\qquad j \in \N.\label{eq:alphai=1}
\end{equation}
In the case when $\mathbf{Q}$ is given by (\ref{eq:alphai>1}), it is a
system of approximate roots, whether or not $A$ is complete. In the
case (\ref{eq:alphai=1}) assume, in addition, that the ring $A$ is
$\mathfrak m$-adically complete. In that
case,
$$
Q_\infty=\lim\limits_{j\to\infty}Q_i^{(j)}
$$
is a well defined element of $A$ and $\mathbf{Q}\cup\{Q_\infty\}$ is a system of
approximate roots.
\end{rek}

We recall some basic facts from Zariski's theory of complete ideals
in regular two-dimen- sional local rings.

Let $(A,\mathfrak{m})$ be a regular 2-dimensional local ring, $x,y$ a
regular system of parameters and let $\nu$ be a valuation centered at
$A$.

\begin{deft} \label{def:blowup}
An ideal $\mathcal{I}$ in a normal ring $B$ is said to be integrally closed
or complete if it contains all the elements $z$ of $B$ satisfying a monic
equation of the form $$z^n+ a_{n-1} z^{n-1} + \cdots + a_0=0$$ where
$a_{n-i} \in \mathcal{I}^{n-i}$.

An ideal $\mathcal{I}$ in $A$ is said to be simple if it cannot be
factored in a non trivial way as a product of two other ideals.

A  \textbf{local blowing up of $A$ with respect to $\nu$ along $\mathfrak{m}$}
is the map $A \to A[\frac{y}{x}]_{\mathfrak{m}_1}$, where
$\mathfrak{m}_1$ is the center of $\nu$ in $A[\frac{y}{x}]$.

For an element $f \in A$, we have $x^{\nu_{\mathfrak{m}}(f)} |\ f$ in
$A[\frac{y}{x}]_{\mathfrak{m}}$. The \textbf{strict transform} of $f$ in
$A[\frac{y}{x}]_{\mathfrak{m}}$ is the element $x^{-\nu_{\mathfrak{m}}(f)}f$.
\end{deft}

\begin{rek}
Any $\nu$-ideal is a complete ideal.
\end{rek}

Now let $I$ be a simple $\mathfrak{m}$-primary $\nu$-ideal. Then

\noi (1) The set $$\mathfrak{m} = I_0 \supset I_1 \supset \cdots \supset
I_\ell = I$$  of simple $\nu$-ideals of $A$ containing $I$
is entirely determined by $I$ (it does not depend on
$\nu$).

\noi (2) Let $A \to A_1$ be the local blowing up with
respect to $\nu$ along $\mathfrak{m}$
and, for $i \geq 1$, let $I'_i$ be the transform of $I_i$ (that is,
$I'_i = x^{-\mu}I_iA_1$ with $\mu= ord_{\mathfrak{m}}I_i$). Then
$$\mathfrak{m}_1 = I'_1 \supset I'_2 \supset \cdots \supset
I'_{\ell-1} = I'$$ is the set .

\noi (3) Iterating this procedure $\ell$-times, we obtain a sequence
of local blowing ups \begin{equation} \label{suite}(A,\mathfrak{m})
  \to (A_1,\mathfrak{m}_1) \to
\cdots \to (A_\ell,\mathfrak{m}_\ell) \end{equation} such that the transform
$I^{(\ell)}$ of $I$ is $\mathfrak{m}_\ell$. For any $f \in A \setminus
I$, the strict transform of $f$ in $A_\ell$
is a unit of $A_\ell$.
\medskip

We recall the following general result from the theory of approximate
roots in regular 2-dimensional local rings (\cite{Spi1}).

 Let $A$ be a 2-dimensional regular local  ring, $\nu$ a valuation on
 $A$. Now let $Q_k$, $k = 1,\dots,g+1$ be the approximate roots of $\nu$
 such that $Q_1,\dots,Q_g \notin I$ and $Q_{g+1} \in I$. Each $I_i$ is
 generated by the
 generalized monomials $\prod Q_j^{\gamma_j}$, $\gamma_j \in \N$, such
 that $\sum \gamma_j \beta_j \geq \nu(I_i)$.

\begin{prop} \label{eclat}
There exist natural numbers  $\ell_1 < \ell_2 < \cdots < \ell_g \leq
\ell$ and a regular system of parameters $x_{\ell_i},y_{\ell_i}$ for
each $i \in \{1,\ldots,g\}$ having the
following properties :

(1) $x_{\ell_i}$ is a monomial of the form
$\dsp{\prod\limits_{j=1}^{i-1}Q_j^{\gamma_j}}$, $\gamma_j \in \N$,

(2) $y_{\ell_i}$ is the strict transform of $Q_i$ in $A_{\ell_i}$,

(3) $Q_1,\dots,Q_{i-1}$ are monomials in $x_{\ell_i},y_{\ell_i}$ times a unit of
$(A_{\ell_i})_{(x_{\ell_i},y_{\ell_i})}$.
\end{prop}

For $\alpha, \beta \in \sper A$, let $Q_1,\dots,Q_s$ be the
approximate roots common to $\alpha$ and $\beta$.

\begin{cor}
If $i \leq s$, both $\nu_\alpha$ and $\nu_\beta$ are centered at
$(x_{\ell_i},y_{\ell_i})$.
\end{cor}

Let $A$ be a 2-dimensional regular local  ring, $\nu$ a valuation on
 $A$. Keep all the above notations.
\smallskip

\noi Convention : below, we adopt the convention that $\alpha_1=1$.

\begin{lem}
For $i\geq 3$, $\nu_{\mathfrak{m}}(Q_i)= \prod_{j=1}^{i-1}\alpha_j$.
\end{lem}

\noi Proof : Let $i=3$, then we can write $Q_3= y^{\alpha_2} + \sum
c_{r s}x^ry^s$ where $c_{r s}$ is a unit in $A$, with $\nu(x^ry^s) \geq
\alpha_2\nu(y)$. As $\nu(y) \geq \nu(x)$, this implies
$\nu_{\mathfrak{m}}(Q_3) = \alpha_2$.

Recall (cf. (\ref{eqn:appr2-1})) that
$$Q_{i+1} = Q_i^{\alpha_i} +
\sum_{\ell=1}^d\left(\sum_{\gamma=(\gamma_1,\ldots,\gamma_{i-1})}
  b_{\ell \gamma} \prod_{r=1}^{i-1} Q_r^{\gamma_r} \right)
Q_i^{\alpha'_i (d-\ell)}.$$

\noi Now to prove the lemma, it suffices to prove that
\begin{equation} \label{eq:suffit}
\alpha_i \nu_{\mathfrak{m}}(Q_i) \leq
\nu_{\mathfrak{m}}\left(\left(\prod_{r=1}^{i-1}
    Q_r^{\gamma_r}\right)Q_i^{\alpha'_i (d-\ell)}\right)
\end{equation}  for all
$\ell$ and $\gamma$ such that $b_{\ell \gamma} \neq 0$.

\noi First remark that, according to the inequalities (\ref{eqn:appr2-4})
and (\ref{eqn:appr2-6}), we
deduce by an easy induction on $i-\ell$
that
\begin{equation} \label{relat-a-b}
\frac{\beta_i}{\prod_{q=\ell}^{i-1} \alpha_q} \geq \beta_\ell.
\end{equation}

\noi We have $\alpha_i \beta_i = \sum_{j=1}^{i-1}
\gamma_j \beta_j + \alpha'_i (d-\ell) \beta_i$, so $$\alpha'_i \ell \beta_i =
\sum_{j=1}^{i-1} \gamma_j \beta_j \leq \sum_{j=1}^{i-1}
\gamma_j \frac{\beta_i}{\prod_{q=j}^{i-1}\alpha_q}$$ by
(\ref{relat-a-b}).

Dividing both sides by $\frac{\beta_i}{\prod_{q=1}^{i-1}\alpha_q}$, we
get
\begin{equation} \label{eq:suffit2}
\alpha'_i \ell \prod_{q=1}^{i-1}\alpha_q \leq \sum_{j=1}^{i-1}
\gamma_j \prod_{q=1}^{j-1}\alpha_q.
\end{equation}

By the induction assumption, the left hand side equals
$\nu_{\mathfrak{m}}(Q_i^{\alpha'_i \ell})$ and the right hand side equals
$\nu_{\mathfrak{m}}(\prod_{j=1}^{i-1}Q_j^{\gamma_j})$. Therefore
inequality (\ref{eq:suffit}) follows from inequality
(\ref{eq:suffit2}).$\Box$
\medskip

 In what follows, we study standard monomials $\dsp{\prod_{j=1}^i
  Q_j^{\gamma_j}}$, with $i < s$, that is monomials such that $0
\leq \gamma_j < \alpha_j$ for $j \in \{2,\ldots,i\}$.

\begin{cor} \label{cor-ppal1}
Consider two standard monomials $\prod_{j=1}^i Q_j^{\gamma_j}$ and
$\prod_{j=1}^i Q_j^{\gamma'_j}$ such that \\
$(\gamma_i,\gamma_{i-1},\ldots,\gamma_1) <_{lex}
(\gamma'_i,\gamma'_{i-1},\ldots,\gamma'_1)$ and having the same
$\nu$-value. We have
$$
\nu_{\mathfrak{m}}\left(\prod_{j=1}^i Q_j^{\gamma_j}\right) >
\nu_{\mathfrak{m}}\left(\prod_{j=1}^i Q_j^{\gamma'_j}\right).
$$
\end{cor}

Let $n =
\nu_{ \mathfrak{m}}(Q_3)$; note that $\alpha_2=n$. Moreover
$[k_2:k]\ |\ n$ and $[k_2:k]=n$ if and only if $\beta_1\ |\ \beta_2$.

\begin{cor}  \label{cor-ppal2} Consider two standard monomials
  $\prod\limits_{j=1}^i Q_j^{\gamma_j}$ and
$\prod\limits_{j=1}^i Q_j^{\gamma'_j}$, with $3 \leq i <s$, such that
$(\gamma_i,\gamma_{i-1},\ldots,\gamma_3) <_{lex}
(\gamma'_i,\gamma'_{i-1},\ldots,\gamma'_3)$.
We have
$$
\nu_{\mathfrak{m}}\left(\prod_{j=3}^i Q_j^{\gamma_j}\right) \leq
\nu_{\mathfrak{m}}\left(\prod_{j=3}^i Q_j^{\gamma'_j}\right)-n.$$
\end{cor}

\noi Proof : Let $j \geq 3$ be the greatest integer such that $\gamma_j <
\gamma'_j$. We have
\begin{eqnarray*} \nu_{\mathfrak{m}}\left(\prod_{j=3}^i Q_j^{\gamma'_j}\right)
-\nu_{\mathfrak{m}}\left(\prod_{j=3}^i Q_j^{\gamma_j}\right) &= &
\sum_{\ell=3}^j \gamma'_\ell \prod_{q=1}^{\ell-1}\alpha_q -
\sum_{\ell=3}^j \gamma_\ell \prod_{q=1}^{\ell-1}\alpha_q\\ & =&
(\gamma'_j-\gamma_j)\prod_{q=1}^{j-1}\alpha_q +
\sum_{\ell=3}^{j-1}(\gamma'_\ell-\gamma_\ell)
\prod_{q=1}^{\ell-1}\alpha_q.
\end{eqnarray*}

\noi \textbf{Claim:} For $j \geq 4$ and $c_\ell < \alpha_\ell$, we have
\begin{equation}\label{eq:claim}
\sum_{\ell=3}^{j-1}c_\ell
\prod_{q=1}^{\ell-1}\alpha_q < \prod_{q=1}^{j-1}\alpha_q.
\end{equation}

\noi Proof of Claim: By induction on $j$. For $j=4$, the inequality is
immediate. Assume the Claim is true for $j-1$. The left hand side of
(\ref{eq:claim}) can be rewritten as
$$
\sum_{\ell=3}^{j-1}c_\ell\prod_{q=1}^{\ell-1}\alpha_q =
\sum_{\ell=3}^{j-2}c_\ell\prod_{q=1}^{\ell-1}\alpha_q +
c_{j-1}\prod_{q=1}^{j-2 }\alpha_q < \prod_{q=1}^{j-2}\alpha_q +
c_{j-1}\prod_{q=1}^{j-2}\alpha_q \leq \prod_{q=1}^{j-1}\alpha_q.
$$

The Claim is proved.

The monomials being standard,
 $0 \leq \gamma_\ell, \gamma'_\ell < \alpha_\ell$, so
 $\gamma'_\ell-\gamma_\ell > -\alpha_\ell$ and applying the Claim, we deduce
that
$$
\sum_{\ell=3}^{j-1}(\gamma'_\ell-\gamma_\ell)
\prod_{q=1}^{\ell-1}\alpha_q > - \prod_{q=1}^{j-1}\alpha_q.
$$

Since $\gamma'_j-\gamma_j \geq 1$, we
get
$$
(\gamma'_j-\gamma_j)\prod_{q=1}^{j-1}\alpha_q +
\sum_{\ell=3}^{j-1}(\gamma'_\ell-\gamma_\ell)
\prod_{q=1}^{\ell-1}\alpha_q > 0.
$$
Each term being an integer divisible by $\alpha_2$, the above
expression is greater or equal to $\alpha_2=n$. $\Box$

\subsection{Real geometric surfaces and their blowings up}\label{realsurfaces}

Let $A$ be a ring and $U$ an open subset of
$\sper(A)$. Let $S_U$ denote the multiplicative set
$$
S_U=\{g\in A\ |\ g(\alpha)\ne0\text{ for all }\alpha\in U\}.
$$
Let $A_U=A_{S_U}$. We have a natural ring homomorphism
$$
\rho_U:A_U \to\prod_{\alpha \in U} A(\alpha).
$$
Define the ring $O_U$ to be
the ring of all maps
$$
f:U \to \coprod_{\alpha \in U}
  A(\alpha)
$$
satisfying the following conditions :

(1) $\forall \alpha \in U$, $f(\alpha) \in A(\alpha)$;

(2) there exists an open covering
\begin{equation}
U=\bigcup\limits_{i\in\Lambda} U_i\label{eq:covering}
\end{equation}
and, for each $i$, an
  element  $f_i \in A_{U_i}$ such that $\forall \beta \in U_i$, we have
  $\rho_{U_i}(f_i)_\beta =f(\beta)$.

The functor which sends $U$ to $O_U$ makes $\sper(A)$ into a locally
ringed space which we will call an \textit{affine real geometric
  space}. This notion is inspired by the notion of real closed spaces defined by
Niels Schwartz (\cite{Sch}).
\medskip

\noi From now till the end of this section we will assume that all our rings are
integral domains.
\begin{rek}
Note that $\iota:A_U \hookrightarrow O_U$ and, if $U$ is connected, this
inclusion becomes an equality. Indeed, consider an element $f\in O_U$, the
open covering (\ref{eq:covering}) and the local representatives $f_i\in
A_{U_i}$ of $f$ as above. Let $K$ denote the common field of fractions of $A$
and all of the $A_U$. Finding an element $g\in A_U$ such that $\iota(g)=f$
amounts to proving that for each $i,j\in\Lambda$ we have
\begin{equation}
f_i=f_j,\label{eq:fi=fj}
\end{equation}
viewed as elements of $K$. By connectedness of $U$, it is sufficient to prove
(\ref{eq:fi=fj}) under the assumption that $U_i\cap U_j\ne\emptyset$. Take a
non-empty basic open subset $V\subset U_i\cap U_j$, defined by finitely many
inequalities $V=\{\alpha\in\sper\ A\ |\ g_1(\alpha)>0,\dots,g_s(\alpha)>0\}$.
Since $V\ne\emptyset$, Propositions 4.3.8 and 4.4.1 (Formal Positivestellensatz)
of \cite{BCR} imply that $V$ contains a point $\alpha$ such that $\mathfrak
p_\alpha=(0)$. Then $A(\alpha)=K$, so the equality
$\rho_{U_i}(f_i)_\alpha=f(\alpha)=\rho_{U_j}(f_j)_\alpha\in A(\alpha)$ implies
that $f_i=f_j$ in $A(\alpha)=K$, as desired.
\end{rek}

\noi \textbf{Notation} To simplify the notation, we will write $A_i$ instead of
$A_{U_i}$.

\begin{deft}
A real geometric space is a locally ringed space $(X,O_X)$ which
admits an open covering $\dsp{X=\bigcup_{i=1}^s \sper(A_i)}$
such that
each $(U_i,O_X|U_i)$ is isomorphic (as locally ringed space) to an
affine real geometric space.
\end{deft}

\begin{deft}
A real geometric surface is a real geometric space $X$ where all $A_i$
can be chosen to be regular 2-dimensional noetherian rings.
\end{deft}

Let $k$ be a field and $z$ an independent variable. Let $A$ be a
regular two-dimensional ring, $x,y$ elements of $A$, $\mathfrak p$ a
maximal ideal of $A$ of height 2, containing $x$. Suppose given an
isomorphism $\iota:\frac A{(x)}\tilde\rightarrow k[z]_\theta$ such that
$y$ mod $(x)$ is sent to $z$ and $\theta $ is a non-zero polynomial in $z$. Let
$g=z^d+\bar a_1z^{d-1}+\dots+\bar a_d$ denote the
monic generator of the ideal $\iota\left(\frac{\mathfrak
    p}{(x)}\right)$. Let $a_i$ be an element of the coset $\iota^{-1}(\bar
a_i)$. Then $\left(x,y^d+a_1y^{d-1}+\dots+a_d\right)$ is a set of
generators of $\mathfrak p$; it induces a regular system of parameters
of $A_{\mathfrak p}$.
\begin{deft} The pair $\left(x,y^d+a_1y^{d-1}+\dots+a_d\right)$ will
  be called a \textbf{privileged system of parameters of }
  $A_{\mathfrak p}$ with respect to the ordered pair $(x,y)$.
\end{deft}

\begin{deft} \label{def-mrgs}
A marked real geometric surface is a real geometric surface $X$ together
with the following additional data:

(1) A finite covering $\dsp{X=\bigcup_{i=1}^s \sper(A_i)}$ where each $A_i$
is a regular 2-dimensional noetherian ring.

(2) For each $i$, a pair of elements $x_i,y_i\in A_i$ and a field
$k_i$, which admits a total ordering.

(3) A subset $\Delta_i\subset\sper\ A_i$, called the
\textbf{privileged} subset of $\sper\ A_i$. Let $z,w$ be independent
variables. We require one of the following to hold:

\qquad(a) There exists an irreducible polynomial $h\in k_i[w]$ and a
homomorphism
$$
\iota:A_i\rightarrow\frac{k_i[z,w]_{\theta_z\theta_w}}{(zh)},
$$
where $\theta_z \in k_i[z,w] \setminus (z,h),\ \theta_w \in
k_i[w] \setminus (h)$,
which maps $x_i$ to $z\mod(zh)$, $y_i$ to $w\mod(zh)$ such that
$\Delta_i$ is the set of points of $\sper\ A_i$ defined by the
vanishing of all the elements of $Ker\ \iota$ (in particular,
$\dsp{\Delta_i\cong\sper\frac{k_i[z,w]_{\theta_z\theta_w}}{(zh)}}$);

\qquad(b) $\Delta_i=\{x_i=0\}$; there is an isomorphism
$\dsp{\iota:\frac{A_i}{(x_i)}\rightarrow k_i[w]_{\theta_w}}$, where
$\theta_w$ is a non-zero polynomial in $k_i[w]$, which sends
$y_i\mod(x_i)$ to $w$; in particular, $\Delta_i\cong\sper\ k_i[w]_{\theta_w}$;

\qquad(c) $\Delta_i=\{x_i=y_i=0\}$; we have $\dsp{\frac{A_i}{(x_i,y_i)}\cong
k_i}$; in particular, $\Delta_i\cong\sper\  k_i$.

(4) For each $i$ and each $\alpha\in\{x_i=0\}\subset\Delta_i$ with $ht\
\mathfrak{p}_\alpha=2$, a regular system of  parameters of
$(A_i)_{\mathfrak{p_\alpha}}$, privileged with respect to $(x_i,h)$ in case (a)
and with respect to $(x_i,y_i)$ in case (b).

(5) In case (a), for each $i$ and each $\alpha\in\{h=0\}\subset\Delta_i$ with
$ht\ \mathfrak{p}_\alpha=2$, a regular system of parameters of
$(A_i)_{\mathfrak{p_\alpha}}$, privileged with respect to $(h,x_i)$.
\end{deft}

\begin{rek}
Let $A$ be a regular 2-dimensional ring, $\mathfrak{m}$ a maximal
ideal of $A$ and $(x,y)$ a regular system of parameters of
$A_{\mathfrak{m}}$. Then $\sper\ A$ is a marked real geometric
surface.
\end{rek}

We now define the notion of blowing up of a real marked geometric
surface. Let $X = \bigcup_i \sper\ A_i$ be a marked real geometric
surface and take a point $\delta \in X$. Assume that $\delta$ belongs
to the privileged set and ht$(\mathfrak{p}_{\delta,i})=2$ in every affine
chart $\sper\ A_i$ containing $\delta$. We want to define the
\textbf{blowing up of $X$ along $\delta$}. First consider the case
$X=\sper\ A$. Let $x,y \in A$ and $k$ be the pair of elements and the
field appearing in the definition of marked real geometric surface.
\medskip

Let $(u,v)$ be the privileged system of regular parameters of
$A_{\mathfrak{p}_\delta}$ given by the definition. It follows from
definition that $(u,v)= (x,y)$ in Case (c), $u=x$ in Case
(a) provided $\delta \in \{x=0\}$ as well as in Case (b), and $u=h$ in
Case (a) if $\delta \in \{h=0\}\setminus \{x=0\}$.
\medskip

A \textbf{blowing up} of $\sper\ A$ along $\mathfrak{p}_\delta$ (or, by abuse of
language, blowing up along $\delta$) is the
marked real geometric surface $X'$ defined as follows. As a
topological space, we put $X'=\sper\
  A'_1\cup\sper\ A'_2$, where $A'_1=A\left[\frac{v}{u}\right]$,
$A'_2=A\left[\frac{u}{v}\right]$ and
$$
\sper\ A'_1\cap\sper\ A'_2=\sper\
A'_1\setminus\left\{\frac{v}{u}=0\right\}=\sper\
A'_2\setminus\left\{\frac{u}{v}=0\right\}.
$$

\noi We have a natural surjective morphism $\pi:X' \rightarrow \sper\
A$.

To define a structure of marked real geometric surface on $X'$, we let
the two elements required in Definition \ref{def-mrgs} (2) be $x'_1=u,
y'_1=\frac vu \in A'_1$ for $\sper\ A'_1$ and $x'_2=v,y'_2= \frac uv
\in A'_2$ for $\sper\ A'_2$. Below, for $q \in \{1,2\}$, we denote the
privileged set of $A'_q$ by $\Delta'_q$ and the field required in the
Definition \ref{def-mrgs} (2) for $\sper\ A'_q$ by $k'_q$. We now
define $\Delta'_q$ and $k'_q$ in the different cases.
\smallskip

\quad $\bullet$ If Case (c) holds for $\sper\ A$: let $k'_q=k$, for $q
\in \{1,2\}$. For $\sper\ A'_1$ the privileged
set is $\Delta'_1 = \{x'_1=0\}$. The existence of a privileged regular
system of parameters required by the Definition \ref{def-mrgs} comes
from the isomorphism $\dsp{\frac{A'_1}{(x'_1)} \cong k[y'_1]}$. For
$\sper\ A'_2$ the situation is entirely analogous.
\smallskip

\quad $\bullet$ If Case (b) holds for $\sper\ A$ : let
$k'_1=\kappa(\mathfrak{p}_\delta)$ and
$\Delta'_1= \{x'_1=0\}$. The existence of a
privileged regular system of parameters required by the Definition
\ref{def-mrgs} comes from $\dsp{\frac{A'_1}{(x'_1)} \tilde\to
  k'_1[y'_1]}$.

 Let $k'_2=k$ and $\Delta'_2= \{x'_2=0\} \cup \{y'_2 =0\}$. By the
 definition of privileged regular system of parameters of
 $A_{\mathfrak{p}_\delta}$, there is an
 irreducible polynomial $v_w \in k[w]$, relatively prime to $\theta_w$
 such that $\iota(\frac{\mathfrak{p}_\delta}{(x)})=(v_w)$. The
 existence of a privileged regular system of parameters at any point
 of $\Delta'_2$, required by the Definition \ref{def-mrgs}, comes from
 $\dsp{\frac{A'_2}{(x'_2y'_2)} \tilde\to
   \frac{k[w,y'_2]_{\theta_w}}{(v_w y'_2)}}$.
\smallskip

\quad $\bullet$ If Case (a) holds, there are three cases to consider :

\qquad (i) $\delta \in \{x=0\}\setminus\{h=0\}$, $\Delta'_q,k'_q$,
$q=1,2$, are given by the same formulas as in Case (b). Let
$A'_3=A_v$. The structure of marked real geometric surface on $\sper\
A'_3$ is induced from that of $\sper\ A$. We have
$k'_3=k$ and $\Delta'_3 =\{x=0\}\cup \{h=0\}$ and $\dsp{\frac{A_v}{(xh)}
\tilde\to \frac{k[z,w]_{v_w\theta_w\theta_z}}{(zh)}}$.

\qquad (ii) $\delta \in \{h=0\}\setminus\{x=0\}$, let
$k'_1=\kappa(\mathfrak{p}_\delta)$ and $k'_2=\frac{k[w]}{(h)}$,
$\Delta'_1=\{x'_1=0\}$, $\Delta'_2=\{x'_2=0\}\cup \{y'_2=0\}$. By the
 definition of privileged regular system of parameters of
 $A_{\mathfrak{p}_\delta}$, there is a polynomial $v_z \in k[z,w]$,
 such that $\dsp{\iota(\frac{\mathfrak{p}_\delta}{(xh)}) =
   \frac{(h,v_z)k[z,w]_{\theta_z\theta_w}}{(xh)}}$. The existence of a
   privileged regular system of parameters comes from the isomorphisms
   $\dsp{\frac{A'_1}{(x'_1)} \tilde\to k'_1[y'_1]}$ and
   $\dsp{\frac{A'_2}{(x'_2y'_2)} \tilde\to
  \frac{k'_2[z,y'_2]_{\theta_z}}{(v_zy'_2)}}$.

Let $A'_3=A_v$. The structure of marked real geometric surface on $\sper\
A'_3$ is induced from that of $\sper\ A$. We have $k'_3=k$, $\Delta'_3
=\{x=0\}\cup \{h=0\}$ and $\dsp{\frac{A_v}{(xh)}
\tilde\to \frac{k[z,w]_{v_z\theta_w\theta_z}}{(zh)}}$.

\qquad (iii) $\delta= \{h=0\} \cap \{x=0\}$, recall that
$u=x,v=h$. Let $k'_1 = \kappa(\mathfrak{p}_\delta)$ and $k'_2=k$. Let
$\Delta'_q = \{x'_q=0\}\cup \{y'_q=0\}$, $q=1,2$.  The
existence of a privileged regular system of parameters comes from the
isomorphisms $\dsp\frac{A'_1}{(x'_1y'_1)} \tilde\to
\frac{k'_1[z,y'_1]_{\theta_z}}{(zy'_1)}$ and
$\dsp\frac{A'_2}{(x'_2y'_2)} \tilde\to \frac{k[w,y'_2]_{\theta_w}}{(hy'_2)}$
(recall that in this case $h=x'_2$).
\medskip

We then define the real marked geometric surface $X'$ to be
$X'=\bigcup\limits_{i=1}^p \sper\,A'_i$ where $p=2$ in cases (b), (c) and (a)
(iii) and
$p=3$ in cases (a) (i) and (ii).
\medskip

\begin{rek}
Note that
$\sper\,A'_3 \subset \sper\,A'_i,\ i=1,2$; but, in the
applications, we need to have the set $\Delta'_3$ defined by fixed
elements $x'_3,y'_3$.
\end{rek}

Let $X=\bigcup\limits_{i=1}^s \sper\, A_i$ be a marked real geometric surface
and $\delta \in X$ belonging to the privileged set and supported in a height 2
ideal $\mathfrak p_{\delta,i}$ in some affine chart $\sper\ A_i$.

If $\delta\notin\sper\, A_i$, let $X'_i=\sper\ A_i$ with the identity
map $X'_i\rightarrow \sper\, A_i$. If $\delta\in\sper\,A_i$,
let $X'_i\rightarrow \sper\,A_i$ be the blowing up of $\sper\,A_i$ along
$\mathfrak{p}_{\delta,i}$. Let $(u,v)$ be the regular system of parameters of
$(A_i)_{\mathfrak{p}_{\delta,i}}$ given by the definition of real marked
geometric
surface. We have $X'_i=\bigcup_{j=1}^p \sper\,A'_{ji}$ where $p=2$ or $3$ as
above.

The marked real geometric surfaces $X'_1,\dots,X'_s$ and the maps
$X'_i\rightarrow \sper\ A_i$ glue together in a natural way to give a
marked real geometric surface
$X'=\bigcup\limits_{i=1}^sX'_i$ and the map $X'\rightarrow X$.

\begin{deft} \label{def:center}
We call $X'$ \textbf{the blowing up of }$X$ along $\delta$ or
the point blowing up of $X$ along $\delta$. The point $\delta$ is called the
\textbf{center} of this blowing up. If $X=\sper\ A$, the blowing up of
$X$ along $\delta$ depends only on the ideal
$\mathfrak{p}_\delta$ and not on the ordering $\leq_\delta$, so we may
speak also about blowing up along $\mathfrak{p}_\delta$.
\end{deft}

\begin{deft} Let $\alpha, \delta$ be two distinct points of the real
  marked surface $\sper\ A$ with
  $$
  ht(\mathfrak{p}_\delta)=2.
  $$
  Let $\pi:X' \to \sper\ A$ be a blowing up along $\delta$.  Let
  $(u,v)$ be the given privileged system of parameters at
  $\delta$. Since $\alpha \neq \delta$, $\{u,v\} \not\subset
  \mathfrak{p}_\alpha$. If $u \notin \mathfrak{p}_\alpha$, the
  \textbf{strict transform} $\alpha'$ of $\alpha$ is defined as
  follows. Let $\mathfrak{p}_{\alpha'}$ be the strict
  transform of $\mathfrak{p}_{\alpha}$ in $A'_1$ and $\leq_{\alpha'}$ be the
  order of $\kappa(\mathfrak{p}_{\alpha'})$ induced by $\leq_\alpha$
  via the natural isomorphism $\kappa(\mathfrak{p}_{\alpha}) \cong
  \kappa(\mathfrak{p}_{\alpha'})$. If $v \notin \mathfrak{p}_\alpha$,
  $\alpha' \in \sper\ A'_2$ is defined similarly.
\end{deft}

On the way to prove the connectedness of $C$ of (\ref{ensC}), we will
now prove a preliminary result on connectedness of a certain type of
subsets (intervals) of the exceptional divisor on a suitable blowing
up of $\sper\ A$.

\begin{rek} \label{rek:ordre}
Fix an order on $k$.
Let $D$ be the set of points $\delta \in \sper(k[z])$ which induce the
given order on $k$. Given two points $\delta_1
\neq \delta_2 \in D$ such that $ht(\mathfrak{p}_{\delta_i})=1$, we view
$\delta_1,\delta_2$ as elements of
the real closure $\sur{k}$ of $k$ with respect to the given order. We may speak
  about the \textbf{interval} $(\delta_1,\delta_2)=\{\delta \in D\ |\ \delta_1
  < z(\delta) < \delta_2 \}$. If $\mathfrak p_\delta=(0)$, we
  compare $\delta_i$ and $z(\delta)$ via the natural embeddings
  $k[z](\delta) \hookrightarrow \sur{k(z)}$ and $\sur{k} \subset
  \sur{k(z)}$.

Now, let $\mathfrak m$ be an ideal of $A$ with $ht\ \mathfrak m=2$ and $\frac
A{\mathfrak m}=k$. Given a blowing up along $\mathfrak{m}$ as above, consider
the open set $\sper(A[\frac{y}{x}])$. The set of points $\delta \in
\sper(A[\frac{y}{x}])$ such that $x(\delta)=0$ and which induce the
given order on $k$ is homeomorphic to $D$.

Finally, let $X$ be a real algebraic surface such that $D \subset
\sper\ k[z] \subset X$. Let $+\infty$ denote the point of $D$ with
support $(0)$ such that $z(+\infty)>c$ for all $c \in k$. Let $\sur\infty$ be
the closed point of $X$ such that $\sur{\infty} \in \sur{\{+\infty\}}$. Assume
there is an open set $\sper\ A_i \subset X$ such that
$\mathfrak{p}_{\sur{\infty}}$ in $A_i$ has height 2. We extend the above
notion of interval to include the case when $\delta_2=\sur{\infty}$
with the obvious meaning assigned to $[\delta_1,\sur{\infty}] =
\bigcup\limits_{\delta > \delta_1}[\delta_1,\delta] \cup \{\sur{\infty}\},\
(\delta_1,\sur{\infty}),...$. Similarly, we may take a closed point
$-\sur{\infty} \in \sur{\{-\infty\}}$. As points of $X$, we have
$\sur{\infty}=-\sur{\infty}$. However, our ordering on $D$ provides us
with a well defined notion of intervals of the
form $(-\sur{\infty},\delta_1),\ [-\sur{\infty},\delta_1)$ and so on.
\end{rek}

\begin{lem} \label{inter-con}
Let $D$ be as in the remark before and $\delta_1 < \delta_2 \in
D$ such that $ht(\mathfrak{p}_{\delta_i})=1$. The closed interval $[\delta_1,\delta_2]$, the semi-open interval
$[\delta_1,\delta_2)$ and the open interval $(\delta_1,\delta_2)$ are
connected.
\end{lem}

\noi Proof : We will prove it for the open case, the closed and the semi-open being
similar. Let $k \hookrightarrow \sur{k}$ be the inclusion of $k$ into
its real closure determined by the given order. This map corresponds
to a morphism $\sper(\sur{k}[z]) \to \sper(k[z])$ which induces a
homeomorphism between $D$ and $\sper(\sur{k}[z])$ sending
$(\delta_1,\delta_2)$ to an interval
$(\sur{\delta}_1,\sur{\delta}_2)$ where
$\sur{\delta}_1,\sur{\delta}_2 \in \sur{k}$. It is
well-known and easy to prove that such an interval is connected - in
the spectral topology (see for instance \cite{BCR}). $\Box$
\medskip

\begin{rek}
Let $\theta \in k[z]$ be a non-zero polynomial. We have natural
homeomorphisms $\sper\ k[z]_\theta \tilde\to \sper\ k[z] \setminus
\{\alpha_1,\ldots,\alpha_t\}$ and $\lambda: D\cap \sper\ k[z]_\theta\tilde\to D\setminus
\{\alpha_1,\ldots,\alpha_t\}$
where $\{\alpha_1,\ldots,\alpha_t\}$ is
the set of points $\alpha_i \in \sper\ k[z]$ such that $\theta \in
\mathfrak{p}_{\alpha_i}$. Let $\delta_1,\delta_2 \in D$ be as
above. Assume that $\alpha_i \notin (\delta_1,\delta_2)$ for all $i
\in \{1,\ldots,t\}$. Then $\lambda((\delta_1,\delta_2))$ is connected
in $D \setminus \{\alpha_1,\ldots,\alpha_t\}$.
\end{rek}

\begin{deft}
Let $\mathfrak{m}$ be a maximal ideal of $A$ of height 2. Let $X' \to
\sper\ A$ be the blowing up along $\mathfrak{m}$. Let
$\mathcal{E}=\{\epsilon \in \sper\ A\ | \ \mathfrak{p}_\epsilon =
\mathfrak{m}\}$. The sets $\pi^{-1}(\epsilon),\ \epsilon \in
\mathcal{E}$ are called the \textbf{components} of
$\pi^{-1}(\mathfrak{m})$.
\end{deft}
Let $(A,\mathfrak{m},k)$ be a regular 2-dimensional local ring and
$(x,y)$ a regular system of parameters. Now consider a sequence
\begin{equation} \label{def:sequence}
 X_t \stackrel{\pi_{t-1}}{\to} \cdots \stackrel{\pi_1}{\to} X_1
 \stackrel{\pi_0}{\to}
\sper\ A
\end{equation}
of point blowings up where the first blowing up $\pi_0:X_1 \to
\sper\ A$ is the blowing up along $\mathfrak{m}$.
\medskip

Fix a point $\epsilon \in \sper\ A$ such that $\mathfrak{p}_{\epsilon}=
\mathfrak{m}$ - this is equivalent to fixing a total
ordering on $k$. For $q \in \{0,\ldots,t-1\}$, let $\eta_q \in X_q$ be
the closed point, compatible with the given order, such that $\pi_q$
is a blowing up along $\eta_q$.

For $i \in\{1,\ldots,t\}$, let $X_i =
\bigcup\limits_{j=1}^{s_i}\sper\ A_{ji}$ be the open affine covering in the
definition of marked real geometric surface.

Let $\rho_i=\pi_0\circ \ldots \circ \pi_{i-1}: X_i \to \sper\ A$.
\medskip

\begin{rek}
The real geometric space $\rho_i^{-1}(\mathfrak{m})$ has the form
$\rho_i^{-1}(\mathfrak{m})= \bigcup\limits_\ell \sper\ B_{i\ell}$ with $B_{i\ell} \cong
k_{i\ell}[z_{i\ell}]$ where $k_{i\ell}$ is a finite algebraic
extension of $k$ and $z_{i\ell}$ is an independant variable.
\end{rek}

\begin{deft}
A subset $E \subset \rho_i^{-1}(\epsilon)$ is a \textbf{component} of
$\rho_i^{-1}(\epsilon)$ if $E$ is either a component of
$\pi_{i-1}^{-1}(\eta_{i-1})$ or a strict transform of a
component of $\rho_{i-1}^{-1}(\epsilon)$ when $i>1$.
\end{deft}

\begin{deft}\label{def-max} Let $\rho_i: X_i \to \sper(A)$. Fix a component $E
  \subset \rho_i^{-1}(\epsilon)$.
Fix an index $j \in \{1,\ldots,s_i\}$. A $j$\textbf{-distinguished} point of
$E$ is a point $\delta \in E$ such that either $\delta \notin \sper\
A_{ji}$ or $\rho_i^{-1}(\{xy=0\}) \supset \{x'y'=0\}$ and
$x'(\delta)=y'(\delta)=0$ where $(x',y')
\in A_{ji}$ is the privileged regular system of parameters at $\delta$ (in
particular, the privileged set of $\sper\ A_{ji}$ is given by
$\{x'=0\}\cup\{y'=0\}$).

A $j$-maximal interval $I$ is a subset $I \subset E$ such that there
exist $j$-distinguished points $\delta_1,\delta_2 \in E$, $\delta_1
\neq \delta_2$, such that

(1) $I=[\delta_1,\delta_2]$ and $I$ is connected;

(2) There are no $j$-distinguished points in $I \setminus
\{\delta_1,\delta_2\}$.

\noi A maximal interval is an interval which is $j$-maximal for some $j$.
\end{deft}

\begin{rek}
Note that a $j$-maximal interval may contain a
$\tilde{\jmath}$-distinguished point, where $j \neq
\tilde{\jmath}$. This occurs if $[\delta_1,\delta_2]$ is a $j$-maximal
interval, $\delta \in (\delta_1,\delta_2)$ and $\exists \tilde{\jmath}
\in \{1,\ldots,s_i\},\ \tilde{\jmath} \neq j$, such that
$(\delta_1,\delta_2)\cap \sper\, A_{\tilde{\jmath}i} \neq \emptyset$ and
$\delta\,
\notin\, \sper\, A_{\tilde{\jmath}i}$.
\end{rek}

\begin{prop}\label{prop-max} Fix a component $E \subset \rho_i^{-1}(\epsilon)$
and a maximal interval $[\delta_1,\delta_2] \subset E$.
Take $q \in \{1~,2\}$. There exists $j \in \{1,\ldots,s_i\}$ such
that $[\delta_1,\delta_2]$ is $j$-maximal and letting $x_i,y_i \in A_{ji}$ be
the elements given by Definition \ref{def-mrgs} we have:

(1)$_i$ $[\delta_1,\delta_2]\setminus \{\delta_q\} \subset \sper\
A_{ji}$,

(2)$_i$ for all $\delta \in [\delta_1,\delta_2]\setminus \{\delta_q\}$ with $ht(\mathfrak p_\delta)=2$,
$x_i$ is a part of the given privileged regular system of parameters of
$(A_{ji})_{\mathfrak{p}_\delta}$,

(3)$_i$ $[\delta_1,\delta_2] \cap \sper\ A_{ji} = \left\{\eta \in \sper\ A_{ji}\ \left|\
 x_i(\eta)=0 \mbox{ and } \sur{\delta_1} \leq y_i(\eta) \leq
\sur{\delta_2}\right.\right\}$ where
$$
\sur{\delta_1}, \sur{\delta_2} \in\sur{k_{ji}} \cup \{-\sur{\infty},\sur{\infty} \},
$$
with the notation of Remark \ref{rek:ordre} and the proof of Lemma \ref{inter-con}.
\end{prop}

\noi Proof: First, let $i=1$. We have $\dsp{X_1=\sper\ A\left[\frac{y}{x}\right]
\cup \sper\ A\left[\frac{x}{y}\right]}$. Denote
$\dsp{A\left[\frac{y}{x}\right]}$ by $A_{11}$ and
$\dsp{A\left[\frac{x}{y}\right]}$ by $A_{21}$. Let $\dsp{x_1=x,
y_1=\frac{y}{x}}$.
Fixing the component $E$ is equivalent to fixing a total order on $k$; this data is already
given. We have
$$
E \cap \sper\ A[y_1] \subset \sper\ k[y_1].
$$
Let the notation be as in Remark \ref{rek:ordre} with $y_1$ playing the role of $z$.

There are exactly two maximal intervals $[0,\sur{\infty}]$ and $[-\sur{\infty},0]$.
Say, for example, $I=[0,\sur{\infty}]$, $q=2$, then $j=1$ satisfies
the conclusion of the Proposition. And similarly for the other
three cases.

Now take $i\ge2$ and suppose the result true for $i-1$. Let $\delta_{p,
  i-1}=\pi_{i-1}(\delta_p)$, $p=1,2$. Let $\eta_{i-1}$ be the center
of the blowing up $\pi_{i-1}$.
First, assume that
\begin{equation}
E \subset \pi_{i-1}^{-1}(\eta_{i-1}).\label{eq:excdivisor}
\end{equation}
Take $\tilde{\jmath} \in\{1,\dots,s_{i-1}\}$ such that $\eta_{i-1}$ belongs to
the privileged set of $\sper\ A_{\tilde{\jmath},i-1}$. Let $(u,v)$ be the given
privileged regular system of parameters at $\eta_{i-1}$. If $j$ is such that
$(\delta_1,\delta_2)\subset\sper\ A_{ji}$ then $A_{ji}$ is one of
$A_{\tilde{\jmath},i-1}[\frac uv]$ or
$A_{\tilde{\jmath},i-1}[\frac vu]$; pick one of these two possible choices $j$
such that $[\delta_1,\delta_2]$ is $j$-maximal.
In this case (1)$_i$ is equivalent to saying that
\begin{equation}
[\delta_1,\delta_2]\ne[-\sur\infty,\sur\infty].\label{eq:maxinterval}
\end{equation}
Now, if we had $[\delta_1,\delta_2]=[-\sur\infty,\sur\infty]$, the point
$x_i=y_i=0$ would be a distinguished point in $(\delta_1,\delta_2)$ (by
definition of distinguished point). This is a contradiction and (1)$_i$ is
proved in the case when (\ref{eq:excdivisor}). (2)$_i$ and (3)$_i$ of the
Proposition follow immediately from the definition of marked real geometric
surface.

>From now on, assume that
\begin{equation}
E \not\subset \pi_{i-1}^{-1}(\eta_{i-1}).\label{eq:notexcdivisor}
\end{equation}

Note that since $[\delta_1,\delta_2]$ is a maximal interval of $E$,
$[\delta_{1,i-1},\delta_{2,i-1}]$ is a maximal interval of
$\pi_{i-1}(E)$. So, by the induction hypothesis, the Proposition holds for
$[\delta_{1,i-1},\delta_{2,i-1}] \subset
\pi_{i-1}(E)$. Take $\tilde{\jmath} \in\{1,\dots,s_{i-1}\}$ which satisfies the
conclusion of the Proposition with $i$ replaced by $i-1$ (in particular,
$[\delta_{1,i-1},\delta_{2,i-1}]$ is $\tilde{\jmath}$-maximal).

If $\eta_{i-1} \notin\ \sper\ A_{\tilde{\jmath},i-1}$, take $j \in
\{1,\ldots,s_i\}$ such that $A_{ji}=A_{\tilde{\jmath},i-1} $. This $j$ satisfies
the conclusion of the Proposition.

Next assume that $\eta_{i-1} \in \sper\ A_{\tilde{\jmath},i-1}$. Take the
elements $u,v\in A_{\tilde{\jmath},i-1}$
which induce the privileged regular system of parameters at $\eta_{i-1}$, given
by the definition of marked real geometric surface.

If $\eta_{i-1} \notin \pi_{i-1}([\delta_1,\delta_2])$, take $j$ such that
$\dsp{A_{ji}=A_{\tilde{\jmath},i-1}}$.

>From now on, assume that $\eta_{i-1} \in
\pi_{i-1}\left([\delta_1,\delta_2]\right)\cap \sper\
A_{\tilde{\jmath},i-1}$. Then
$$
\eta_{i-1}\in\{\delta_{1,i-1},\delta_{2,i-1}\} :
$$
if not,
$\pi_{i-1}^{-1}(\eta_{i-1}) \cap (\delta_1,\delta_2)$ would be a
$j$-distinguished point in $(\delta_1,\delta_2)$, which is impossible. The
intersection is taken as subsets of the topological space $X_i$; if
$\eta_{i-1} \notin\{\delta_{1,i-1},\delta_{2,i-1}\}$, this intersection is not
empty and consists of a single point. Let $j\in
\{1,\ldots,s_i\}$ be such that $A_{ji}$ is one of $A_{\tilde{\jmath},i-1}[\frac
uv]$ or
$A_{\tilde{\jmath},i-1}[\frac vu]$; pick one of these two possible choices $j$
such that $[\delta_1,\delta_2]$ is $j$-maximal.

In all the cases the index $j$ chosen in this way satisfies the conclusion of
the Proposition.
$\Box$
\medskip

\subsection{A proof of the Pierce-Birkhoff conjecture for regular 2-dimen-
sional rings.}\label{proofPB}

Let $A$ be a 2-dimensional regular local  ring, $\nu$ a valuation on
 $A$.imensional ring. In this section, we prove that
$A$ is a Pierce-Birkhoff ring (\cite{Mad1}).
Our proof is based on Madden's unpublished preprint (\cite{Mad2}), but
there are some differences. Here, we have tried to present a proof
which should be a pattern for a general proof of the conjecture in any
dimension.
1
\begin{thm} Let $A$ be a 2-dimensional regular ring, then $A$ is a
  Pierce-Birkhoff ring.
\end{thm}

Actually, we prove that $A$ satisfies the Definable Connectedness
Conjecture and also, in the special case where $A$ is
excellent, the Connectedness Conjecture.
\medskip

We start with some results which do not assume that $A$ is excellent
and which are needed in the proof of both of the above
versions of the Connectedness Conjecture. Let $\alpha, \beta \in \sper\
A$. By Remark \ref{rek:spezial}, we may assume that neither of $\alpha,
\beta$ is a specialization of the other.
\medskip

There are two possibilities : either $ht(<\alpha,\beta>)
= 1$ or $ht(<\alpha,\beta>) = 2$.
\medskip

\subsubsection{The case of height 1.}

Let $\delta$ be the most general common specialization of $\alpha$ and
$\beta$ and let $\mathfrak{p} = \sqrt{<\alpha,\beta>}$ be the support
of $\delta$. Then $A_{\mathfrak{p}}$ is a discrete valuation ring; take an
element $t \in A$ whose image in $A_\mathfrak{p}$ is a regular parameter of
$A_\mathfrak{p}$. Since $ht(\mathfrak{p})=1$ and
neither of $\alpha, \beta$ is a specialization of the other, we have
$\mathfrak{p}_\alpha = \mathfrak{p}_\beta=(0)$. There are only two
orders on $A$ which induce the given order on $A/\mathfrak{p}$ : one
with $t>0$ and one with $t<0$. Since $\alpha \neq \beta$,
$<\alpha,\beta> = \mathfrak{p}$ : of course, any element $g$ of
$\mathfrak{p}$ can be written as $g=t^\gamma \frac{a}{b}$, $a,b \notin
\mathfrak{p}$. As $t \in <\alpha,\beta>$, if $\gamma \geq 2$,
$\nu_\alpha(g) = \nu_\beta(g) > \nu_\alpha(t)$ so $g \in
<\alpha,\beta>$ and if $\gamma=1$, $g$ changes sign between $\alpha$
and $\beta$, so again $g \in <\alpha,\beta>$.
\medskip

Now let $f_1,\ldots,f_r \notin <\alpha,\beta>= \mathfrak{p}$, so
$f_i(\delta) \neq 0$ for $i\in \{1,\dots,r\}$. As $\delta \in
\sur{\{\alpha\}}$ and $\delta \in \sur{\{\beta\}}$, we conclude that
$\alpha$ and $\beta$ belong to the same connected component of $\sper\
A \setminus \{f_1\cdots f_r =0\}$.

\subsubsection{The case of height 2.}\label{Ctoquadrant} Now assume

$ht(<\alpha,\beta>) = 2$,
that is $\mathfrak{m}= \sqrt{<\alpha,\beta>}$ is maximal. By
Proposition \ref{loca1}, replacing $A$ by $A_{\mathfrak{m}}$ does
not change the problem, so we may assume that $A$ is local with
maximal ideal $\mathfrak{m}$.

Let $g \in \N$ be such that $Q_1,\ldots,Q_g
\notin <\alpha,\beta>$, $Q_{g+1} \in <\alpha,\beta>$ be the approximate
roots common to $\nu_\alpha$ and $\nu_\beta$ as in section
\ref{subs-proof}.
\medskip

Let $(x,y)$ be a regular system of parameters of $A$ such that $\nu_\alpha(x) =
\nu_\alpha(\mathfrak{m})$ and $\nu_\beta(x) = \nu_\beta(\mathfrak{m})$.

Let $\pi:A \to A'$ be a local blowing up with respect to $\nu_\alpha$ and
denote by $k'$ the residue field of $A'$. Recall from (\cite{Zar}, Appendix 5)
that the \textbf{weak
  transform} $I' \subset A'$ of an ideal $I \subset A$ is defined by
$I'=x^{-a}IA'$ where $a = \nu_{\mathfrak{m}}(I)$.

\begin{prop}
We assume that $\pi$ is also a local blowing up with respect to
$\nu_\beta$. Let $\alpha'$ and $\beta'$ be the strict transforms of
$\alpha$ and $\beta$. Then the separating ideal $<\alpha',\beta'>$ is
equal to the weak transform of $<\alpha,\beta>$.
\end{prop}

\noi Proof : Since by hypothesis, $\alpha',\ \beta'$ are both centered
at a maximal ideal $\mathfrak{m}'$, we have $<\alpha,\beta>
\subsetneqq \mathfrak{m}$. In particular, $x \notin <\alpha,\beta>$,
hence $x$ does not change sign between $\alpha$ and $\beta$. Then $f \in
A$ changes sign between $\alpha$ and $\beta$ if and only if $x^{-a}f$
changes sign between $\alpha'$ and $\beta'$.

Since $<\alpha,\beta>$ is generated by elements changing sign between
$\alpha$ and $\beta$, its weak transform is generated by elements
which change sign between $\alpha'$ and $\beta'$; hence the weak
transform of $<\alpha,\beta>$ is contained in $<\alpha',\beta'>$.

To prove the opposite inclusion, let $I'=<\alpha',\beta'>$ and let $I$
be the inverse transform of $I'$, that is the unique complete ideal of
$A$ whose weak transform is $I'$ (\cite{Zar}, Appendix 5, p. 388). It
remains to prove that $I \subseteq <\alpha,\beta>$.

In order to do this, it suffices to find an element $z \in I$ which
changes sign between $\alpha$ and $\beta$ and such that $\nu_\alpha(z)
= \nu_\alpha(I)$.

Let $J_+$ be the greatest $\nu_\alpha$-ideal of $A'$ whose
$\nu_\alpha$-value is strictly greater than $\nu_\alpha(I)$. Note that
$\dsp{\frac{IA'}{J_+\cap IA'}}$ is a $k'$-vector space. Let
$b_1,\ldots,b_\ell$, $b_j=\prod_{r=1}^iQ_r^{\gamma_{jr}}$, where $i$
is the maximal index of the approximate roots $Q_s$ involved, be a set of
elements of $I$ which induces a basis of $\dsp{\frac{IA'}{J_+\cap IA'}}$, each
monomial being standard. Moreover, since $x$ divides $y$ in $A'$, if
$\nu_\alpha(x)=\nu_\alpha(y)$, we may assume $\gamma_{j2}=0$ for all $j$ and
$b_1$ is the unique monomial which maximizes the vector
$(\gamma_{i1},\gamma_{i-1,1}, \ldots, \gamma_{31})$ in the lexicographical
ordering.

Let $a=\nu_\mathfrak{m}(I)$. Let $\tilde{z} \in I'$ be such that
$\nu_\alpha(\tilde{z}) = \nu_\alpha(I')$ and $\tilde{z}$ changes sign between
$\alpha'$ and $\beta'$. Let $z^\dag=x^a\tilde{z}$. Then $z^\dag \in IA'$ and
$\nu_\alpha(z^\dag) = \nu_\alpha(IA') = \nu_\alpha(I)$.  Write
$z^\dag=\sum_{j=1}^\ell z_jb_j$. We may assume $z_1=1$. Denote by
$\sur{z}_j$ the image of $z_j$ in the residue field $k'$.

First, suppose $\nu_\alpha(x) < \nu_\alpha(y)$. Then $k'=k$. For each $j \in
\{1,\ldots,\ell\}$, let $w_j$ be a representative of $\sur{z}_j$ in $A$. Put
$z=\sum_{j=1}^\ell w_jb_j$ .

Next, suppose $\nu_\alpha(x)=\nu_\alpha(y)$, since $b_1$ is the unique
monomial which maximizes the vector $(\gamma_{i1},\gamma_{i-1,1}, \ldots,
\gamma_{31})$, by the corollary (\ref{cor-ppal2}), we have, for $j\geq 2$,
$$
\nu_{\mathfrak{m}}\left(\prod_{r=3}^iQ_r^{\gamma_{jr}}\right) \leq
\nu_{\mathfrak{m}}\left(\prod_{r=3}^iQ_r^{\gamma_{1r}}\right)-n \leq a-n.
$$
Write $z^\dag=b_1+\sum_{j=2}^\ell (z_jx^n)(x^{-n}b_j)$. Write $\dsp{\sur{z}_j=
\sum_{t=0}^{n-1} c_t \sur{\left(\frac{y}{x}\right)}^t}$ where $c_t \in k$. So
letting $a_t$ be an element of $A$ such that $\sur{a}_t = c_t$ and $v_j \in
A$ be the element $v_j= \sum_{t=0}^{n-1}a_ty^tx^{n-t}$, we have
\begin{equation} \label{ineg} \nu_\alpha(v_j-z_jx^n)>
  n\nu_\alpha(x).
\end{equation}

\begin{lem}
For $j \geq 2$, $x^{-n}b_j \in A$.
\end{lem}

\noi Proof of lemma : By Corollary \ref{cor-ppal1},
$\nu_{\mathfrak{m}}(b_j) > \nu_{\mathfrak{m}}(b_1)$ and by Corollary
\ref{cor-ppal2},
$$
\nu_{\mathfrak{m}}\left(\prod_{r=3}^iQ_r^{\gamma_{jr}}\right)
\leq \nu_{\mathfrak{m}}\left(\prod_{r=3}^iQ_r^{\gamma_{1r}}\right)-n.
$$
Now $\gamma_{j1}= \nu_{\mathfrak{m}}(b_j)-\nu_{\mathfrak{m}}
\left(\prod\limits_{r=3}^iQ_r^{\gamma_{jr}}\right)
> \nu_{\mathfrak{m}}(b_1) - \nu_{\mathfrak{m}}
\left(\prod\limits_{r=3}^iQ_r^{\gamma_{1r}}\right) +n \geq n$. $\Box$
 \smallskip

Put $z = b_1+ \sum_{j=2}^\ell v_j(x^{-n}b_j)$  We have $z \in IA'\cap A = I$
(because $I$ is a contracted ideal).

In both cases, $\nu_\alpha(x)=\nu_\alpha(y)$ and
$\nu_\alpha(x)<\nu_\alpha(y)$, since $z^\dag$ changes sign
between $\alpha'$ and $\beta'$ and in view of (\ref{ineg}), $z$
changes sign between $\alpha$ and $\beta$. This ends the proof of the
proposition. $\Box$

\begin{rek} If $B= R[x,y]$. Let
$$
B \to R[x_1,y_1] \to \cdots \to R[x_\ell,y_\ell]
$$
be a sequence of blowings up induced by (\ref{suite}), where
we take
$I=<\alpha,\beta>$.
Let $C_\ell$ be the preimage of
  $C$ (see (\ref{ensC})) in $\sper\ R[x_\ell,y_\ell]$. By proposition
  (\ref{eclat}), there
  exist monomials
$\omega_1,\dots,\omega_s ,\epsilon_1,\dots,\epsilon_s, \theta_1,\dots,
\theta_t$, $\lambda_1,\dots$, $\lambda_t$ in $x_\ell,y_\ell$ such that
$$
C_\ell=\left\lbrace \delta \in \sper\ R[x_\ell,y_\ell] \ \left| \ \begin{array}{c}
\nu_\delta(\omega_k) < \nu_\delta(\epsilon_k), \ k \in \{1,\dots,s\}
\\ \nu_\delta(\theta_j)= \nu_\delta(\lambda_j), \ j \in \{1,\dots,t\}
\\ sgn_\delta(x_\ell) = sgn_\alpha(x_\ell), sgn_\delta(y_\ell) =
sgn_\alpha(y_\ell)  \end{array} \right. \right\rbrace.
$$

\noi By connectedness theorem (\cite{LMSS}), $C_\ell$ is connected, hence
so is $C$. This completes the proof of the Connectedness Conjecture for
$R[x,y]$ and so provides a new proof of the classical Pierce-Birkhoff
Conjecture in dimension 2.
\end{rek}

Let $A$ be a regular 2-dimensional local ring with regular
parameters $(x,y)$.
Consider the set $C'$ defined by the
inequalities (\ref{eq:dominant1})
\begin{equation} \label{eq:dominant2}
\left|\sum_{j=1}^{m_i} b_{ji}\mathbf{Q}^{\theta_{ji}}\right| >_\delta
n_i|\mathbf{Q}^{\epsilon_{j'i}}| \
\forall i \in \{1,\ldots,r\}, \ \forall j' \in \{1,\ldots,n_i\}
\end{equation}
and the two sign conditions appearing in (\ref{ensC}).
\medskip

Consider the sequence (\ref{suite}) of local blowings up with
$I=<\alpha,\beta>$. Let $C_\ell'$ be the preimage of $C'$ in $\sper\
A_\ell$. Rather than prove connectedness of $C'_\ell$, we will prove that
$\alpha^{(\ell)}$ and $\beta^{(\ell)}$ lie in the same connected component of
$C' _\ell$; this will imply that $\alpha$ and $\beta$ lie in the same
connected component of $C'$. Let $\epsilon$ denote the common
specialization of $\alpha^{(\ell)}$ and $\beta^{(\ell)}$. By
definition of (\ref{suite}), we have $\mathfrak
p_\epsilon=\mathfrak{m}_\ell$. Let $U$ be the subset of $C'_\ell$
consisting of all the generizations of $\epsilon$ lying in $C'_\ell$. It is
sufficient to prove that $\alpha^{(\ell)}$ and $\beta^{(\ell)}$ lie in the
same connected component of $U$.
\medskip

\noi There are two cases to consider.

\textbf{Case 1.} Only one component of the exceptional divisor (that
is the inverse image $\rho_{\ell-1}^{-1}(\mathfrak{m})$) passes
through $\eta_\ell$.
\smallskip

\textbf{Case 2.} Two components of the exceptional divisor pass
through $\eta_\ell$. 
\smallskip

\noi Let $(x_\ell,y_\ell)$ be a regular system of parameters of
$(A_\ell)_{\mathfrak{m}_\ell}$ such that the local equation of the
exceptional divisor at $\eta_\ell$ is $x_\ell=0$ in case 1 and $x_\ell
y_\ell=0$ in case 2.
\medskip

By Zariski's theory of complete ideals, for any $f \in A \setminus
<\alpha,\beta>$, the strict transform of $f$ in $A_\ell$ is a unit. In
other words, $f$ has the form $f=x_\ell^n v$ in case 1
(resp. $f=x_\ell^n y_\ell^m v$ in case 2) where $v$ denotes a unit in
$(A_\ell)_{\mathfrak{m}_\ell}$.

The inequalities (\ref{eq:dominant2}), appearing in
the definition of $C'$, hold on all of $U$. The set $U$ is defined inside
the set of generizations of $\epsilon$ in $\sper\ A_\ell$ either by specifying
sgn($x_\ell$) or by specifying both sgn($x_\ell$) and sgn($y_\ell$).

\begin{lem}
Let $E$ be an irreducible component of the exceptional divisor passing
through $\eta_\ell$, defined by $x_\ell=0$. There exists $f \in A
\setminus <\alpha,\beta>$ such that $f=x_\ell^nv$, $v$ is a unit of
$(A_\ell)_{\mathfrak{m}_\ell}$ and $n$ is odd.
\end{lem}

\noi Proof : Let $j \in \{1, \ldots, \ell-1\}$ be such that $E$ is the
strict transform in $X_\ell$ of $\pi_{j-1}^{-1}(\eta_{j-1})$. Let
$\nu_j$ be the divisorial valuation corresponding to $E$; this
valuation is defined as follows : for each $f \in A_\ell$, write
$f=x_\ell^ng$ such that $x_\ell \nmid g$ in
$(A_\ell)_{\mathfrak{m}_\ell}$, then $\nu_j(f)=n$.

Let $\mathfrak{m}=\mathfrak{p}_0 \supset \cdots \supset
\mathfrak{p}_j$ be the complete list of simple $\nu_j$-ideals given by
Zariski's theory of complete ideals. Note that, since $j<\ell$,
$\mathfrak{p}_j \supset <\alpha,\beta>$.

It follows from Zariski's factorization theorem for complete ideals
that $\nu_j(A\setminus \{0\})$ is generated by
$\nu_j(\mathfrak{p}_0),\ldots,\nu_j(\mathfrak{p}_j)$. Since the value
group of $\nu_j$ is $\Z$, the semigroup $\nu_j(A\setminus \{0\})$
contains all the sufficiently large integers. Hence one of
$\nu_j(\mathfrak{p}_0),\ldots,\nu_j(\mathfrak{p}_j)$ is odd.
$\Box$
\medskip

The lemma shows that $x_\ell$ does not change sign between
$\alpha^{(\ell)}$ and $\beta^{(\ell)}$ in Case 1 (resp. neither
$x_\ell$ nor $y_\ell$ change sign between $\alpha^{(\ell)}$ and
$\beta^{(\ell)}$ in Case 2).
\medskip

Let
$$
\tilde{U} = \left\{\delta \in \sper\ A_\ell\ \left|\
sgn(x_\ell(\delta))= sgn(x_\ell(\alpha)), \epsilon\in\overline{\{\delta\}}\right.\right\}
$$
in Case 1 and
$$
\tilde{U} = \left\{\delta \in \sper\ A_\ell\ \left|\
sgn(x_\ell(\delta))= sgn(x_\ell(\alpha)), sgn(y_\ell(\delta))=
sgn(y_\ell(\alpha)), \epsilon\in\overline{\{\delta\}}\right.\right\}
$$
in Case 2. The above reasoning shows that $\alpha^{(\ell)}, \beta^{(\ell)} \in
\tilde{U} \subset U$.
\medskip

To prove the Definable Connectedness Conjecture (resp. the
Connectedness Conjecture for excellent $A$), it remains to prove the
definable connectedness of $\tilde{U}$ (resp. connectedness of
$\tilde{U}$ whenever $A$ is excellent).
\medskip

We are now ready to prove the above two versions of the Connectedness
Conjecture.

\subsection{Proof of the Connectedness Conjecture in the case of an
  excellent regular 2-dimensional ring.}\label{connectedexcellent}

\begin{thm}
Let $A$ be an excellent regular local 2-dimensional ring. Let $C
\subset \sper\ A$ be the subset satisfying the conditions of
(\ref{ensC}). Then $\alpha$ and $\beta$ belong to the same connected
component of $C$.
\end{thm}

\noi Proof: Let $\epsilon, \ell$ and $\tilde{U}$ as above. By the above
considerations, it is sufficient to prove that $\tilde{U}$ is
connected. Thus it remains to prove the following lemma.

\begin{lem}
Let $A$ be an excellent regular n-dimensional local ring, $x_1,\dots,x_n$
regular parameters of $A$. Fix a subset $T \subset \{1,\ldots,n\}$ and
let $D = \{\delta \in \sper\ A \ | \
x_i(\delta)> 0, \ i \in T \mbox{ and }\epsilon \in \sur{\{\delta\}}
\}$. Then $D$ is connected.
\end{lem}

\noi Proof: The point $\epsilon$ determines an order on $k$. Let $R$
denote the real closure of $k$ relative to this order. Consider the
natural homomorphisms
\begin{equation}
A \to \hat{A}=k[[X_1,\ldots,X_n]] \stackrel{\sigma}{\to} R[[X_1,\ldots,X_n]]
\end{equation}
where $\sigma$ is induced by $\epsilon$.

Let $\hat{\epsilon}$ denote the point of $\sper\ \hat{A}$ such that
$\mathfrak{p}_{\hat{\epsilon}} =(X_1,\ldots,X_n)$ and $\leq_{\hat{\epsilon}}$
is the total ordering of $k$ given by $\epsilon$.

Following (\cite{And}, proposition 8.6), $D$ is connected if and only
if $$\hat{D} = \{\delta \in \sper\ k[[X_1,\dots,X_n]] \ | \
X_i(\delta)>0, \ i \in T, \hat{\epsilon} \in \sur{\{\delta \}} \}$$ is
connected (this is where we are using the fact that $A$ is
excellent). Moreover, $\hat{D}$ is the image of
$$
\tilde{D} =\{\delta \in \sper\ R[[X_1,\dots,X_n]] \ | \
X_i(\delta)> 0, \ i \in T \}
$$ under the natural map induced by $\sigma$
$$
\sper\ R[[X_1,\ldots,X_n]] \to \sper\ k[[X_1,\ldots,X_n]].
$$
Thus it suffices to prove that $\tilde{D}$ is
connected.

By (\cite{And}, proposition 8.6), $\tilde{D}$ is connected if and only
if the set
$$
D^\dag = \{\delta \in \sper\ R[X_1,\dots,X_n]_{(X_1,\dots,X_n)} \ | \
X_i(\delta)> 0, \ i \in T,\  \delta \mbox{ is centered at
}(X_1,\dots,X_n) \}
$$
is connected.

We have the following natural homomorphisms
$$
\xymatrix{R[X_1,\dots,X_n] \ar[d]_-\psi \ar[r]^-\phi &
R[X_1,\dots,X_n]_{X_1\cdots X_n}
  \\ R[X_1,\dots,X_n]_{(X_1,\dots,X_n)}}
$$ and the corresponding maps of real spectra
$$
\xymatrix{\sper\ R[X_1,\dots,X_n]_{X_1,\dots,X_n} \ar[r]^-{\phi^*} &
\sper\ R[X_1,\dots,X_n]
  \\ & \ar[u]_-{\psi^*} \sper\ R[X_1,\dots,X_n]_{(X_1,\dots,X_n)}}
.$$

Define
$$
D_0 = \{\delta \in \sper\ R[X_1,\dots,X_n] \ | \
X_i(\delta)> 0, \ i \in T,\  \delta \mbox{ is centered at
}(X_1,\dots,X_n) \}
$$
and
$$
D_{loc}= \{\delta \in \sper\ R[X_1,\dots,X_n]_{X_1 \cdots X_n} \ | \
X_i(\delta)> 0, \ i \in T,\  \phi^*(\delta) \mbox{ is centered at
}(X_1,\dots,X_n) \}.
$$
Now the maps $\phi^*$ and $\psi^*$ induce homeomorphisms
\begin{eqnarray}
\phi^*|_{D_{loc}}:D_{loc}&\cong&D_0\quad\text{ and }\\
\psi^*|_{D^\dag}:D^\dag&\cong&D_0.
\end{eqnarray}

Thus it suffices to prove that $D_{loc}$ is connected. But
$$
D_{loc} =
\bigcap_{N \in \N} D_N
$$ where

$$
D_N = \left\{\delta \in \sper\ R[X_1,\dots,X_n]_{X_1 \cdots X_n} \ \left| \
\frac{1}{N} \geq X_i(\delta) \geq  0, \ i \in T \right.\right\}.
$$

By Proposition 7.5.1. of \cite{BCR}, each $D_N$ is a non-empty closed
connected subset of \\ $\sper\
R[X_1,\dots,X_n]_{X_1 \cdots X_n}$, hence $D_{loc}$ is connected by
(\cite{LMSS}, lemma 7.1). $\Box$
\medskip

The lemma proves that any ``quadrant'' is connected, $\tilde{U}$ is
a quadrant, hence it is connected. This completes the proof of the
Connectedness Conjecture for any excellent 2-dimensional ring $A$.

\begin{rek}
The above proof is a special case of the following general
principle. Let $A$ be an excellent regular local ring with regular
parameters $x=(x_1,\dots,x_n)$ whose residue field $k$ is equipped with a
total ordering. Let $R$ be the real closure of $k$. We
have natural morphisms
$$
\xymatrix{\  \sper\ A &  \sper\ R[[X_1,\dots,X_n]] \ar[d]_-\pi \ar[l]^-\phi
  \\ & \sper\ R[X_1,\dots,X_n]_{(X_1,\dots,X_n)}}
$$
Let $D \subset \sper\ A$ be a constructible set such that all the
elements of $A$ appearing in the definition of $D$ belong to $A \cap
R[X_1,\dots,X_n]_{(X_1,\dots,X_n)}$. Let $\hat{D} = \phi^{-1}(D)$, let
$U$ be the subset of all points of $\sper\
R[X_1,\dots,X_n]_{(X_1,\dots,X_n)}$ centered at the origin. Let
$D_{pol}$ be the subset of $U$ defined by the same formulae as $D$. By
(\cite{And}, proposition 8.6), to show that $D$ is connected, it is
enough to prove that $D_{pol}$ is connected.

In many cases, this principle applies also to nested intersection $\dsp{D
  = \bigcap_{N \in \N} D_N}$ of constructible sets defined by elements of
$A \cap R[X_1,\dots,X_n]_{(X_1,\dots,X_n)}$.

This allows to transpose all the results of (\cite{LMSS}) from the
case of polynomial rings to that of arbitrary excellent regular local
rings.
\end{rek}

\subsection{Proof of the Definable Connectedness Conjecture for regular
2-dimensional local rings.}\label{proof}
Next we prove the Definable Connectedness Conjecture, hence the
Pierce-Birkhoff Conjecture, without the excellence hypothesis on $A$.

\begin{thm} \label{principale} Let $(A,\mathfrak{m},k)$ be a regular
2-dimensional local ring, $(x,y)$ a regular system of parameters of $A$. The
sets
\begin{eqnarray}
U&=&\{\delta \in \sper\ A \ | \
x(\delta)>0, \epsilon\in\overline{\{\delta\}}\} \\
V&=&\{\delta \in \sper\ A \ | \
x(\delta)>0,y(\delta)>0,
\epsilon\in\overline{\{\delta\}}\}
\end{eqnarray}
are definably connected.
\end{thm}

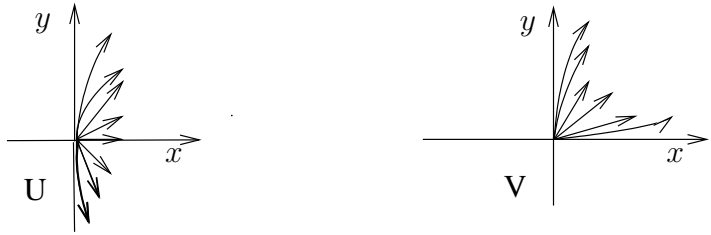
\begin{figure}[htbp]
\begin{center}
\input{UV2.pstex_t}
\end{center}
\caption{The sets $U$ and $V$}
\end{figure}

\noi Proof : We argue by contradiction. Let $\Omega$ be either $U$ or
$V$. Write $\Omega= F \coprod G$, $F=\bigcup
F_i$, $G=\bigcup G_i$ where $\{F_i\},\ \{G_i\}$ are finite collections
of basic open sets. Each $F_i$ and
$G_i$ is defined by finitely many inequalities of the form $g>0$, $g
\in A$. Let $g_3,\ldots,g_r \in A$ be the list of elements
of $A$, appearing in the definition of all of $F_i$ and $G_i$ and let
$g_1=x$, $g_2=y$. A proof of the
Theorem will be given after a
few auxiliary definitions and results.
\medskip

Let $\sper\, A \leftarrow X_1 \leftarrow \cdots \leftarrow X_t$ be a sequence
of point blowings up. Let $X_t= \bigcup\limits_{j=1}^s \sper\, A_{jt}$ be the
open covering of $X_t$, given by the definition of real geometric surface.

\begin{deft}  We say that a collection
  $\{h_1,\ldots,h_r\}$ of elements of $A$ are simultaneously locally
  monomial in $X_t$ if for all $j \in \{1,\ldots,s\}$ and any maximal
  ideal $\mathfrak{m}' \subset A_{jt}$, there exists a regular system
  of parameters $(x',y')$ of $A':=(A_{jt})_{\mathfrak{m}'}$ such that
  all of $h_1,\ldots, h_r$ are monomials in $(x',y')$ multiplied by
  units of $(A_{jt})_{\mathfrak{m}'}$.
\end{deft}

Let $g_1,\ldots,g_r \in A$ be as above. By standard results on
resolution of singularities, there exists a sequence $\sper\, A
\leftarrow X_1 \leftarrow \cdots \leftarrow X_t$ of point blowings up
such that $g_1,\ldots, g_r$ are simultaneously locally monomial in
$X_t$. Denote by $\rho_t:X_t \to \sper\,A$ the composition of all the
morphisms in that sequence (with the notations following
(\ref{def:sequence})).

Let
$\Omega^{(t)}=\rho_t^{-1}(\Omega)$, $F^{(t)}=\rho_t^{-1}(F)$, $G^{(t)}=\rho_t^{-1}(G)$, $U^{(t)} = \rho_t^{-1}(U)$.

Take a point $\delta \in \rho_t^{-1}(\epsilon)$, let
$A',\mathfrak{m}',x',y',A_{jt}$ be as in the definition of
simultaneously locally monomial.

\begin{deft}
We say that $\delta$ is a special point of $\rho_t^{-1}(\epsilon)$ if
$ht(\mathfrak{p}_\delta)=2$ and
$$
\{x'y'=0\} = \rho_t^{-1}(\epsilon)\cup \{g_1\cdots g_r=0\}
$$
locally near $\delta$.
\end{deft}

 Given a special point $\delta \in \rho_t^{-1}(\epsilon)$ and
$(u',v')$ a regular system of parameters at $\delta$,
let $$C(\delta,u',v') = \{\gamma \in X_t \ | \ u'(\gamma)>0,\
v'(\gamma) >0,\ \delta \in \sur{\{\gamma\}} \}.$$

\begin{lem}   
Take a point $\xi \in \rho_t^{-1}(\epsilon)$, not lying on the strict
transform of $\{x=0\}$. Take $j \in \{1,\ldots,s_i\}$ such that
$\rho_t^{-1}(\epsilon)$ is contained in the privileged set of $\sper\
A_{jt}$ near $\xi$. Let $x_{jt},y_{jt} \in A_{jt}$ be the elements
given in Definition \ref{def-mrgs}. Assume that the privileged set is
given by $\{x_{jt}=0\}$ and is homeomorphic to $\sper\
k'[z]_{\theta_z}$, where $\theta_z$ is a non-zero polynomial, with
$k'$ finite over $k$ and that ht($\mathfrak{p}_\xi$)=2. Let $(x',y')$
be as in the definition of simultaneously locally monomial where we
take $\mathfrak{m}'=\mathfrak{p}_\xi$ (we may assume $x'=x_{jt}$). We
view $k'$ as an ordered field via the inclusion $k' \subset
A_{jt}(\xi)$. Let $$E=\{
\delta \in \sper\ A_{jt} \ | \ x_{jt}(\delta)=0 \mbox{ and } k' \subset
A_{jt}(\delta) \text{ is an inclusion of ordered fields } \}.$$

Take special points $\delta_1,\delta_2 \in E$ such that the intervals
$(\delta_1,\xi)$ and
$(\xi,\delta_2)$ are connected and contain no
special points.

For $i \in \{1,2\}$, let $(x',v_i')$ be a regular system of parameters
at $\delta_i$ such that $\{v'_i>0\} \cap (\delta_1,\delta_2) \neq
\emptyset$.

Then the set
$$ D(\delta_1,\delta_2)=
C(\delta_1,x',v'_1) \cup C(\delta_2,x',v'_2)
\cup
\{\delta \in U^{(t)} \ | \ x'(\delta) >0,\ \sur{\{\delta\}}
\cap (\delta_1,\delta_2) \neq \emptyset \}
$$
is contained either in $F^{(t)}$ or in $G^{(t)}$.
\end{lem}

\noi Proof : First, assume $\xi$ is not special. Then there are no special
points in $(\delta_1,\delta_2)$. Let $F_\dag=\sur{F^{(t)}} \cap
[\delta_1,\delta_2]$ and $G_\dag = \sur{G^{(t)}} \cap
[\delta_1,\delta_2]$. Then
 $F_\dag$, $G_\dag$ are relatively closed in $[\delta_1,\delta_2]$ and
 $[\delta_1,\delta_2]$ is connected (Lemma \ref{inter-con})
 , so $F_\dag
 \cap G_\dag \neq \emptyset$.

Take a point $\eta \in F_\dag \cap G_\dag$. Replacing $\eta$ by
its specialization, we may assume that $ht(\mathfrak p_\eta)=2$. For each $i
\in \{1,\ldots,r\}$, locally near $\eta$, write $g_i =
x'^a g'_i$ if $\eta \notin \{\delta_1,\delta_2\}$ and $g_i= x'^a
y'^bg'_i$ if $\eta = \delta_\ell,\  \ell \in \{1,2\}$ with $y'=v'_i$,
where, in both cases, $g'_i$ is invertible locally near $\eta$.

Take an open set $W$, containing $\eta$, such that for all $\delta \in
W$ and all $i \in \{1,\ldots,r\}$, we have
\begin{equation}\label{eq:signes}sgn(g'_i(\delta))=sgn(g'_i(\eta)).
\end{equation}

Since $\eta \in \sur{F^{(t)}} \cap \sur{G^{(t)}}$, there exist $\delta \in
F^{(t)} \cap W$, $\gamma \in G^{(t)}\cap W$ and an $i \in \{1,\ldots,r\}$ such
that $g_i$ changes sign between $\delta$ and $\gamma$.

Since $x'$ (resp. $x',y'$) does not change sign between $\gamma$ and
$\delta$ this contradicts (\ref{eq:signes}).

Therefore $\xi$ must be special.
Let $\delta \in
D(\delta_1,\delta_2)$
be the unique point such that
$$
x'(\delta)>0, y'(\delta)=0.
$$
We have $\{\delta\} = D(\delta_1,\delta_2)\cap \{y'=0\}$. Then $$
D(\delta_1,\delta_2) = \{\delta\} \coprod D(\delta_1,\xi) \coprod
D(\xi,\delta_2).$$

Let $\delta_- \in D(\delta_1,\xi)$ be the unique point such that
$x'(\delta_-)>0, y'(\delta_-)<0$ and $|y'(\delta_-)|^N<
|x'(\delta_-)|, \forall N \in \N$. Then $\delta \in
\sur{ \{\delta_- \} }$, in particular, \begin{equation} \delta \in
\sur{D(\delta_1,\xi)}. \label{blubs1}
\end{equation}

Similarly \begin{equation} \delta \in
\sur{D(\xi,\delta_2)}. \label{blubs2}
\end{equation}

By the previous case, each of $D(\delta_1,\xi)$, $D(\xi,\delta_2)$ is
contained either in $F^{(t)}$ or $G^{(t)}$.

Without loss of generality, assume that $D(\delta_1,\xi) \subset
F^{(t)}$. By (\ref{blubs1}) and the relative
closedness of $F^{(t)}$, we have
$\delta \in F^{(t)}$. By (\ref{blubs2}) and the relative closedness
of $G^{(t)}$, we have $D(\xi,\delta_2) \subset F^{(t)}$, so
$D(\delta_1,\delta_2) \subset F^{(t)}$ as desired. $\Box$
\medskip

\begin{cor}\label{loc-conn}
Let $[\delta_1,\delta_2]$ be a maximal interval. Then
$D(\delta_1,\delta_2)$ is entirely contained either in $F^{(t)}$ or
in $G^{(t)}$.
\end{cor}

\noi Proof : This follows from the preceding lemma by induction on the
number of special points inside $[\delta_1,\delta_2]$.
\medskip

  In order to address the global connectedness, we need a notion of
\textbf{signed dual graph} associated to a sequence of point
blowings up of a point $\epsilon \in \sper\, A$ and a subset $W$ of
$\sper\, A$.
\vspace{-.3cm}

\noi For each maximal interval $I$ (see Definition \ref{def-max}), take
$\sper\ A_{jt} \subset X_t$ such that $I \setminus\sper\ A_{jt}$ is
either empty or
consists of one distinguished point (such an $A_{jt}$ exists by
Proposition \ref{prop-max}). When necessary, we will denote this $j$
by $j(I)$. Let $x_{jt},y_{jt}\in A_{jt}$ be the elements
given in the Definition \ref{def-mrgs}. By Proposition
\ref{prop-max}, we have $I \cap \sper\ A_{jt} \subset \{x_{jt}=0\}$.
\medskip

Let $W^{(t)}=\rho_t^{-1}(W)$. Let $I$ a maximal interval, denote by
$I^\circ$ its interior, and $s \in \{+,-\}$, let
\begin{equation}
  W(I,s) = \{\delta \in \sper\, A_{jt}\ |\
sgn(x_{jt}(\delta)) = s, \, \sur{\{\delta\}} \cap
I^\circ \neq \emptyset \}.
\end{equation}

\begin{deft} \label{def:admissible}
Consider a pair $(I,s)$ as above. We say that $(I,s)$ is admissible if
\begin{equation}
 W^{(t)} \cap \sper\ A_{jt} \supset W(I,s) \neq \emptyset.
\end{equation}
\end{deft}

Consider two admissible pairs $(I,s)$, $(\tilde{I},\tilde{s})$. We say
that these two pairs are \textbf{equivalent} if the following conditions hold :

(a) $I \cap  \sper\ A_{jt} \cap \sper\ A_{\tilde{\jmath}t} = \tilde{I} \cap
\sper\ A_{jt} \cap \sper\ A_{\tilde{\jmath}t}$,

(b) the sets $$\{\delta \in \sper\ A_{jt} \cap \sper\ A_{\tilde{\jmath}t}\
|\ sgn(x_{jt}(\delta)) = s \}$$ and $$\{\delta \in \sper\ A_{jt} \cap
\sper\ A_{\tilde{\jmath}t}\ |\ sgn(x_{\tilde{\jmath}t}(\delta)) = \tilde{s} \}$$
coincide in a neighbourhood of $I \cap  \sper\ A_{jt} \cap \sper\
A_{\tilde{\jmath}t}$.

Given two equivalent admissible pairs $(I,s)$ and $(\tilde I,\tilde s)$, the set of endpoints of $I$ coincides with the set of endpoints of $\tilde I$ (viewed as points of the marked real geometric surface $X_t$). In this way, given an equivalence class of admissible pairs $\{(I,s)\}$, it makes sense to talk about endpoints of $\{(I,s)\}$.
\begin{deft} \label{def:graphe}
1. A vertex of the signed dual graph $\Gamma_t$ associated to $X_t$ and $W$
is an equivalence class of
an admissible
pair $(I,s)$, which we will still denote, by
abuse of notation, by $(I,s)$.

2. By definition, two distinct vertices $(I,s)$ and $(\tilde{I}, \tilde{s})$ of
$\Gamma_t$ are connected by an edge of $\Gamma_t$ if the following
conditions hold :

(a) $I$ and $\tilde{I}$ share a common endpoint
$\xi$ and suppose that
$\tilde{I} \notsubset \{x_{ji}=0\}$;

(b) we have
$$W^{(t)} \cap
\{\delta \in X_t \ |\
sgn(x_{jt}(\delta))=s, \ sgn(x_{\tilde{\jmath}t}(\delta))=\tilde{s} ,\
\sur{\{\delta\}} \ni
 \xi \} \neq \emptyset.$$
\end{deft}

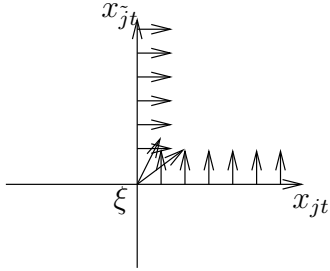
\begin{figure}[h!]
\begin{center}
\input{edge.pstex_t}
\caption{This figure represents an edge of $\Gamma_t$
    connecting two vertices $(I,s)$ and $(\tilde{I},\tilde{s})$.
Here $I=[0,\infty]$, $\tilde{I}=[0,\infty]$,
    $s=\tilde{s}=+$.}
\end{center}
\end{figure}

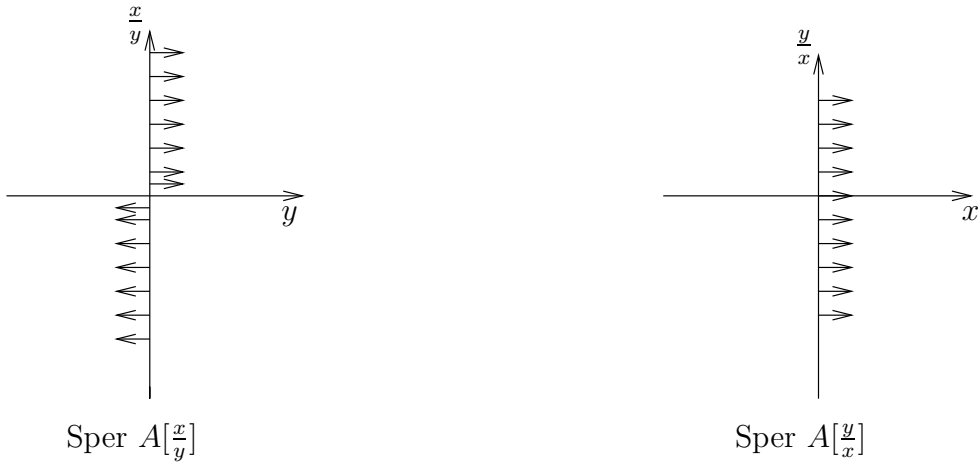
\begin{figure}[h!]
\begin{center}
\input{U1.pstex_t}
\caption{This figure shows the set $U^{(1)}$ in the affine charts} \label{u1}
\end{center}
\end{figure}

\noi Example:
If $W = U$ or $W=V$ then
$\Gamma_1$ consists of one vertex and no edges (see
Fig. (\ref{u1}) for a picture of $U^{(1)}$; the case of $V^{1}$ is similar but easier).
\medskip

\begin{prop}\label{bamboo} If $W=U$ or $W=V$, the graph $\Gamma_i$ is a bamboo, that is, a connected, simply connected graph every one of whose vertices belongs to at most two edges.
\end{prop}

\noi Proof: By induction on $i$. For $i=1$, the graph consisting of
one vertex is connected and satisfies the conclusion of the Proposition. The induction step follows from the next
Lemma, which describes the transformation law from
$\Gamma_i$ to $\Gamma_{i+1}$ in the case when $W=U$ or $W=V$.
\medskip

Consider the point blowing up $\pi_i:X_{i+1} \to X_i$.  Let $\xi$ be
the center of the blowing up; recall that, by
definition of blowing up in the category of real marked geometric surfaces,
$\xi$ belongs to the distinguished set of $X_i$. Let $\sper\, A_{ji}$ be an
affine chart of $X_i$ containing $\xi$. Let $\mathfrak{p}_\xi$ be the
support of $\xi$ in $A_{ji}$. Let $k_{ji}$ be the field of Definition
\ref{def-mrgs} (2). Let $E_1,\ldots,E_p$ be the components of the set
$\{x_{ji}=0\} \cap \rho_i^{-1}(\epsilon)$. Picking a component $E_q,\ q \in
\{1,\ldots,p\}$ amounts to fixing a total order on $k_{ji}$, which induces the
order on $k$ given by $\epsilon$. For $q \in \{1,\ldots,p\}$, let
$\{\xi_1^{(q)},\ldots,\xi_\ell^{(q)}\}$ be the set of points of $E_q$
supported at $\mathfrak{p}_\xi$.  For each $q \in \{1,\ldots,p\}$, the total
order on $k_{ji}$ corresponding to $E_q$ induces a total order on the set
$\{\xi_1^{(q)},\ldots,\xi_\ell^{(q)}\}$. Renumbering
$\{\xi_1^{(q)},\ldots,\xi_\ell^{(q)}\}$, we may assume
$\xi_1^{(q)} < \xi_2^{(q)} < \cdots < \xi_\ell^{(q)}$.

It follows from the definition of distinguished that one of the points
$\xi_t^{(q)}$ is $j$-distinguished if and only
if
all of them are.

Fix a pair $(q,t)$,  $q \in \{1,\ldots,p\}$,  $t \in \{1,\ldots,\ell\}$. Two
cases are possible :

\noi $\bullet$ Case 1 : There exist $a= (I,s),\ b=(\tilde{I},\tilde{s})$ two
vertices of
$\Gamma_i$ connected by an edge $(a,b)$ such that $\xi_t^{(q)}$ is the point
common to $I$ and $\tilde{I}$ (note that the pair $a,b$ is not, in general
uniquely determined by $\xi_t^{(q)}$). In particular, the points $\xi_t^{(q)}$
are $j$-distinguished. In this case, we have $p=1$, so we may denote our points
by $\xi_1,\ldots,\xi_\ell$. Let $x_{ji},\tilde{x}_{ji}$ be a privileged regular
system of parameters at the points $\xi_t$.

\noi $\bullet$ Case 2 : We are not in Case 1.

- Case 2.1. : None of the points $\xi_t^{(q)}$ is $j$-distinguished. Let
$(x_{ji},y')$ be a regular system of parameters of the local ring
$A_{\mathfrak{p}_\xi}$. The set $\pi_i^{-1}(\xi_t^{(q)})$ is covered by two
affine charts : $\sper\, A_{ji}[\frac{x_{ji}}{y'}]$ and $\sper\,
A_{ji}[\frac{y'}{x_{ji}}]$. Let $x'_{ji} = \frac{x_{ji}}{y'}$.

- Case 2.2 : The point $\xi$ is $j$-distinguished and lies on the strict
transform of $\{x=0\}$ or $\{y=0\}$. In this case, $p=\ell=1$.
\medskip

Next, we study the neighbourhood of $\pi_i^{-1}(\xi_t^{(q)})$ for each $q \in
\{1,\ldots,p\},\ t \in \{1,\ldots,\ell\}$  and analyze the changes from
$\Gamma_i$ to $\Gamma_{i+1}$ induced by the blowing-up $\pi_i$ locally on the
part of $\Gamma_i$ which represents a neighbourhood of $\xi_t^{(q)}$. Since
$\pi_i$ induces an isomorphism outside the points $\xi_t^{(a)}$, the rest of
the graph $\Gamma_i$ remains unchanged under the blowing-up $\pi_i$.

In the statement of the following lemma, we refer to the cases 1 and 2 defined
above.

\begin{lem}\label{graphtransformation}

$\bullet$ Case 1 : Fix $t \in \{1,\ldots,\ell \}$. For each pair of vertices
$a$,
$b$ as above, remove the edge $(a,b)$ and add a new vertex $c$ and two new edges
$(a,c)$ and $(b,c)$. The graph $\Gamma_{i+1}$ is obtained from $\Gamma_i$ by
successively performing the above operation for each of $\xi_1,\ldots,\xi_\ell$.
\smallskip

$\bullet$ Case 2 : Consider a vertex $a= (I,s)$ such that $\xi_t^{(q)} \in I$
for some $t \in \{1,\ldots,\ell \}$ and $q \in \{1,\ldots, p \}$. Write
$I=[\delta_1,\delta_2]$ (again, the vertex $a$ is not, in general, uniquely
determined by $\xi_t^{(q)}$).

 - Case 2.1: Take $\lambda \in \{0,\ldots,\ell-1\}$ and $\omega
\in \{1,\ldots,\ell -\lambda\}$ such that
\begin{equation}
 [\delta_1,\delta_2] \cap
\{\xi_t^{(q)}\ |\ t \in \{1,\ldots,\ell\},\ q \in \{1,\ldots,p\}\} =
\{\xi_{\lambda+1}^{(q)},\xi_{\lambda+2}^{(q)},\ldots,\xi_{\lambda+\omega}^{(q)}
\}
\end{equation}
for some $q \in \{1,\ldots,p \}$. Replace $a$ by a bamboo with $2\omega+1$
vertices. More precisely, we distinguish three cases :

\hspace{.6cm} (a) If $a$ belongs to two edges $(a,b),\ (a,c)$ of
$\Gamma_i$, remove $a$ and the two edges $(a,b)$, $(a,c)$. Introduce
the bamboo
$$\stackrel{b\, \bullet \text{\!---}\underbrace{\bullet\text{\!---}\bullet\text{
\!---}\bullet\text{ \!--- }\bullet \cdots \cdots \bullet\text{\!---}\bullet}\,
\text{\!---}\bullet\,  c}{2\omega+1} $$

\hspace{.6cm} (b) If $a$ belongs to only one edge $(a,b)$, remove $a$ and the
edge $(a,b)$ and introduce the bamboo
$$\stackrel{b\, \bullet \text{\!---}\underbrace{\bullet\text{\!---}\bullet\text{
\!---} \bullet\text{ \!--- } \bullet \cdots \cdots
\bullet\text{\!---}\bullet}}{2\omega+1} $$

\hspace{.6cm} (c) If $a$ belongs to no edges (in other words, if $i=1$) then
$$\Gamma_2=\stackrel{\underbrace{\bullet\text{\!---}\bullet\text{
\!---} \bullet\text{ \!--- } \bullet \cdots \cdots
\bullet\text{\!---}\bullet}}{2\omega+1}$$
is a chain of $2\omega+1$ vertices and $2\omega$ edges.

The graph $\Gamma_{i+1}$ is obtained from $\Gamma_i$ by performing
successively the above operation for each vertex $a$ as above.

 - Case 2.2:

\hspace{.6cm}(a) $i=1$ and $W=U$, then
$\Gamma_2=\bullet$\!---$\bullet$\!---$\bullet$ is a chain of three vertices and
two edges

\hspace{.6cm}(b) $i>1$ or $W=V$, then each vertex $a=(I,s)$ such that
$\xi \in I$ is an endpoint
of $\Gamma_i$. For each such vertex $a$,
we add a new vertex $b$ and a new edge $(a,b)$.

\end{lem}

\noi Proof: Case 1 : Let $\delta_1,\ \delta_2$ be points of $\sper\ A_{ji}$ such
that $I=[\delta_1,\xi_t]$,
$\tilde{I} = [\xi_t,\delta_2]$. Let
$\delta_1'=\pi_i^{-1}(\delta_1)$,
$\delta_2'=\pi_i^{-1}(\delta_2)$. Let
$x_{ji}, x_{\tilde{\jmath},i} \in A_{ji}$ be as in the Definition
\ref{def:admissible} applied to $(I,s)$ and $(\tilde{I},\tilde{s})$, respectively. The pair
$(x_{ji},x_{\tilde{\jmath},i})$ forms a regular system of parameters at $\xi_t$.
Let $x'_{ji} = \frac{x_{ji}}{x_{\tilde{\jmath},i}}$ and $x'_{\tilde{\jmath},i}
= \frac{x_{\tilde{\jmath},i}}{x_{ji}}$.

Let
$$
\xi_a \in \left\{x'_{ji}=0\right\}\cap \pi_i^{-1}(\xi_t) \subset
\sper\, A_{ji}[x'_{ji}]
$$
and
$$
\xi_b \in \left\{x'_{\tilde{\jmath}i}=0\right\} \cap \pi_i^{-1}(\xi_t)
\subset \sper\, A_{ji}[x'_{\tilde{\jmath},i}];
$$
note that these conditions characterize $\xi_a$ and $\xi_b$ uniquely.
Let
$J=[\xi_a,\xi_b]$, viewed as a maximal interval of $\sper\, A_{ji}[x'_{ji}]$. Let
$\sigma = s \cdot \tilde{s}$. 

Let $a_{i+1}, b_{i+1}, c_{i+1}$ be the vertices of $\Gamma_{i+1}$
defined by $a_{i+1}=([\xi_a,\delta_1'],\sigma)$,
$b_{i+1}=([\xi_b,\delta'_2],\sigma)$, $c_{i+1}=(J,\tilde{s})$.
We have to verify that those three pairs are admissible; first, we will show
the admissibility of $([\xi_a,\delta'_1],\sigma)$.

Since $(I,s)$ is admissible, we know that
$$\emptyset \neq  W(I,s) \subset W^{(i)} \cap \sper\,
A_{ji}$$ and we need to show that
\begin{equation} \label{eq:admiss}
\emptyset \neq W([\xi_a,\delta'_1],\sigma) \subset W^{(i+1)}
\cap \sper\, A_{ji}[x'_{ji}].
\end{equation}
To prove (\ref{eq:admiss}), note that $\pi_i$ induces an isomorphism outside the set
$\{\xi_{t'}\ |\ 1 \leq t'\leq \ell\}$; in particular, it induces an isomorphism
of a neighbourhood of the open interval $(\xi_a,\delta'_1)$ onto a neighbourhood of $(\xi_t,
\delta_1)$. Moreover, the fact that $a$ and $b$ are connected by an edge of
$\Gamma_i$ implies that $sgn(x_{\tilde{\jmath},i}(\delta)) = \tilde s$ for
$\delta \in W(I,s)$. Hence $\delta' \in \pi_i^{-1 }(W(I,s))$ if and only if
$\sur{\{\delta'\}} \cap (\xi_a,\delta'_1) \neq \emptyset$ and
$sgn(x'_{ji}(\delta'))= s \cdot \tilde s$. In other words, $W([\xi_a,
\delta'_1],\sigma) = \pi_i^{-1}(W(I,s))$. This proves
(\ref{eq:admiss}), so $([\xi_a,\delta'_1],\sigma)$ is
admissible. By symmetry, the pair $([\xi_b,\delta'_2],\sigma)$ is also
admissible.

To prove the admissibility of $(J,\tilde s)$, we note that
$x_{\tilde \jmath,i} =0$ is the local equation of the exceptional divisor in
$\sper\, A_{ji}[x'_{ji}]$ and hence $$\pi_i(W(J,\tilde s)) = \{\delta \in
\sper\, A_{ji}
\ |\ \sur{\{\delta\}} \ni \xi_t, \ sgn(x_{\tilde
\jmath, i}(\delta)) = \tilde s,\  \ sgn(x_{ji}(\delta)) = s\}.$$ Now the fact
that $a$ and $b$ are connected by an edge of $\Gamma_i$ (see Definition
\ref{def:graphe} (b)) implies that $$ \emptyset \neq W(J,\tilde s) \subset
W^{(i+1)} \cap \sper\, A_{ji}[x'_{ji}],$$ so $(J, \tilde s)$ is admissible.

To check that $a_{i+1}$ and $c_{i+1}$ are connected by an edge of
$\Gamma_{i+1}$, consider the set
$$
\{\delta' \in \sper\, A_{ji}[x'_{ji}]\ |\
\sur{\{\delta'\}} \ni \xi_a,\ sgn(x_{\tilde \jmath,i}(\delta'))= \tilde s,\
sgn(x'_{ji}(\delta'))=\sigma \}.
$$

We have \footnote{By $\delta$ tangent to
$\{x_{ji}=0\}$, we mean $\delta$ such that $\forall N \in \N,\
N|x_{ji}(\delta)| < |x_{\tilde \jmath,i}(\delta)|$ }

\begin{multline*}
 \pi_i(\{\delta' \in \sper\, A_{ji}[x'_{ji}]\ |\
\sur{\{\delta'\}} \ni \xi_a,\ sgn(x_{\tilde \jmath,i}(\delta'))= \tilde s,\
sgn(x'_{ji}(\delta'))=\sigma \}) \\
=  \{\delta \in \sper\, A_{ji}\ |\
\sur{\{\delta\}} \ni \xi_t,\ sgn(x_{\tilde \jmath,i}(\delta))= \tilde s,\
sgn(x_{ji}(\delta)) =s,\ \delta \mbox{ tangent} \mbox{ to } \{x_{ji}=0\} \} \\
\subset
\{\delta \in \sper\, A_{ji}\ |\
\sur{\{\delta\}} \ni \xi_t,\ sgn(x_{\tilde \jmath,i}(\delta))= \tilde s,\
sgn(x_{ji}(\delta)) =s \}\subset W^{(i)},
\end{multline*}
where the last inclusion comes from the fact that $a$ and $b$ are connected by an
edge in $\Gamma_i$.

Hence
\begin{multline}
 \emptyset \neq \{\delta' \in \sper\, A_{ji}[x'_{ji}]\ |\
\sur{\{\delta'\}} \ni \xi_a,\ sgn(x_{\tilde \jmath,i}(\delta'))= \tilde s,\
sgn(x'_{ji}(\delta'))=
\sigma
\} \\ \subset W^{(i+1)} \cap \sper\, A_{ji}[x'_{ji}],
\end{multline}
which proves that $a_{i+1}$ is connected to $c_{i+1}$. By symmetry, $b_{i+1}$
is also connected to $c_{i+1}$.
\smallskip

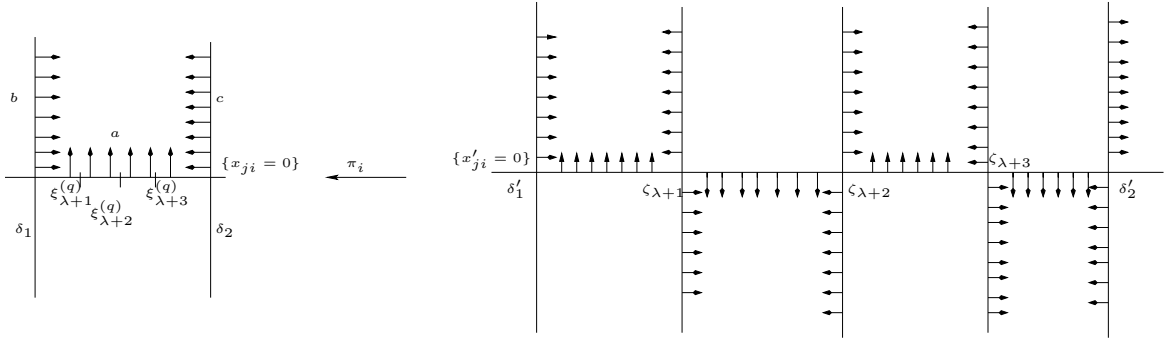
\begin{figure}[h!]
\begin{center}
\input{global4.pstex_t}
\caption{This figure shows, in the Case 2.1 (a), with $\omega=3$, the
transformation of the dual graph under the blowing up $\pi_i$.}
\label{Usuite}
\end{center}
\end{figure}
Case 2.1 (a) Recall that $(x_{ji},y')$ is the chosen regular system of
parameters at
$\mathfrak p_\xi$.

Let $\delta'_\tau=\pi_i^{-1}(\delta_ \tau)$, $\tau \in  \{1,2\}$.
Let $\zeta_t=\pi_i^{-1}(\xi_t^{(q)}) \cap \{x'_{ji}=0\}$, $t \in \{\
\lambda+1,\dots,\lambda+\omega \}$.
The new distinguished points in the open interval $(\delta'_1,\delta'_2)$ are
$\zeta_t$, $t \in \{\ \lambda+1,\dots,\lambda+\omega \}$. The components of
$\pi_i^{-1}(\mathfrak{p}_\xi)$ are $\pi_i^{-1}(\xi_{\lambda+1}^{(q)}), \ldots,
\pi_i^{-1}(\xi_{\lambda+\omega}^{(q)})$. For $t \in \{
\lambda+1,\dots,\lambda+\omega \}$, let us denote the interval
$[-\infty,+\infty] \subset \pi_i^{-1}(\xi_t^{(q)}) \cap \sper\
A_{ji}[\frac{y'}{
x_{ji}
}]$ by $[-\infty,\infty]_t$.

Now, there are $2\omega+1$ maximal intervals in
$\pi_i^{-1}((\delta_1,\delta_2))$. They are : $[\delta'_1,\zeta_{\lambda+1}]$,
$[\zeta_{\lambda+\omega},\delta'_2]$,
$[\zeta_t,\zeta_{t+1}]$, $t \in \{ \lambda+1,\ldots,\lambda+\omega-1\}$ and
$[-\infty,\infty]_t$, $t \in \{ \lambda+1,\dots,\lambda+\omega \}$.

Without loss of generality, we may assume that $y'(\delta_1) >0$. Each of this
maximal intervals gives rise to an admissible pair as follows.

The intervals $[\zeta_t,\zeta_{t+1}] \subset \{x'_{ji}=0\}$ give rise to
admissible pairs $([\zeta_t,\zeta_{t+1}], (-1)^t\cdot s)$.

We have admissible pairs $([\delta'_1,\zeta_{\lambda+1}], s)$ and
$([\zeta_{\lambda + \omega},\delta'_2],(-1)^\omega \cdot s)$. Finally, the
intervals $[-\infty,\infty]_t$ give rise to admissible pairs
$([-\infty,\infty]_t,s)$.

To see that the pair $([\zeta_t,\zeta_{t+1}], (-1)^t\cdot s)$ is admissible, we
use the fact that $\pi_i$ is an isomorphism from a neighbourhood of the open
interval $(\zeta_t,\zeta_{t+1})$ to a neighbourhood of the open interval
$(\xi_t^{(q)},\xi_{t+1}^{(q)})$. Since $y'(\delta_1)>0$ and since $y'$ changes
sign once at each point $\xi_t^{(q)}$ the sign of $y'$ on
$(\xi_t^{(q)},\xi_{t+1}^{(q)})$ is $(-1)^t$. Hence
\begin{multline}
\pi_i(\{\delta'\in X_{i+1} \ |\ \sur{\{\delta'\}} \cap (\zeta_t,\zeta_{t+1})
\neq \emptyset,\ sgn(x'_{ji}(\delta'))= (-1)^t \cdot s \}) \\
= \{\delta \in X_i\ |\
\sur{\{\delta\}} \cap (\xi_t^{(q)},\xi_{t+1}^{(q)}) \neq \emptyset,\
sgn(x_{ji}(\delta))= s \}.
\end{multline}

This proves the admissibility of $([\zeta_t,\zeta_{t+1}],(-1)^t \cdot s)$. The
proof that $([\delta'_1,\zeta_{\lambda+1}],s)$ and $
([\zeta_{\lambda+\omega},\delta'_2],(-1)^\omega \cdot s)$ are admissible is
similar and we omit it.

To prove the admissibility of $([-\infty,\infty]_t,s)$, note that
$$\pi_i^{-1}(\{ \delta \in X_i\ |\ \sur{\{\delta\}} \ni \xi_t^{(q)},\
sgn(x_{ji}(\delta)) =s\}) \supset W([-\infty,\infty]_t,s),$$ where the notation
$
W([-\infty,\infty]_t,s)$ is applied to the affine chart $\sper\
A_{ji}[\frac{y'}{x_{ji}}]$ and the element $x_{ji} \in \sper\
A_{ji}[\frac{y'}{x_{ji}}]$.

We claim that the graph $\Gamma_{i+1}$ contains a bamboo consisting of the
above $2\omega+1$ vertices, arranged in the following order :
\begin{multline}
 ([\delta'_1,\zeta_{\lambda+1}],s),\
([-\infty,\infty]_{\lambda+1},s),\ ([\zeta_{\lambda+1},\zeta_{\lambda+2}],-s),\
([-\infty,\infty]_{\lambda+2},s),\\
([\zeta_{\lambda+2},\zeta_{\lambda+3}],s),\ldots,
([-\infty,\infty]_{\lambda+\omega},s),\
([\zeta_{\lambda+\omega},\delta'_2],(-1)^\omega s).
\end{multline}

\noi We discuss a sample of edge of this bamboo, for example,
$([\zeta_t,\zeta_{t+1}],(-1)^t s),\ ([-\infty,\infty]_{t+1},s)$. The existence
of the other edges can be proved in a similar way.

The two maximal intervals $([\zeta_t,\zeta_{t+1}]$ and $[-\infty,\infty]_{t+1}$ have a common endpoint, namely, $\zeta_{t+1}$. We must show that
$$W^{(i+1)} \cap \{ \delta' \in X_{i+1}\ |\ sgn(x'_{ji}(\delta')) = (-1)^t s,\
sgn(x_{ji}(\delta'))=s,\ \sur{\{\delta'\}}
\ni
\zeta_{t+1} \}
 \neq \emptyset
.$$
The image of this set under $\pi_i$ is
$$ W^{(i)} \cap \{ \delta \in X_i\ |\ sgn(x_{ji}(\delta)) =
 s,\ \sur{\{\delta\}}
\ni
\xi_{t+1}^{(q)},\ \delta \mbox{ tangent to
}\{x_{ji}=0\} \}$$
and the result follows.

This proves 2.1(a). The cases 2.1(b) and (c) are similar but
easier.
\smallskip

Case 2.2: (a) Let $a=(I,s)$. Then $x_{ji}=y$. Put $y'=\frac{x}{y}$; $(x_{ji},y')$ is a regular system of
parameters at $\xi$. Let $A_{11}=A[x_{ji},y']$. The point $\xi \in \sper\ A_{11}$ is the unique
point such that supp$(\xi$)$=(x_{ji},y')$ and which induces the given
order on $k$. Let $A_{12}=A_{11}[x'_{ji},y']$ where
$x'_{ji}=\frac{x_{ji}}{y'}$.

Let $I' \subset \{x'_{ji}=0\}$ be the 1-maximal interval given by
$-\infty \leq y' \leq +\infty$ and $\tilde{I}' \subset \{y'=0\}$
the 1-maximal interval given by $ 0\leq x_{ji}' \leq +\infty$.

Now the vertices of $\Gamma_2$ are $(\tilde{I}',+),\ (I',+),\ (\tilde{I}',-)$
with the edges clearly defined.

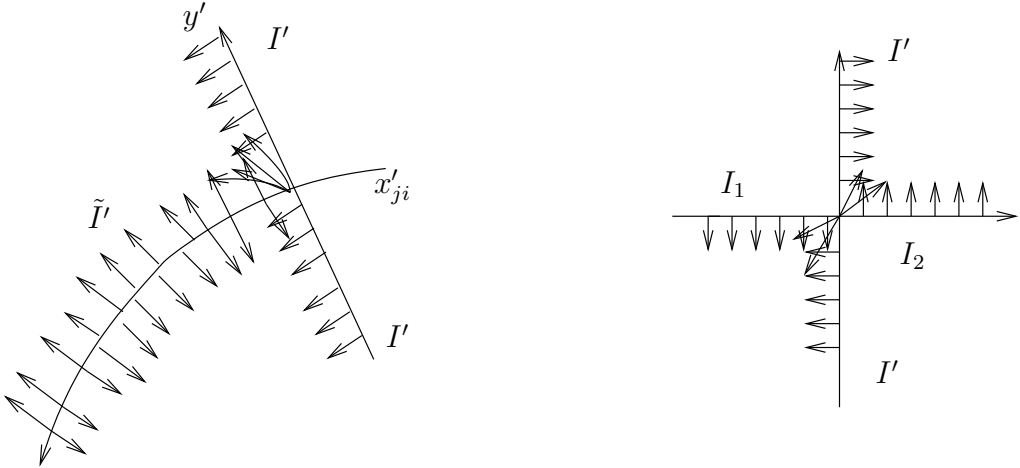
\begin{figure}[h!]
\begin{center}
\input{Usuite.pstex_t}
\caption{This figure shows the set $U^{(2)}$ in the cases 2.2.a and
  2.1.c respectively}
\label{Usuite}
\end{center}
\end{figure}

(b) Let $a=(I,s)$ be a vertex such that $\xi \in I$; the vertex $a$ is an endpoint of
$\Gamma_i$. Suppose that $\xi \in \sper\, A_{ji}$. Let $(x_{ji},y')$ be a
regular system of parameters at $\xi$. Let $A_{j,i+1}=A_{ji}[x'_{ji},y']$ where
$x'_{ji}=\frac{x_{ji}}{y'}$. Without loss of generality, assume that
$x_{ji} >0, y' >0$ on $W^{(i)}$.

Let $I' \subset \{x'_{ji}=0\}$ be the strict transform of $I$ in $\sper\
A_{j,i+1}$. Then $I'$ is an $(i+1)$-maximal interval. Let $\tilde{I}' \subset
\{y'=0\}$  be the $(i+1)$-maximal interval given by $ 0\leq x_{ji}' \leq
+\infty$. Now the new vertex $b$ added to $\Gamma_{i+1}$ is
$(\tilde{I}',+)$. It is connected by an edge to $a$ which is
represented in $\sper\ A_{j,i+}$ by $(I',+)$. This completes the proof of Lemma \ref{graphtransformation} and with it Proposition \ref{bamboo}.
$\Box$
\medskip

Let us finish the proof of
Theorem \ref{principale}. To each
vertex $(I=[\delta_1,\delta_2],s)$ of
$\Gamma_t$ we associate the
set $D(\delta_1,\delta_2) \subset U^{(t)}$ which by Corollary
\ref{loc-conn} is entirely contained in $F^{(t)}$ or $G^{(t)}$. This
defines a partition
$\Gamma_F=\left\{(I,s)\ \left|\  D(\delta_1,\delta_2) \subset
F^{(t)}\right.\right\}$, $\Gamma_G=\left\{(I,s)\ \left|\  D(\delta_1,\delta_2) \subset
G^{(t)} \right.\right\}$ of the set of vertices of $\Gamma_t$. Assume that $\Gamma_F \neq
\emptyset$ and $\Gamma_G \neq \emptyset$. Since $\Gamma_t$ is
connected, there exist $a=([\delta_{1a},\delta_{2a}],s_a) \in \Gamma_F$,
$b=([\delta_{1b},\delta_{2b}],s_b) \in \Gamma_G$ such that
$(a,b)$ is an edge of $\Gamma_t$. Then $D(\delta_{1a},\delta_{2a})
\subset F^{(t)}$, $D(\delta_{1b},\delta_{2b}) \subset G^{(t)}$ and
$D(\delta_{1a},\delta_{2a}) \cap D(\delta_{1b},\delta_{2b}) \neq \emptyset$.
This is a contradiction. This concludes the proof of
Theorem
\ref{principale}. $\Box$


\end{document}

%% file: UV2.pstex_t
\begin{picture}(0,0)%
\includegraphics{UV2.pstex}%
\end{picture}%
\setlength{\unitlength}{3947sp}%
\begingroup\makeatletter\ifx\SetFigFont\undefined%
\gdef\SetFigFont#1#2#3#4#5{%
  \reset@font\fontsize{#1}{#2pt}%
  \fontfamily{#3}\fontseries{#4}\fontshape{#5}%
  \selectfont}%
\fi\endgroup%
\begin{picture}(4436,1467)(3377,-3504)
\put(6608,-2221){\makebox(0,0)[lb]{\smash{{\SetFigFont{12}{14.4}{\rmdefault}{\mddefault}{\updefault}{\color[rgb]{0,0,0}$y$}%
}}}}
\put(3563,-2206){\makebox(0,0)[lb]{\smash{{\SetFigFont{12}{14.4}{\rmdefault}{\mddefault}{\updefault}{\color[rgb]{0,0,0}$y$}%
}}}}
\put(4382,-3053){\makebox(0,0)[lb]{\smash{{\SetFigFont{12}{14.4}{\rmdefault}{\mddefault}{\updefault}{\color[rgb]{0,0,0}$x$}%
}}}}
\put(7531,-3053){\makebox(0,0)[lb]{\smash{{\SetFigFont{12}{14.4}{\rmdefault}{\mddefault}{\updefault}{\color[rgb]{0,0,0}$x$}%
}}}}
\end{picture}%

%% file: edge.pstex_t
\begin{picture}(0,0)%
\includegraphics{edge.pstex}%
\end{picture}%
\setlength{\unitlength}{3947sp}%
\begingroup\makeatletter\ifx\SetFigFont\undefined%
\gdef\SetFigFont#1#2#3#4#5{%
  \reset@font\fontsize{#1}{#2pt}%
  \fontfamily{#3}\fontseries{#4}\fontshape{#5}%
  \selectfont}%
\fi\endgroup%
\begin{picture}(1899,1746)(4789,-6823)
\put(6601,-6436){\makebox(0,0)[lb]{\smash{{\SetFigFont{12}{14.4}{\rmdefault}{\mddefault}{\updefault}{\color[rgb]{0,0,0}$x_{jt}$}%
}}}}
\put(5401,-5236){\makebox(0,0)[lb]{\smash{{\SetFigFont{12}{14.4}{\rmdefault}{\mddefault}{\updefault}{\color[rgb]{0,0,0}$x_{\tilde{j}t}$}%
}}}}
\put(5476,-6436){\makebox(0,0)[lb]{\smash{{\SetFigFont{12}{14.4}{\rmdefault}{\mddefault}{\updefault}{\color[rgb]{0,0,0}$\xi$}%
}}}}
\end{picture}%

%% file: U1.pstex_t
\begin{picture}(0,0)%
\includegraphics{U1.pstex}%
\end{picture}%
\setlength{\unitlength}{3947sp}%
\begingroup\makeatletter\ifx\SetFigFont\undefined%
\gdef\SetFigFont#1#2#3#4#5{%
  \reset@font\fontsize{#1}{#2pt}%
  \fontfamily{#3}\fontseries{#4}\fontshape{#5}%
  \selectfont}%
\fi\endgroup%
\begin{picture}(6099,2869)(2014,-4859)
\put(8026,-3361){\makebox(0,0)[lb]{\smash{{\SetFigFont{12}{14.4}{\rmdefault}{\mddefault}{\updefault}{\color[rgb]{0,0,0}$x$}%
}}}}
\put(6976,-2311){\makebox(0,0)[lb]{\smash{{\SetFigFont{12}{14.4}{\rmdefault}{\mddefault}{\updefault}{\color[rgb]{0,0,0}$\frac{y}{x}$}%
}}}}
\put(3751,-3361){\makebox(0,0)[lb]{\smash{{\SetFigFont{12}{14.4}{\rmdefault}{\mddefault}{\updefault}{\color[rgb]{0,0,0}$y$}%
}}}}
\put(2776,-2161){\makebox(0,0)[lb]{\smash{{\SetFigFont{12}{14.4}{\rmdefault}{\mddefault}{\updefault}{\color[rgb]{0,0,0}$\frac{x}{y}$}%
}}}}
\put(2401,-4786){\makebox(0,0)[lb]{\smash{{\SetFigFont{12}{14.4}{\rmdefault}{\mddefault}{\updefault}{\color[rgb]{0,0,0}$\mbox{Sper }A[\frac{x}{y}]$}%
}}}}
\put(6601,-4786){\makebox(0,0)[lb]{\smash{{\SetFigFont{12}{14.4}{\rmdefault}{\mddefault}{\updefault}{\color[rgb]{0,0,0}$\mbox{Sper }A[\frac{y}{x}]$}%
}}}}
\end{picture}%

%% file: global4.pstex_t
\begin{picture}(0,0)%
\includegraphics{global4.pstex}%
\end{picture}%
\setlength{\unitlength}{1658sp}%
\begingroup\makeatletter\ifx\SetFigFont\undefined%
\gdef\SetFigFont#1#2#3#4#5{%
  \reset@font\fontsize{#1}{#2pt}%
  \fontfamily{#3}\fontseries{#4}\fontshape{#5}%
  \selectfont}%
\fi\endgroup%
\begin{picture}(17349,5049)(364,-4723)
\put(17026,-2536){\makebox(0,0)[lb]{\smash{{\SetFigFont{5}{6.0}{\rmdefault}{\mddefault}{\updefault}{\color[rgb]{0,0,0}$\delta'_2$}%
}}}}
\put(7876,-2536){\makebox(0,0)[lb]{\smash{{\SetFigFont{5}{6.0}{\rmdefault}{\mddefault}{\updefault}{\color[rgb]{0,0,0}$\delta'_1$}%
}}}}
\put(9901,-2536){\makebox(0,0)[lb]{\smash{{\SetFigFont{5}{6.0}{\rmdefault}{\mddefault}{\updefault}{\color[rgb]{0,0,0}$\zeta_{\lambda+1}$}%
}}}}
\put(12976,-2536){\makebox(0,0)[lb]{\smash{{\SetFigFont{5}{6.0}{\rmdefault}{\mddefault}{\updefault}{\color[rgb]{0,0,0}$\zeta_{\lambda+2}$}%
}}}}
\put(15076,-2086){\makebox(0,0)[lb]{\smash{{\SetFigFont{5}{6.0}{\rmdefault}{\mddefault}{\updefault}{\color[rgb]{0,0,0}$\zeta_{\lambda+3}$}%
}}}}
\put(7051,-2086){\makebox(0,0)[lb]{\smash{{\SetFigFont{5}{6.0}{\rmdefault}{\mddefault}{\updefault}{\color[rgb]{0,0,0}$\{x'_{ji}=0\}$}%
}}}}
\put(1951,-1711){\makebox(0,0)[lb]{\smash{{\SetFigFont{5}{6.0}{\rmdefault}{\mddefault}{\updefault}{\color[rgb]{0,0,0}$a$}%
}}}}
\put(451,-1186){\makebox(0,0)[lb]{\smash{{\SetFigFont{5}{6.0}{\rmdefault}{\mddefault}{\updefault}{\color[rgb]{0,0,0}$b$}%
}}}}
\put(3526,-1186){\makebox(0,0)[lb]{\smash{{\SetFigFont{5}{6.0}{\rmdefault}{\mddefault}{\updefault}{\color[rgb]{0,0,0}$c$}%
}}}}
\put(1051,-2611){\makebox(0,0)[lb]{\smash{{\SetFigFont{5}{6.0}{\rmdefault}{\mddefault}{\updefault}{\color[rgb]{0,0,0}$\xi^{(q)}_{\lambda+1}$}%
}}}}
\put(1651,-2911){\makebox(0,0)[lb]{\smash{{\SetFigFont{5}{6.0}{\rmdefault}{\mddefault}{\updefault}{\color[rgb]{0,0,0}$\xi^{(q)}_{\lambda+2}$}%
}}}}
\put(2476,-2611){\makebox(0,0)[lb]{\smash{{\SetFigFont{5}{6.0}{\rmdefault}{\mddefault}{\updefault}{\color[rgb]{0,0,0}$\xi^{(q)}_{\lambda+3}$}%
}}}}
\put(3526,-3136){\makebox(0,0)[lb]{\smash{{\SetFigFont{5}{6.0}{\rmdefault}{\mddefault}{\updefault}{\color[rgb]{0,0,0}$\delta_2$}%
}}}}
\put(526,-3136){\makebox(0,0)[lb]{\smash{{\SetFigFont{5}{6.0}{\rmdefault}{\mddefault}{\updefault}{\color[rgb]{0,0,0}$\delta_1$}%
}}}}
\put(3601,-2161){\makebox(0,0)[lb]{\smash{{\SetFigFont{5}{6.0}{\rmdefault}{\mddefault}{\updefault}{\color[rgb]{0,0,0}$\{x_{ji}=0\}$}%
}}}}
\put(5476,-2161){\makebox(0,0)[lb]{\smash{{\SetFigFont{5}{6.0}{\rmdefault}{\mddefault}{\updefault}{\color[rgb]{0,0,0}$\pi_i$}%
}}}}
\end{picture}%

%% file: Usuite.pstex_t
\begin{picture}(0,0)%
\includegraphics{Usuite.pstex}%
\end{picture}%
\setlength{\unitlength}{3947sp}%
\begingroup\makeatletter\ifx\SetFigFont\undefined%
\gdef\SetFigFont#1#2#3#4#5{%
  \reset@font\fontsize{#1}{#2pt}%
  \fontfamily{#3}\fontseries{#4}\fontshape{#5}%
  \selectfont}%
\fi\endgroup%
\begin{picture}(6399,2958)(2164,-4948)
\put(7801,-3661){\makebox(0,0)[lb]{\smash{{\SetFigFont{12}{14.4}{\rmdefault}{\mddefault}{\updefault}{\color[rgb]{0,0,0}$I_2$}%
}}}}
\put(7726,-2386){\makebox(0,0)[lb]{\smash{{\SetFigFont{12}{14.4}{\rmdefault}{\mddefault}{\updefault}{\color[rgb]{0,0,0}$I'$}%
}}}}
\put(7651,-4411){\makebox(0,0)[lb]{\smash{{\SetFigFont{12}{14.4}{\rmdefault}{\mddefault}{\updefault}{\color[rgb]{0,0,0}$I'$}%
}}}}
\put(3826,-2311){\makebox(0,0)[lb]{\smash{{\SetFigFont{12}{14.4}{\rmdefault}{\mddefault}{\updefault}{\color[rgb]{0,0,0}$I'$}%
}}}}
\put(4576,-4186){\makebox(0,0)[lb]{\smash{{\SetFigFont{12}{14.4}{\rmdefault}{\mddefault}{\updefault}{\color[rgb]{0,0,0}$I'$}%
}}}}
\put(2701,-3436){\makebox(0,0)[lb]{\smash{{\SetFigFont{12}{14.4}{\rmdefault}{\mddefault}{\updefault}{\color[rgb]{0,0,0}$\tilde{I'}$}%
}}}}
\put(3301,-2161){\makebox(0,0)[lb]{\smash{{\SetFigFont{12}{14.4}{\rmdefault}{\mddefault}{\updefault}{\color[rgb]{0,0,0}$y'$}%
}}}}
\put(6676,-3211){\makebox(0,0)[lb]{\smash{{\SetFigFont{12}{14.4}{\rmdefault}{\mddefault}{\updefault}{\color[rgb]{0,0,0}$I_1$}%
}}}}
\put(4501,-3211){\makebox(0,0)[lb]{\smash{{\SetFigFont{12}{14.4}{\rmdefault}{\mddefault}{\updefault}{\color[rgb]{0,0,0}$x'_{ji}$}%
}}}}
\end{picture}%

%% file: approximatelastversion.bbl
\begin{thebibliography}{99}
\bibitem{AM1} S. Abhyankar and T.T. Moh, {\em Newton--Puiseux expansion and
generalized Tschirnhausen transformation I.} Reine Agew. Math. 260, 47--83
(1973).
\bibitem{AM2} S. Abhyankar and T.T. Moh, {\em Newton--Puiseux expansion and
generalized Tschirnhausen transformation II.} Reine Agew. Math. 261, 29--54
(1973).
\bibitem{And} C. Andradas, L. Bröcker, J.M. Ruiz, {\em Constructible
    Sets in Real Geometry}, Springer (1996).
\bibitem{Al} D. Alvis, B. Johnston, J.J. Madden, {\em Local structure of
the real spectrum of a surface, infinitely near points and separating
ideals.} Preprint.
\bibitem{Baer} R. Baer, Uber nicht-archimedisch geordnete K\"{o}rper
  (Beitrage zur Algebra). Sitz. Ber. Der Heidelberger Akademie, 8
  Abhandl. (1927).
\bibitem{BiPi} G. Birkhoff and R. Pierce, {\em Lattice-ordered rings.}
Annales Acad. Brasil Ci{\^e}nc. 28, 41--69 (1956).
\bibitem{BCR} J. Bochnak, M. Coste, M.-F. Roy, {\em G\'eom\'etrie
alg\'ebrique r\'eelle.} Springer--Verlag, Berlin 1987.
\bibitem{CT} S.D. Cutkosky, B. Teissier, {\em Semi-groups of
    valuations on local rings}, Michigan Mat. J., Vol. 57, pp. 173-193
  (2008).
\bibitem{De} C. N. Delzell, {\em On the Pierce--Birkhoff conjecture over
ordered fields}, Rocky Mountain J. Math.
\bibitem{Fuc} L. Fuchs, {\em Telweise geordnete algebraische
    Strukturen.} Vandenhoeck and Ruprecht, 1966.
\bibitem{GAS} Herrera Govantes, F. J.; Olalla Acosta, M. A.; Spivakovsky, M.
{\em Valuations in algebraic field extensions}, Journal of Algebra ,
Vol. 312 , N. 2 , pp. 1033--1074 (2007).
\bibitem{GT} R. Goldin and B. Teissier, {\em Resolving singularities of
plane analytic branches with one toric morphism.}
\bibitem{GOST} F. J. Herrera Govantes, M. A. Olalla Acosta, M. Spivakovsky, B.
Teissier, {\em Extending a valuation centered in a local domain to
the formal completion.} Preprint.
\bibitem{HenIsb} M. Henriksen and J. Isbell, {\em Lattice-ordered rings
and function rings.} Pacific J. Math. 11, 533--566 (1962).
\bibitem{Kap} I. Kaplansky, {\em Maximal fields with valuations I.} Duke
Math. J., 9:303--321 (1942).
\bibitem{Kap2} I. Kaplansky, {\em Maximal fields with valuations II.} Duke
Math. J., 12:243--248 (1945).
\bibitem{Krull} W. Krull, Allgemeine Bewertungstheorie, J. Reine
  Angew. Math. 167, 160--196 (1932).
\bibitem{Kuo1} T.C. Kuo, {\em Generalized Newton--Puiseux theory and
Hensel's lemma in $\C[[x$, $y]]$.} Canadian J. Math., (6) XLI, 1101--1116
(1989).
\bibitem{Kuo2} T.C. Kuo, {\em A simple algorithm for deciding primes in
$\C[[x$, $y]]$.} Canadian J. Math., 47 (4), 801--816 (1995).
\bibitem{LJ} M. Lejeune-Jalabert, {\em Th{\'e}se d'Etat.} Universit{\'e}
Paris 7 (1973).
\bibitem{LMSS} F. Lucas, J.J. Madden, D. Schaub and M.
Spivakovsky, {\em On connectedness of sets in the real spectra of polynomial
rings}, Manuscripta Math. 128, 505-547, 2009.
\bibitem{LSS} F. Lucas, D. Schaub and M. Spivakovsky, {\em On the
    Pierce-Birkhoff conjecture in dimension 3}, in preparation.
\bibitem{Mac1} S. MacLane, {\em A construction for prime ideals as
absolute values of an algebraic field.} Duke Math. J. 2, 492--510 (1936).
\bibitem{Mac2} S. MacLane, {\em A construction for absolute values in
polynomial rings.} Transactions of the AMS 40, 363--395 (1936).
\bibitem{Mac3} S. MacLane and O.F.G. Schilling, {\em Zero-dimensional
branches of rank one on algebraic varieties.} Ann. of Math. 40, 3 (1939).
\bibitem{Mad1} J. J. Madden, {\em Pierce--Birkhoff rings}. Arch. Math. 53,
565--570 (1989).
\bibitem{Mad2} J. J. Madden, preprint.
\bibitem{Mah} L. Mah\'e, {\em On the Pierce--Birkhoff conjecture.} Rocky
Mountain J. Math. 14, 983--985 (1984).
\bibitem{Mar1} M. Marshall, {\em Orderings and real places of commutative
rings.} J. Alg. 140, 484--501 (1991).
\bibitem{Mar2} M. Marshall {\em The Pierce--Birkhoff conjecture for
curves.} Can. J. Math. 44, 1262--1271 (1992).
\bibitem{Mat} H. Matsumura {\em Commutative Algebra.} Benjamin/Cummings
Publishing Co., Reading, Mass., 1970.
\bibitem{Pre} A. Prestel {\em Lectures on formally real fields}, Lecture Notes in
Math., Springer--Verlag---Berlin,
Heidelberg, New York, 1984.
\bibitem{PD}  A. Prestel, C.N. Delzell {\em Positive Polynomials}, Springer monographs in
mathematics, 2001.
\bibitem{PC} S. Priess-Crampe, {\em Angeordnete strukturen: gruppen, k\"orper,
projektive Ebenen}, Springer -- Verlag --- Berlin,
Heidelberg, New York, 1983.
\bibitem{Sch} N. Schwartz, {\em Real closed spaces.} Habilitationsschrift,
M\"unchen 1984.
\bibitem{Spi1} M. Spivakovsky, {\em Valuations in function fields of
surfaces.} Amer. J. Math 112, 1, 107--156 (1990).
\bibitem{Spi2} M. Spivakovsky, {\em A solution to Hironaka's polyhedra
game.} Arithmetic and Geometry, Vol II, Papers dedicated to
I. R. Shafarevich on the occasion of his sixtieth birthday, M. Artin and
J. Tate, editors, Birkh{\"a}user, 1983, pp. 419--432.
\bibitem{Te} B. Teissier, {\em Valuations, deformations and toric
geometry.}, Proceedings of the Saskatoon Conference and Workshop on
valuation theory, Vol. II, F-V. Kuhlmann, S. Kuhlmann, M. Marshall,
editors, Fields Institute Communications, 33 (2003), 361-459.
\bibitem{V1} M. Vaqui\'e, {\em Famille admise associ\'ee \`a une
    valuation de $K[x]$}, S\'eminaires et Congr\`es 10, edited by
  Jean-Paul Brasselet and Tatsuo Suwa, 2005, pp. 391--428.
\bibitem{V2} M. Vaqui\'e, {\em Extension d'une valuation}, Trans. Amer. Math. Soc.  359  (2007),  no. 7, 3439--3481.
\bibitem{V3} M. Vaqui\'e, {\em Alg\`ebre gradu\'ee associ\'ee \`a une valuation de $K[x]$}, Adv. Stud. Pure Math., 46, Math. Soc. Japan, Tokyo, 2007.
\bibitem{V4} M. Vaqui\'e, {\em Famille admissible de valuations et
    défaut d'une extension}, J. Algebra  311  (2007),  no. 2,
  859--876.
\bibitem{V5} M. Vaqui\'e, {\em Valuations},  539--590, Progr. Math.,
  181, Birkhäuser, Basel, 2000.
\bibitem{W} S. Wagner, {\em On the Pierce-Birkhoff Conjecture for
    Smooth Affine Surfaces over Real Closed Fields}, ArXiv:0810.4800, 2009.
\bibitem{Zar} O. Zariski, P. Samuel {\em Commutative Algebra}, Vol. II, Springer Verlag.
\end{thebibliography}
